\documentclass [10pt] {article}

 \usepackage [cp1251] {inputenc}
 \usepackage [english] {babel}
 \usepackage {indentfirst}
 \usepackage {hyperref}
 \usepackage [curve] {xypic}
 \usepackage {tikz}
 \usepackage {mathtools}
 \usepackage {amsthm}
 \usepackage {amssymb}
 \usepackage {stmaryrd}                          
 \usepackage [mathscr] {eucal}
 \usepackage [LGR,OT1] {fontenc}                 

 \DeclareFontEncoding {LS1}{}{}                  
 \DeclareFontSubstitution {LS1}{stix2}{m}{n}
 \DeclareSymbolFont {stix2-symbols}{LS1}{stix2scr}{m}{n}         
 \DeclareSymbolFont {stix2-operators}{LS1}{stix2}{m}{n}          

 \DeclareMathAlphabet \sfi{OT1}{lmss}{m}{sl}
 \DeclareMathAlphabet \sfbi{OT1}{lmss}{bx}{sl}
 \DeclareMathAlphabet \gsfb{LGR}{lmss}{bx}{n}
 \DeclareMathAlphabet \gsfbi{LGR}{lmss}{bx}{sl}

 \let \til\~     

 \def \?{{?}}
 \def \atom{{\bullet}}
 \def \x{\times}                 

 \def \O{{\displaystyle\raisebox{0.99ex}{$\scriptscriptstyle\boldsymbol<$}\mkern-4.1mu{|}}}
 \def \o{{\scriptstyle\raisebox{0.51ex}{$\scriptscriptstyle<$}\mkern-4.4mu{|}}}
 \def \0{{\mathchoice\O\O\o\o}}

 \DeclareMathSymbol \sbighash{\mathord}{stix2-symbols}{"9F}      
 \def \shash{{\mathchoice
 {\vcenter{\hbox{$\scriptstyle\sbighash$}}}
 {\vcenter{\hbox{$\scriptstyle\sbighash$}}}
 {\vcenter{\hbox{$\scriptscriptstyle\sbighash$}}}
 {\vcenter{\hbox{$\scriptscriptstyle\sbighash$}}}
 }}
 \def \uc{\shash}                
 \def \huc{\widehat\shash{}}

 \def \ZZ{{\mathbb Z}}

 \def \Mg{{\boldsymbol{\mathrm{Mg}}}}            
 \def \sMg{{\boldsymbol{\mathrm{sMg}}}}
 \def \Ab{{\boldsymbol{\mathrm{Ab}}}}

 \def \F{{\mathscr F}}           
 \def \PP{{\mathscr P}}          
 \def \RR{{\mathscr R}}          
 \def \AA{{\mathscr A}}          
 \def \bd{\partial}              

 \def \id{{\mathrm{id}}}
 \def \pro{{\mathrm{pr}}}        
 \def \inc{{\mathrm{in}}}        
 \def \Inc{{\boldsymbol{\mathsf{in}}}}
 \def \ins{{\mathrm{in}}}        
 \def \Ins{{\boldsymbol{\mathsf{in}}}}

 \def \d{\triangle}
 \def \D{{\gsfb D}}
 \def \E{{\boldsymbol{\mathsf E}}}

 \def \b{{\gsfb b}}

 \def \oC{{\mathrm C}}           
 \def \vC{{\check{\mathrm C}}}
 \def \hC{{\hat{\mathrm C}}}
 \def \OC{{\boldsymbol{\mathsf C}}}
 \def \VC{{\boldsymbol{\check{\mathsf C}}}}
 \def \Vc{{\boldsymbol{\check{\mathsf c}}}}
 \def \HC{{\boldsymbol{\hat{\mathsf C}}}}

 \def \bS{{\overline\Sigma}}     
 \def \HS{{\boldsymbol{\hat{\gsfb S}}}}
 
 \def \K{{\sfi K}}               
 \def \L{{\sfi L}}
 \def \M{{\sfi M}}

 \def \A{{\sfbi A}}              
 \def \B{{\sfbi B}}
 \def \T{{\sfbi T}}
 \def \tT{\skew{3}\widetilde\T}
 \def \U{{\sfbi U}}
 \def \V{{\sfbi V}}
 \def \W{{\sfbi W}}
 \def \Z{{\sfbi Z}}
 \def \tZ{\skew{3}\widetilde\Z}
 \def \e{{\sfbi e}}
 \def \f{{\sfbi f}}
 \def \g{{\sfbi g}}
 \def \h{{\sfbi h}}
 \def \i{{\sfbi i}}
 \def \k{{\sfbi k}}
 \def \p{{\sfbi p}}
 \def \q{{\sfbi q}}
 \def \r{{\sfbi r}}
 \def \v{{\sfbi v}}
 \def \dg{{\gsfbi d}}    
 \def \xg{{\gsfbi x}}    
 \def \rh{{\gsfbi r}}    
 \def \sg{{\gsfbi s}}    
 \def \tg{{\gsfbi t}}    

 \def \<{\langle}                
 \def \>{\rangle}

 \def \LQ{$\scriptscriptstyle<$}                 
 \def \`{{\mathchoice
 {\raisebox{0.20ex}\LQ}
 {\raisebox{0.20ex}\LQ}
 {\raisebox{0.10ex}\LQ}
 {\raisebox{0.05ex}\LQ}
 }}

 \def \RQ{$\scriptscriptstyle>$}                 
 \def \'{{\mathchoice
 {\raisebox{0.20ex}\RQ}
 {\raisebox{0.20ex}\RQ}
 {\raisebox{0.10ex}\RQ}
 {\raisebox{0.05ex}\RQ}
 }}

 \def \({\llbracket}             
 \def \){\rrbracket}
 \def \[{\lfloor}                
 \def \]{\rceil}

 \def \bigbou{\bigvee}                           

 \newcommand* {\bds} [3] {                       
 \vcenter {
 \hrule height #1
 \vspace {#2}
 \hbox{$#3$}
 \vspace {#2}
 \hrule height #1
 }
 }

 \def \sbigOplus{{\mathchoice                    
 {\bds{0.20ex}{0.23ex}{\displaystyle\bigoplus}}
 {\bds{0.15ex}{0.25ex}{\textstyle\bigoplus}}
 {\bds{0.13ex}{0.21ex}{\scriptstyle\bigoplus}}
 {\bds{0.09ex}{0.15ex}{\scriptscriptstyle\bigoplus}}
 }}
 \def \bigOplus{\mathop{\sbigOplus}}

 \def \sbigCop{{\mathchoice                      
 {\bds{0.15ex}{0.25ex}{\textstyle\coprod}}
 {\bds{0.15ex}{0.25ex}{\textstyle\coprod}}
 {\bds{0.13ex}{0.21ex}{\scriptstyle\coprod}}
 {\bds{0.09ex}{0.15ex}{\scriptscriptstyle\coprod}}
 }}
 \def \bigCop{\mathop{\sbigCop}}

 \def \sbigBou{{\mathchoice                      
 {\bds{0.15ex}{0.25ex}{\textstyle\bigvee}}
 {\bds{0.15ex}{0.25ex}{\textstyle\bigvee}}
 {\bds{0.13ex}{0.21ex}{\scriptstyle\bigvee}}
 {\bds{0.09ex}{0.15ex}{\scriptscriptstyle\bigvee}}
 }}
 \def \bigBou{\mathop{\sbigBou}}

 \newcommand* {\bxd} [5] {                       
 \vcenter {
 \vspace{#4}
 \hbox{\hspace{#3}\tikz{
 \node [inner sep=0ex,minimum size=#1,line width=#2,draw] {$#5$}
 }\hspace{#3}}
 \vspace{#4}
 }
 }

 \def \sbigecol{{\mathchoice                     
 {\bxd{3.1ex}{0.15ex}{0.15ex}{0.23ex}{}}
 {\bxd{2.2ex}{0.12ex}{0.13ex}{0ex}{}}
 {\bxd{1.5ex}{0.10ex}{0.13ex}{0ex}{}}
 {\bxd{1.1ex}{0.07ex}{0.09ex}{0ex}{}}
 }}
 \def \bigecol{\mathop\sbigecol}

 \def \sbigeCop{{\mathchoice                     
 {\bxd{3.1ex}{0.15ex}{0.15ex}{0.23ex}{\displaystyle\boldsymbol\sqcup}}
 {\bxd{2.2ex}{0.12ex}{0.13ex}{0ex}{\scriptstyle\boldsymbol\sqcup}}
 {\bxd{1.5ex}{0.10ex}{0.13ex}{0ex}{\scriptscriptstyle\boldsymbol\sqcup}}
 {\bxd{1.1ex}{0.07ex}{0.09ex}{0ex}{\scriptscriptstyle\sqcup}}
 }}
 \def \bigeCop{\mathop\sbigeCop}

 \def \sbigeBou{{                                
 \mathchoice
 {\bxd{3.1ex}{0.15ex}{0.15ex}{0.23ex}{\displaystyle\boldsymbol\vee}}
 {\bxd{2.2ex}{0.12ex}{0.13ex}{0ex}{\scriptstyle\boldsymbol\vee}}
 {\bxd{1.5ex}{0.10ex}{0.13ex}{0ex}{\scriptscriptstyle\boldsymbol\vee}}
 {\bxd{1.1ex}{0.07ex}{0.09ex}{0ex}{\scriptscriptstyle\vee}}
 }}
 \def \bigeBou{\mathop\sbigeBou}

 \let \Im\undefined
 \DeclareMathOperator \Im{Im}
 \DeclareMathOperator \Fix{Fix}

 \def \bou{\vee}         
 \def \Bou{\mathbin{\overline{\underline\vee}}}                  

 \def \secol{{\mathchoice                                        
 {\bxd{1.47ex}{0.10ex}{0.13ex}{0ex}{}}
 {\bxd{1.47ex}{0.10ex}{0.13ex}{0ex}{}}
 {\bxd{1.13ex}{0.08ex}{0.13ex}{0ex}{}}
 {\bxd{0.88ex}{0.07ex}{0.14ex}{0ex}{}}
 }}

 \def \seBou{                                                    
 {\bxd{1.47ex}{0.10ex}{0.13ex}{0ex}{\scriptscriptstyle\boldsymbol\vee}}
 }
 \def \eBou{\mathbin{\seBou}}

 \def \cro{\times}       
 \def \sma{\wedge}       

 \def \ssm{{\mathchoice         
 {{\scriptstyle\wedge}}
 {{\scriptstyle\wedge}}
 {{\scriptscriptstyle\wedge}}
 {{\scriptscriptstyle\wedge}}
 }}
 \def \sm{\mathbin\ssm}

 \def \hsm{\leftthreetimes}      

 \def \shs{{\mathchoice         
 {{\scriptstyle\leftthreetimes}}
 {{\scriptstyle\leftthreetimes}}
 {{\scriptscriptstyle\leftthreetimes}}
 {{\scriptscriptstyle\leftthreetimes}}
 }}
 \def \hs{\mathbin\shs}

 \def \cat{\mathbin*}    
 \def \mee{\wedge}       

 \def \sOplus{{\mathchoice                                       
 {\bds{0.06ex}{0.29ex}{\displaystyle\oplus}}
 {\bds{0.06ex}{0.29ex}{\textstyle\oplus}}
 {\bds{0.05ex}{0.25ex}{\scriptstyle\oplus}}
 {\bds{0.05ex}{0.17ex}{\scriptscriptstyle\oplus}}
 }}
 \def \Oplus{\mathbin{\sOplus}}

 \def \le{\leqslant}
 \def \ge{\geqslant}

 \DeclareMathSymbol \approxident{\mathrel}{stix2-operators}{"FC}                 

 \def \xto{\xrightarrow}
 \def \xfrom{\xleftarrow}

 \def \-{\overline}              
 \def \_{\underline}
 \def \~{\widetilde}

 \def \+{\boldsymbol}            

 \renewcommand* {\=} [2] {       
 \prescript{#1}{}#2
 }

 \renewcommand* {\%} [2] {       
 \overset{#2}#1
 }

 \renewcommand* {\|} [2] {       
 \mathrel{{#1}|_{#2}}
 }

 \newcommand* {\head} [1] {
 \subsubsection* {#1}
 }

 \newcommand* {\subhead} [1] {
 \addvspace\medskipamount
 \noindent {\it #1\/}
 }

 \newenvironment* {claim} [1] []
 {\begin{trivlist}\item [\hskip\labelsep {\bf #1}] \it}
 {\end{trivlist} }

 \newenvironment* {demo} [1] []
 {\begin{trivlist}\item [\hskip\labelsep {\it #1}] }
 {\end{trivlist} }

 \newenvironment* {dem}
 {\begin{trivlist}\item}
 {\end{trivlist} }

 \makeatletter                   
 \def \QED{\displaymath@qed}
 \makeatother

\begin {document}

 \title {\large\bf
         Homotopy similarity of maps.
         Strong similarity}

 \author {\normalsize\rm
          S.~S.~Podkorytov}

 \date {}

 \maketitle

 \vspace {-2\bigskipamount}
 
 \begin {abstract} \noindent
 Given based cellular spaces $X$ and $Y$,
 $X$ compact,
 and
 an integer $r\ge0$,
 we
 define a relation $\%\approx r$ on $[X,Y]$
 and
 argue for the conjecture that
 it always coincides with the $r$-similarity $\%\sim r$.
 \end {abstract}


 \head {\S~1. Introduction}


 This paper continues \cite{sim}.
 We adopt notation and conventions thereof.
 Let $X$ and $Y$ be cellular spaces,
 $X$ compact.
 For each $r\ge0$,
 we define a relation $\%\approx r$,
 called the strong $r$-similarity,
 on the set $[X,Y]$.
 We will need it in our next paper \cite{sim-3-com}.
 We conjecture that
 strong $r$-similarity always coincides with
 $r$-similarity $\%\sim r$.
 It follows immediately from the definition that
 it implies $r$-similarity and
 gets nonstrictly stronger as $r$ grows.
 We prove that
 the strong $r$-similarity is an equivalence
 provided $X$ is a suspension
 (\S~8).
 The main results are as follows.
 Strong $1$-similarity coincides with $1$-similarity
 (Theorem~14.2).
 (We believe that
 $1$-similarity can be given a homological characterization
 similar to that of homotopy invariants of order at most $1$
 \cite{str}.)
 If $X=S^1$,
 the strong $r$-similarity coincides with the $r$-similarity
 (\S~24).
 All $(r+1)$-fold Whitehead products are strongly
 $r$-similar to zero
 (Theorem~27.2).


 \head {\S~2. Definition of strong similarity}


 \subhead {Augmentation.}
 For a set $W$,
 introduce the homomorphism
 $$
 \epsilon:
 \<W\>
 \to
 \ZZ,
 \qquad
 \`w\'
 \mapsto 1,
 $$
 the {\it augmentation}.
 An ensemble $S\in\<W\>$ is called {\it affine\/} if
 $\epsilon(S)=1$.

 \subhead {Unbased maps.}
 Given unbased spaces
 $U$
 and
 $V$,
 we let
 $V^{(U)}$ be the unbased space of unbased maps $U\to V$.
 Introduce the unbased map
 $$
 \Xi^U:V\to V^{(U)},
 \qquad
 \Xi^U(v):u\mapsto v.
 $$


 \subhead {Combining product of ensembles.}
 Given
 a coproduct of unbased spaces
 $$
 U
 =
 \coprod_{i\in(m)}
 U_i
 $$
 (hereafter,
 $(m)=\{1,\dotsc,m\}$)
 and
 an unbased space $V$,
 we have
 the operation of combining
 $$
 \prod_{i\in(m)}
 V^{(U_i)}
 \to
 V^{(U)},
 \qquad
 (w_i)_{i\in(m)}
 \mapsto
 \bigCop_{i\in(m)}
 w_i,
 $$
 and
 the $\ZZ$-multilinear operation
 \begin {equation} \label {bigeCop}
 \bigeCop_{i\in(m)}:
 \prod_{i\in(m)}
 \<V^{(U_i)}\>
 \to
 \<V^{(U)}\>,
 \qquad
 \bigeCop_{i\in(m)}
 \`w_i\'
 =
 \`
 \bigCop_{i\in(m)}
 w_i
 \',
 \end {equation}
 which we call the {\it combining product}.


 \subhead {Simplex and its faces.}
 Fix a nonempty finite set $E$.
 Let $\PP_\x(E)$ be the set of nonempty subsets $F\subseteq E$.
 Let $\AA(E)$ be the set of subsets
 $A\subseteq\PP_\x(E)$ such that
 all $F\in A$ are disjoint
 ({\it layouts\/}).

 Let $\Delta E$ be the simplex spanned by $E$.
 For $F\in\PP_\x(E)$,
 $\Delta F\subseteq\Delta E$ is a face.
 For $A\in\AA(E)$,
 put
 $$
 \Delta[A]
 =
 \coprod_{F\in A}
 \Delta F
 \subseteq
 \Delta E.
 $$


 \subhead {Fissile ensembles.}
 Given an unbased space $V$,
 we call an ensemble $S\in\<V^{(\Delta E)}\>$ {\it fissile\/}
 if,
 for each $A\in\AA(E)$,
 \begin {equation} \label {fis}
 S|_{\Delta[A]}
 =
 \bigeCop_{F\in A}
 S|_{\Delta F}
 \end {equation}
 in $\<V^{(\Delta[A])}\>$.

 An ensemble of the form $\`w\'$ is fissile.
 A fissile ensemble is affine
 (take $A=\varnothing$ in the definition).
 An affine ensemble $S$ is fissile if
 it satisfies \eqref{fis}
 for all $A$
 with $|A|=2$.
 Given an unbased space $\~V\supseteq V$,
 we have $\<\~V^{(\Delta E)}\>\supseteq\<V^{(\Delta E)}\>$;
 the ensemble $S$ is fissile
 as an element of $\<\~V^{(\Delta E)}\>$
 if and only if
 it is fissile
 as an element of $\<V^{(\Delta E)}\>$.
 

 \subhead {Spaces of maps.}
 Let
 $X$
 and
 $Y$
 be spaces,
 $X$ compact Hausdorff.
 Then
 $Y^X$ is the space of maps $X\to Y$;
 its basepoint is the constant map $\0^X_Y$.
 Given a map $a:X\to Y$,
 let $Y^X_a\subseteq Y^X$ be the path component containing $a$.


 \subhead {The filtration $\<(Y^X)^{(U)}\>^{(s)}_X$.}
 Let $U$ be an unbased space.
 Introduce the space
 $$
 U\hsm X
 =
 (U\cro X)/(U\cro\{\0_X\}).
 $$
 We have the projection
 $$
 U\cro X
 \to
 U\hsm X,
 \qquad
 (u,x)
 \mapsto
 u\hs x.
 $$
 The bijection
 $$
 \uc^X:
 (Y^X)^{(U)}
 \to
 Y^{U\hsm X},
 \qquad
 \uc^X(w):u\hs x\mapsto w(u)(x),
 $$
 induces the isomorphism
 $$
 \<\uc^X\>:
 \<(Y^X)^{(U)}\>
 \to
 \<Y^{U\hsm X}\>.
 $$
 The filtration of $\<Y^{U\hsm X}\>$
 (see \cite{sim})
 induces a filtration of $\<(Y^X)^{(U)}\>$:
 $$
 \<(Y^X)^{(U)}\>^{(s)}_X
 =
 \<\uc^X\>^{-1}
 (
 \<Y^{U\hsm X}\>^{(s)}
 ).
 $$


 \subhead {Strong similarity.}
 Let
 $X$
 and
 $Y$
 be cellular spaces,
 $X$ compact,
 and
 $a,b:X\to Y$ be maps.
 We say that $a$ is {\it strongly $r$-similar\/} to $b$,
 $$
 a\%\approx rb,
 $$
 if,
 for any nonempty finite set $E$,
 there exists a fissile ensemble
 $S\in\<(Y^X_a)^{(\Delta E)}\>\subseteq\<(Y^X)^{(\Delta E)}\>$
 such that
 $$
 \`\Xi^{\Delta E}(b)\'
 -
 S
 \in
 \<(Y^X)^{(\Delta E)}\>^{(r+1)}_X.
 $$
 We have $a\%\approx r a$
 (put $S=\`\Xi^{\Delta E}(a)\'$).
 Clearly,
 $a\%\approx rb$ implies $a\%\sim rb$
 (take $E=\{\atom\}$).
 We prove below
 (Theorem~6.1)
 that the relation $\%\approx r$ is homotopy invariant.
 

 \head {\S~3. Filtrations $\<(Y^X)^{(U)}\>^{(s)}_X$ and
 $\<(Y^X)^T\>^{(s)}_X$.}


 \subhead {Naturality of the filtration
 $\<(Y^X)^{(U)}\>^{(s)}_X$.}

 \begin {claim} [3.1. Lemma.]
 Let
 $X$
 and
 $Y$
 be spaces,
 $X$
 compact Hausdorff,
 $U$
 and
 $\~U$
 be unbased spaces,
 and
 $k:\~U\to U$ be an unbased map.
 Then the homomorphism
 $$
 \<(Y^X)^{(k)}\>:
 \<(Y^X)^{(U)}\>
 \to
 \<(Y^X)^{(\~U)}\>
 $$
 takes
 $\<(Y^X)^{(U)}\>^{(s)}_X$
 to
 $\<(Y^X)^{(\~U)}\>^{(s)}_X$.
 \end {claim}

 \begin {demo} [Proof.]
 We have the commutative diagram
 $$
 \xymatrix {
 \<(Y^X)^{(U)}\>
 \ar[rr]^-{\<\uc^X\>}
 \ar[d]_-{\<(Y^X)^{(k)}\>}
 &&
 \<Y^{U\hsm X}\>
 \ar[d]^-{\<Y^{k\hsm\id_X}\>}
 \\
 \<(Y^X)^{(\~U)}\>
 \ar[rr]^-{\<\uc^X\>}
 &&
 \<Y^{\~U\hsm X}\>.
 }
 $$
 By the definition of $\<(Y^X)^{(U)}\>^{(s)}_X$,
 $\<\uc^X\>$ takes
 it
 to
 $\<Y^{U\hsm X}\>^{(s)}$.
 By \cite[Lemma~2.1]{sim-1-cir},
 $\<Y^{k\hsm\id_X}\>$
 takes
 the latter
 to
 $\<Y^{\~U\hsm X}\>^{(s)}$.
 By commutativity of the diagram,
 $\<(Y^X)^{(k)}\>$
 takes
 $\<(Y^X)^{(U)}\>^{(s)}_X$
 to
 $\<\uc^X\>^{-1}(\<Y^{\~U\hsm X}\>^{(s)})$,
 which is $\<(Y^X)^{(\~U)}\>^{(s)}_X$
 by the definition of the latter.
 \qed
 \end {demo}


 \subhead {A technical lemma.}

 \begin {claim} [3.2. Lemma.]
 Let
 $X$,
 $Y$,
 and $\~X$
 be spaces,
 $X$
 and
 $\~X$
 compact Hausdorff,
 and
 $k:\~X\to X$ be a surjective map.
 Then the homomorphism
 $$
 \<Y^k\>:
 \<Y^X\>
 \to
 \<Y^{\~X}\>
 $$
 satisfies
 \begin {equation} \label {preimage}
 \<Y^X\>^{(s)}
 =
 \<Y^k\>^{-1}(\<Y^{\~X}\>^{(s)}).
 \end {equation}
 \end {claim}

 \begin {demo} [Proof.]
 By \cite[Lemma~2.1]{sim-1-cir},
 $\<Y^k\>$ preserves the filtration,
 which yields the inclusion $\subseteq$ in \eqref{preimage}.
 Check the inclusion $\supseteq$.
 Take $V\in\<Y^k\>^{-1}(\<Y^{\~X}\>^{(s)})$
 and
 show that
 $V\in\<Y^X\>^{(s)}$.
 Take $R\in\F_{s-1}(X)$.
 We should check that
 $V|_R=0$.
 We have $R=k(Q)$
 for some $Q\in\F_{s-1}(\~X)$.
 Since $\<Y^k\>(V)\in\<Y^{\~X}\>^{(s)}$,
 we have $\<Y^k\>(V)|_Q=0$.
 We have the commutative diagram
 $$
 \xymatrix @C1.4ex {
 {\scriptstyle
 V
 }
 &
 \<Y^X\>
 \ar[rr]^-{\<Y^k\>}
 \ar[d]_-{\?|_R}
 &&
 \<Y^{\~X}\>
 \ar[d]^-{\?|_Q}
 &
 {\scriptstyle
 \<Y^k\>(V)
 }
 \\
 {\scriptstyle
 V|_R
 }
 &
 \<Y^R\>
 \ar[rr]^-{\<Y^h\>}
 &&
 \<Y^Q\>,
 &
 {\scriptstyle
 0
 }
 }
 $$
 where $h=k|_{Q\to R}$.
 Since $h$ is surjective,
 $\<Y^h\>$ is injective.
 Thus
 $V|_R=0$.
 \qed
 \end {demo}


 \subhead {The filtration $\<(Y^X)^T\>^{(s)}_X$.}
 Let
 $X$,
 $Y$,
 and
 $T$
 be spaces,
 $X$
 compact Hausdorff.
 We
 have the inclusion $\<(Y^X)^T\>\subseteq\<(Y^X)^{(T)}\>$
 and
 define the subgroups $\<(Y^X)^T\>^{(s)}_X\subseteq\<(Y^X)^T\>$
 by putting
 \begin {equation} \label {YXT}
 \<(Y^X)^T\>^{(s)}_X
 =
 \<(Y^X)^T\>
 \cap
 \<(Y^X)^{(T)}\>^{(s)}_X
 \subseteq
 \<(Y^X)^{(T)}\>.
 \end {equation}
 We have
 the projection
 $$
 T\cro X
 \to
 T\sma X,
 \qquad
 (t,x)
 \mapsto
 t\sm x,
 $$
 and
 the bijection
 $$
 \huc^X:
 (Y^X)^{T}
 \to
 Y^{T\sma X},
 \qquad
 \huc^X(v):
 t\sm x
 \mapsto
 v(t)(x).
 $$

 \begin {claim} [3.3. Lemma.]
 One has
 $$
 \<(Y^X)^T\>^{(s)}_X
 =
 \<\huc^X\>^{-1}(\<Y^{T\sma X}\>^{(s)}).
 $$
 \end {claim}

 \begin {demo} [Proof.]
 We have
 the projection
 $$
 k:
 T\hsm X
 \to
 T\sma X,
 \qquad
 t\hs x
 \mapsto
 t\sm x,
 $$
 and
 the commutative diagram
 $$
 \xymatrix {
 \<(Y^X)^T\>
 \ar[r]^-{\<\huc^X\>}
 \ar[d]_-{\inc}
 &
 \<Y^{T\sma X}\>
 \ar[d]^-{\<Y^k\>}
 \\
 \<(Y^X)^{(T)}\>
 \ar[r]^-{\<\uc^X\>}
 &
 \<Y^{T\hsm X}\>.
 }
 $$
 By the definitions,
 $$
 \<(Y^X)^T\>^{(s)}_X
 =
 \inc^{-1}(\<(Y^X)^{(T)}\>^{(s)}_X)
 $$
 and
 $$
 \<(Y^X)^{(T)}\>^{(s)}_X
 =
 \<\uc^X\>^{-1}(\<Y^{T\hsm X}\>^{(s)}).
 $$
 By Lemma~3.2,
 $$
 \<Y^{T\sma X}\>^{(s)}
 =
 \<Y^k\>^{-1}(\<Y^{T\hsm X}\>^{(s)}).
 $$
 The desired equality follows
 by the diagram.
 \qed
 \end {demo}


 \head {\S~4. Primitive transforms}


 Let
 $V$
 and
 $\~V$
 be unbased spaces
 and
 $g:V\to\~V$ be an unbased map.
 For an unbased space $U$,
 we have the induced function $g^{(U)}:V^{(U)}\to\~V^{(U)}$.

 \begin {claim} [4.1. Lemma.]
 Let $E$ be a nonempty finite set.
 Consider the homomorphism
 $$
 \<g^{(\Delta E)}\>:
 \<V^{(\Delta E)}\>
 \to
 \<\~V^{(\Delta E)}\>.
 $$
 Then,
 for any fissile ensemble $S\in\<V^{(\Delta E)}\>$,
 the ensemble $\<g^{(\Delta E)}\>(S)$ is fissile.
 \end {claim}

 \begin {demo} [Proof.]
 Take $A\in\AA(E)$.
 We have the commutative diagram
 $$
 \xymatrix @C1.4ex @R1.4ex {
 {\scriptstyle
 (S)_{F\in A}
 }
 \ar@{|->}[rrrrrr]^-{(1)}
 \ar@{|->}[ddd]_-{(2)}
 &&&&&&
 {\scriptstyle
 (\<g^{(\Delta E)}\>(S))_{F\in A}
 }
 \ar@{|->}[ddd]^-{(3)}
 \\
 &
 \prod\limits_{F\in A}
 \<V^{(\Delta E)}\>
 \ar[rrrr]^{
 \prod\limits_{F\in A}
 \<g^{(\Delta E)}\>
 }
 \ar[dd]_-{
 \prod\limits_{F\in A}
 \?|_{\Delta F}
 }
 &&&&
 \prod\limits_{F\in A}
 \<\~V^{(\Delta E)}\>
 \ar[dd]^-{
 \prod\limits_{F\in A}
 \?|_{\Delta F}
 }
 &
 \\ \\
 {\scriptstyle
 (S|_{\Delta F})_{F\in A}
 }
 \ar@{|->}[dd]_-{(4)}
 &
 \prod\limits_{F\in A}
 \<V^{(\Delta F)}\>
 \ar[rrrr]^{
 \prod\limits_{F\in A}
 \<g^{(\Delta F)}\>
 }
 \ar[dd]_-{
 \bigeCop\limits_{F\in A}
 }
 &&&&
 \prod\limits_{F\in A}
 \<\~V^{(\Delta F)}\>
 \ar[dd]^-{
 \bigeCop\limits_{F\in A}
 }
 &
 {\scriptstyle
 (\<g^{(\Delta E)}\>(S)|_{\Delta F})_{F\in A}
 }
 \ar@{|->}[dd]^-{(5)}
 \\ \\
 {\scriptstyle
 S|_{\Delta[A]}
 }
 &
 \<V^{(\Delta[A])}\>
 \ar[rrrr]^-{\<g^{(\Delta[A])}\>}
 &&&&
 \<\~V^{(\Delta[A])}\>
 &
 {\scriptstyle
 \bigeCop\limits_{F\in A}
 \<g^{(\Delta E)}\>(S)|_{\Delta F}
 }
 \\ \\
 &
 \<V^{(\Delta E)}\>
 \ar[uu]^-{\?|_{\Delta[A]}}
 \ar[rrrr]^-{\<g^{(\Delta E)}\>}
 &&&&
 \<\~V^{(\Delta E)}\>.
 \ar[uu]_-{\?|_{\Delta[A]}}
 &
 \\
 {\scriptstyle
 S
 }
 \ar@{|->}[uuu]^-{(6)}
 \ar@{|->}[rrrrrr]^-{(8)}
 &&&&&&
 {\scriptstyle
 \<g^{(\Delta E)}\>(S)
 }
 \ar@{|->}[uuu]_-{(7)}
 }
 $$
 The sending (4) is fissility of $S$.
 The sendings (1), (2), (3), (5), (6), and (8) are obvious.
 The sending (7) follows.
 It is fissility of $\<g^{(\Delta E)}\>(S)$.
 \qed
 \end {demo}


 \subhead {Primitivity.}
 Let
 $X$,
 $Y$,
 $\~X$,
 and
 $\~Y$
 be spaces,
 $X$
 and
 $\~X$
 compact Hausdorff,
 and
 $g:Y^X\to\~Y^{\~X}$
 be an unbased map
 (a {\it transform}).
 We suppose that
 the transform $g$ is
 {\it primitive\/}:
 for each point $p\in\~X$,
 there is
 a point $k(p)\in X$
 and
 an unbased map $h^p:Y\to\~Y$
 such that
 $$
 g(d)(p)
 =
 h^p(d(k(p))),
 \qquad
 d\in Y^X.
 $$

 \begin {claim} [4.2. Lemma.]
 For
 an unbased space $U$,
 the homomorphism
 $$
 \<g^{(U)}\>:
 \<(Y^X)^{(U)}\>
 \to
 \<(\~Y^{\~X})^{(U)}\>
 $$
 takes $\<(Y^X)^{(U)}\>^{(s)}_X$
 to $\<(\~Y^{\~X})^{(U)}\>^{(s)}_{\~X}$.
 \end {claim}

 \begin {demo} [Proof.]
 We may assume that
 $k(\0_{\~X})=\0_X$
 and
 $h^{\0_{\~X}}(\0_Y)=\0_{\~Y}$.
 We have the
 (possibly discontinuous)
 function
 $$
 K
 =
 \id\hsm k:
 U\hsm\~X
 \to
 U\hsm X.
 $$
 For $Q\in\F_{s-1}(U\hsm\~X)$,
 we have $K(Q)\in\F_{s-1}(U\hsm X)$.
 We have
 the function
 $$
 H:
 Y^{K(Q)}
 \to
 \~Y^Q,
 \qquad
 H(v):
 u\hs p\mapsto h^p(v(K(u\hs p))),
 \
 u\hs p\in Q,
 $$
 and
 the commutative diagram
 $$
 \xymatrix {
 \<(Y^X)^{(U)}\>
 \ar[r]^-{\<\uc^X\>}
 \ar[d]_-{\<g^{(U)}\>}
 &
 \<Y^{U\hsm X}\>
 \ar[r]^-{\?|_{K(Q)}}
 &
 \<Y^{K(Q)}\>
 \ar[d]^-{\<H\>}
 \\
 \<(\~Y^{\~X})^{(U)}\>
 \ar[r]^-{\<\uc^{\~X}\>}
 &
 \<\~Y^{U\hsm\~X}\>
 \ar[r]^-{\?|_Q}
 &
 \<\~Y^Q\>.
 }
 $$
 By the definition of $\<(Y^X)^{(U)}\>^{(s)}_X$,
 it goes to zero
 under the composition in the upper row.
 Thus
 its image under $\<g^{(U)}\>$
 goes to zero
 under the composition in the lower row.
 Since $Q$ was taken arbitrarily,
 this image is contained in $\<(\~Y^{\~X})^{(U)}\>^{(s)}_{\~X}$
 by the definition of the latter.
 \qed
 \end {demo}

 Suppose that
 $X$,
 $Y$,
 $\~X$,
 and
 $\~Y$
 are cellular.

 \begin {claim} [4.3. Lemma.]
 Let $a,b:X\to Y$ be maps such that
 $a\%\approx r b$.
 Then
 $g(a)\%\approx r g(b)$.
 \end {claim}

 \begin {demo} [Proof.]
 Take a nonempty finite set $E$.
 We have a fissile ensemble $S\in\<(Y^X_a)^{(\Delta E)}\>$ such that
 $$
 \`\Xi^{\Delta E}(b)\'
 -
 S
 \in
 \<(Y^X)^{(\Delta E)}\>^{(r+1)}_X.
 $$
 Consider the homomorphism
 $$
 \<g^{(\Delta E)}\>:
 \<(Y^X)^{(\Delta E)}\>
 \to
 \<(\~Y^{\~X})^{(\Delta E)}\>.
 $$
 We have
 \begin {multline*}
 \`\Xi^{\Delta E}(g(b))\'
 -
 \<g^{(\Delta E)}\>(S)
 =
 \hfill
 \text{(since $\Xi^{\Delta E}(g(b))=g^{(\Delta E)}(\Xi^{\Delta E}(b))$)}
 \qquad
 \\
 \qquad
 =
 \<g^{(\Delta E)}\>(
 \`\Xi^{\Delta E}(b)\'
 -
 S
 )
 \in
 \hfill
 \text {(by Lemma~4.2)}
 \hfill
 \in
 \<(\~Y^{\~X})^{(\Delta E)}\>^{(r+1)}_{\~X}.
 \end {multline*}
 By Lemma~4.1,
 the ensemble $\<g^{(\Delta E)}\>(S)$ is fissile.
 Since $g$ is continuous,
 it takes $Y^X_a$ to $\~Y^{\~X}_{g(a)}$.
 Thus
 $$
 \<g^{(\Delta E)}\>(S)
 \in
 \<(\~Y^{\~X}_{g(a)})^{(\Delta E)}\>.
 $$
 We are done.
 \qed
 \end {demo}


 \head {\S~5. Compositions and smash products}


 \subhead {Compositions.}
 Let
 $X$,
 $Y$,
 $\~X$,
 and
 $\~Y$
 be cellular spaces,
 $X$
 and
 $\~X$
 compact.

 \begin {claim} [5.1. Corollary.]
 Let
 $k:\~X\to X$
 and
 $h:Y\to\~Y$
 be maps
 and
 $a,b:X\to Y$ be maps such that $a\%\approx rb$.
 Then
 $a\circ k\%\approx rb\circ k$ in $Y^{\~X}$
 and
 $h\circ a\%\approx rh\circ b$ in $\~Y^X$.
 \end {claim}

 \begin {demo} [Proof.]
 The transforms
 $$
 Y^X
 \to
 Y^{\~X},
 \qquad
 d
 \mapsto
 d\circ k,
 $$
 and
 $$
 Y^X
 \to
 \~Y^X,
 \qquad
 d
 \mapsto
 h\circ d,
 $$
 are primitive.
 By Lemma~4.3,
 they preserve strong $r$-similarity.
 \qed
 \end {demo}

 \begin {claim} [5.2. Corollary.]
 Let
 $k:\~X\to X$
 and
 $h:Y\to\~Y$
 be maps
 and
 $a:X\to Y$ be a map such that $\0\%\approx ra$.
 Then
 $\0\%\approx ra\circ k$ in $Y^{\~X}$
 and
 $\0\%\approx rh\circ a$ in $\~Y^X$.
 \end {claim}

 \begin {dem}
 Follows from Corollary~5.1.
 \end {dem}


 \subhead {Smash products.}
 Let
 $X$,
 $Y$,
 and
 $T$
 be cellular spaces,
 $X$
 and
 $T$
 compact.

 \begin {claim} [5.3. Corollary.]
 Let
 $a,b:X\to Y$ be maps such that $a\%\approx rb$.
 Then
 the maps
 $$
 a\sma\id_T,
 b\sma\id_T:
 X\sma T
 \to
 Y\sma T
 $$
 satisfy $a\sma\id_T\%\approx rb\sma\id_T$.
 \end {claim}

 \begin {demo} [Proof.]
 The transform
 $$
 Y^X
 \to
 (Y\sma T)^{X\sma T},
 \qquad
 d
 \mapsto
 d\sma\id_T,
 $$
 is primitive.
 By Lemma~4.3,
 it preserves strong $r$-similarity.
 \qed
 \end {demo}

 \begin {claim} [5.4. Corollary.]
 Let
 $a:X\to Y$ be a map such that $\0\%\approx ra$.
 Then
 the map
 $$
 a\sma\id_T:
 X\sma T
 \to
 Y\sma T
 $$
 satisfies $\0\%\approx ra\sma\id_T$.
 \end {claim}

 \begin {dem}
 Follows from Corollary~5.3.
 \end {dem}


 \head {\S~6. Homotopy invariance}


 Let
 $X$
 and
 $Y$
 be cellular spaces,
 $X$
 compact.

 \begin {claim} [6.1. Theorem.]
 Let maps $a,b,\~a,\~b:X\to Y$ satisfy
 $$
 \~a\sim a\%\approx rb\sim\~b.
 $$
 Then $\~a\%\approx r\~b$.
 \end {claim}

 \begin {demo} [Proof.]
 We
 crop $Y$
 and
 assume it compact.
 By \cite[Corollary~4.2]{sim},
 we can continuously associate
 to each path $v:[0,1]\to Y$
 an unbased homotopy $E_t(v):Y\to Y$, $t\in[0,1]$, such that
 $E_0(v)=\id$
 and
 $E_t(v)(v(0))=v(t)$.
 Let $h_t:X\to Y$, $t\in[0,1]$, be a homotopy such that
 $h_0=b$
 and
 $h_1=\~b$.
 For $x\in X$,
 introduce the path $v_x:[0,1]\to Y$,
 $t\mapsto h_t(x)$.
 We have
 $v_x(0)=h_0(x)=b(x)$
 and
 $v_x(1)=h_1(x)=\~b(x)$.
 Introduce the homotopy
 $$
 H_t:
 X\cro Y
 \to
 Y,
 \
 t\in[0,1],
 \qquad
 H_t(x,y)
 =
 E_t(v_x)(y).
 $$
 We have
 $$
 H_0(x,y)
 =
 E_0(v_x)(y)
 =
 y
 $$
 and
 $$
 H_1(x,b(x))
 =
 E_1(v_x)(b(x))
 =
 E_1(v_x)(v_x(0))
 =
 v_x(1)
 =
 \~b(x).
 $$

 Consider the primitive transforms
 $$
 g_t:
 Y^X
 \to
 Y^X,
 \
 t\in[0,1],
 \qquad
 g_t(d):x\mapsto H_t(x,d(x)).
 $$
 We have
 $d=g_0(d)\sim g_1(d)$,
 $d\in Y^X$,
 and
 $g_1(b)=\~b$.

 We have
 $$
 \~a
 \sim
 a
 \sim
 g_1(a)
 \%\approx r
 g_1(b)
 =
 \~b,
 $$
 where $\%\approx r$ holds by Lemma~4.3.
 By definition,
 the relation $\%\approx r$ tolerates homotopy of its left
 argument.
 Thus
 $\~a\%\approx r\~b$.
 \qed
 \end {demo}

 Using Theorem~6.1,
 we define the relation of strong $r$-similarity
 on the set $[X,Y]$
 by the rule
 $$
 [a]
 \%\approx r
 [b]
 \quad
 \Leftrightarrow
 \quad
 a
 \%\approx r
 b.
 $$


 \head {\S~7. More combining products}


 Let
 $X_1$,
 $X_2$,
 $Y$
 be spaces,
 $X_i$
 compact Hausdorff,
 and
 $U$ be an unbased space.
 We have the $\ZZ$-bilinear operations
 \begin {equation} \label {eBou}
 \eBou:
 \<Y^{X_1}\>
 \cro
 \<Y^{X_2}\>
 \to
 \<Y^{X_1\bou X_2}\>,
 \qquad
 \`d_1\'\eBou\`d_2\'
 =
 \`d_1\Bou d_2\',
 \end {equation}
 and
 \begin {multline*}
 \eBou_U:
 \<(Y^{X_1})^{(U)}\>
 \cro
 \<(Y^{X_2})^{(U)}\>
 \to
 \<(Y^{X_1\bou X_2})^{(U)}\>,
 \qquad
 \`w_1\'\eBou_U\`w_2\'
 =
 \`w\',
 \\
 w(u)
 =
 w_1(u)\Bou w_2(u):
 X_1\bou X_2
 \to
 Y,
 \qquad
 u\in U,
 \end {multline*}
 (combining products,
 cf.\ \eqref{bigeCop}).

 \begin {claim} [7.1. Lemma.]
 Let
 $E$ be a nonempty finite set
 and
 $S_i\in\<(Y^{X_i})^{(\Delta E)}\>$,
 $i=1,2$,
 be fissile ensembles.
 Then
 the ensemble
 $$
 S_1\eBou_{\Delta E}S_2
 \in
 \<(Y^{X_1\bou X_2})^{(\Delta E)}\>
 $$
 is fissile.
 \end {claim}

 \begin {demo} [Proof.]
 Take $A\in\AA(E)$.
 We have the commutative diagram
 $$
 \xymatrix @R1.4ex {
 \prod\limits_{F\in A}
 (
 \<(Y^{X_1})^{(\Delta E)}\>
 \cro
 \<(Y^{X_2})^{(\Delta E)}\>
 )
 \ar[rr]^-{
 \prod\limits_{F\in A}
 \eBou_{\Delta E}
 }
 \ar[dd]_-{
 \prod\limits_{F\in A}
 (
 \?|_{\Delta F}
 \cro
 \?|_{\Delta F}
 )
 }
 &&
 \prod\limits_{F\in A}
 \<(Y^{X_1\bou X_2})^{(\Delta E)}\>
 \ar[dd]^-{
 \prod\limits_{F\in A}
 \?|_{\Delta F}
 }
 \\ \\
 \prod\limits_{F\in A}
 (
 \<(Y^{X_1})^{(\Delta F)}\>
 \cro
 \<(Y^{X_2})^{(\Delta F)}\>
 )
 \ar[rr]^-{
 \prod\limits_{F\in A}
 \eBou_{\Delta F}
 }
 \ar@{=}[d]
 &&
 \prod\limits_{F\in A}
 \<(Y^{X_1\bou X_2})^{(\Delta F)}\>
 \ar[ddd]^-{
 \bigeCop\limits_{F\in A}
 }
 \\
 \prod\limits_{F\in A}
 \<(Y^{X_1})^{(\Delta F)}\>
 \cro
 \prod\limits_{F\in A}
 \<(Y^{X_2})^{(\Delta F)}\>
 \ar[dd]_-{
 \bigeCop\limits_{F\in A}
 \cro
 \bigeCop\limits_{F\in A}
 }
 &&
 \\ \\
 \<(Y^{X_1})^{(\Delta[A])}\>
 \cro
 \<(Y^{X_2})^{(\Delta[A])}\>
 \ar[rr]^-{
 \eBou_{\Delta[A]}
 }
 &&
 \<(Y^{X_1\bou X_2})^{(\Delta[A])}\>
 \\ \\
 \<(Y^{X_1})^{(\Delta E)}\>
 \cro
 \<(Y^{X_2})^{(\Delta E)}\>
 \ar[uu]^-{
 \?|_{\Delta[A]}
 \cro
 \?|_{\Delta[A]}
 }
 \ar[rr]^-{
 \eBou_{\Delta E}
 }
 &&
 \<(Y^{X_1\bou X_2})^{(\Delta E)}\>
 \ar[uu]_-{
 \?|_{\Delta[A]}
 }
 }
 $$
 with the sendings
 $$
 \xymatrix @R1.4ex {
 {\scriptstyle
 ((S_1,S_2))_{F\in A}
 }
 \ar@{|->}[rr]^-{(1)}
 \ar@{|->}[dd]_-{(2)}
 &&
 {\scriptstyle
 (S_1\eBou_{\Delta E}S_2)_{F\in A}
 }
 \ar@{|->}[dd]^-{(3)}
 \\ \\
 {\scriptstyle
 ((S_1|_{\Delta F},S_2|_{\Delta F}))_{F\in A}
 }
 \ar@{=}[d]
 &&
 {\scriptstyle
 ((S_1\eBou_{\Delta E}S_2)|_{\Delta F})_{F\in A}
 }
 \ar@{|->}[ddd]^-{(5)}
 \\
 {\scriptstyle
 ((S_1|_{\Delta F})_{F\in A},(S_2|_{\Delta F})_{F\in A})
 }
 \ar@{|->}[dd]_-{(4)}
 &&
 \\ \\
 {\scriptstyle
 (S_1|_{\Delta[A]},S_2|_{\Delta[A]})
 }
 &&
 {\scriptstyle
 \bigeCop\limits_{F\in A}
 (S_1\eBou_{\Delta E}S_2)|_{\Delta F}
 }
 \\ \\
 {\scriptstyle
 (S_1,S_2)
 }
 \ar[uu]^-{(6)}
 \ar[rr]^{(8)}
 &&
 {\scriptstyle
 S_1\eBou_{\Delta E}S_2.
 }
 \ar[uu]_-{(7)}
 }
 $$
 The sending (4) holds by fissility of
 $S_1$
 and
 $S_2$.
 The sendings (1), (2), (3), (5), (6), and (8) are obvious.
 The sending (7) follows.
 Thus $S_1\eBou_{\Delta E}S_2$ is fissile.
 \qed
 \end {demo}

 \begin {claim} [7.2. Lemma.]
 We have
 $$
 \<(Y^{X_1})^{(U)}\>^{(p)}_{X_1}
 \eBou_U
 \<(Y^{X_2})^{(U)}\>^{(q)}_{X_2}
 \subseteq
 \<(Y^{X_1\bou X_2})^{(U)}\>^{(p+q)}_{X_1\bou X_2}.
 $$
 \end {claim}

 \begin {demo} [Proof.]
 Take ensembles
 \begin {equation} \label {Zi}
 W_1\in\<(Y^{X_1})^{(U)}\>^{(p)}_{X_1},
 \qquad
 W_2\in\<(Y^{X_2})^{(U)}\>^{(q)}_{X_2}.
 \end {equation}
 We have the commutative diagram
 $$
 \xymatrix @C1.4ex {
 \<(Y^{X_1})^{(U)}\>
 \cro
 \<(Y^{X_2})^{(U)}\>
 \ar[rrr]^-{\eBou_U}
 \ar[d]_-{
 \<\uc^{X_1}\>
 \cro
 \<\uc^{X_2}\>
 }
 &&&
 \<(Y^{X_1\bou X_2})^{(U)}\>
 \ar[d]^-{
 \<\uc^{X_1\bou X_2}\>
 }
 \\
 \<Y^{U\hsm X_1}\>
 \cro
 \<Y^{U\hsm X_2}\>
 \ar[rr]^-{\eBou}
 &&
 \<Y^{(U\hsm X_1)\bou(U\hsm X_2)}\>
 \ar@{=}[r]
 &
 \<Y^{U\hsm(X_1\bou X_2)}\>
 }
 $$
 with sendings
 $$
 \xymatrix @C1.4ex {
 {\scriptstyle
 (W_1,W_2)
 }
 \ar@{|->}[rrr]
 \ar@{|->}[d]
 && &
 {\scriptstyle
 W_1\eBou_UW_2
 }
 \ar@{|->}[d]
 \\
 {\scriptstyle
 (
 \<\uc^{X_1}\>(W_1),
 \<\uc^{X_2}\>(W_2)
 )
 }
 \ar@{|->}[rr]
 &&
 {\scriptstyle
 \<\uc^{X_1}\>(W_1)
 \eBou
 \<\uc^{X_2}\>(W_2)
 }
 \ar@{=}[r]
 &
 {\scriptstyle
 \<\uc^{X_1\bou X_2}\>
 (W_1\eBou_UW_2).
 }
 }
 $$
 It follows from \eqref{Zi} that
 $$
 (
 \<\uc^{X_1}\>(W_1),
 \<\uc^{X_2}\>(W_2)
 )
 \in
 \<Y^{U\hsm X_1}\>^{(p)}
 \cro
 \<Y^{U\hsm X_2}\>^{(q)}.
 $$
 Thus,
 by \cite[Lemma~3.1]{sim-1-cir},
 $$
 \<\uc^{X_1}\>(W_1)
 \eBou
 \<\uc^{X_2}\>(W_2)
 \in
 \<Y^{(U\hsm X_1)\bou(U\hsm X_2)}\>^{(p+q)}.
 $$
 Equivalently,
 $$
 \<\uc^{X_1\bou X_2}\>
 (W_1\eBou_UW_2)
 \in
  \<Y^{U\hsm(X_1\bou X_2)}\>^{(p+q)}.
 $$
 Thus
 $$
 W_1\eBou_UW_2
 \in
  \<(Y^{X_1\bou X_2})^{(U)}\>^{(p+q)}_{X_1\bou X_2}.
 \QED
 $$
 \end {demo}

 Let
 $X_1$,
 $X_2$,
 and
 $Y$
 be cellular spaces,
 $X_i$
 compact.

 \begin {claim} [7.3. Corollary.]
 Let maps $a_i:X_i\to Y$,
 $i=1,2$,
 satisfy $\0\%\approx ra_i$.
 Then
 the map
 $$
 a_1\Bou a_2:
 X_1\bou X_2
 \to
 Y
 $$
 satisfies $\0\%\approx ra_1\Bou a_2$.
 \end {claim}

 \begin {demo} [Proof.]
 Take a nonempty finite set $E$.
 We have fissile ensembles $S_i\in\<(Y^{X_i}_\0)^{(\Delta E)}\>$,
 $i=1,2$,
 such that
 $$
 \`\Xi^{\Delta E}(a_i)\'-S_i
 \in
 \<(Y^{X_i})^{(\Delta E)}\>^{(r+1)}_{X_i}.
 $$
 By Lemma~7.1,
 the ensemble
 $$
 S_1\eBou_{\Delta E}S_2
 \in
 \<(Y^{X_1\bou X_2}_\0)^{(\Delta E)}\>
 $$
 is fissile.
 We have
 \begin {multline*}
 \`\Xi^{\Delta E}(a_1\Bou a_2)\'
 -
 S_1\eBou_{\Delta E}S_2
 =
 \\
 =
 \`\Xi^{\Delta E}(a_1)\'
 \eBou_{\Delta E}
 \`\Xi^{\Delta E}(a_2)\'
 -
 S_1\eBou_{\Delta E}S_2
 =
 \\
 =
 (\`\Xi^{\Delta E}(a_1)\'-S_1)
 \eBou_{\Delta E}
 \`\Xi^{\Delta E}(a_2)\'
 +
 {}
 \\
 +
 S_1
 \eBou_{\Delta E}
 (\`\Xi^{\Delta E}(a_2)\'-S_2)
 \in
 \<(Y^{X_1\bou X_2})^{(\Delta E)}\>^{(r+1)}_{X_1\bou X_2},
 \end {multline*}
 where $\in$ holds by Lemma~7.2.
 We are done.
 \qed
 \end {demo}


 \head {\S~8. Strong similarity for an admissible couple}


 Let $X$ and $Y$ be cellular spaces,
 $X$ compact.
 Let $X$ be equipped with maps
 $\mu:X\to X\bou X$ (comultiplication) and
 $\nu:X\to X$ (coinversion).
 The set $Y^X$ carries the operations
 $$
 (a,b)
 \mapsto
 (a\cat b:X\xto{\mu}X\bou X\xto{a\Bou b}Y)
 $$
 and
 $$
 a
 \mapsto
 (a^\dag:X\xto{\nu}X\xto{a}Y).
 $$
 We suppose that
 $(X,\mu,\nu;Y)$ is an admissible couple
 in the sense of \cite{sim-1-cir},
 that is,
 the set $[X,Y]$ is a group with
 the multiplication
 $$
 [a][b]=[a\cat b],
 $$
 the inversion
 $$
 [a]^{-1}=[a^\dag],
 $$
 and
 the identity $1=[\0^X_Y]$.
 We are mainly interested in the case of
 $X=\Sigma T$ with standard
 $\mu$
 and
 $\nu$.

 We proceed parallelly to \cite{sim-1-cir}.
 The subsets
 $$
 [X,Y]^{((r+1))}
 =
 \{\,\+a\in[X,Y]\mid1\%\approx r\+a\,\}
 $$
 form the filtration
 $$
 [X,Y]
 =
 [X,Y]^{((1))}
 \supseteq
 [X,Y]^{((2))}
 \supseteq
 \dotso.
 $$

 \begin {claim} [8.1. Theorem.]
 $[X,Y]^{((r+1))}\subseteq[X,Y]$ is a normal subgroup.
 \end {claim}

 \begin {demo} [Proof.]
 Take $a,b:X\to Y$,
 $\0\%\approx ra,b$.
 Check that
 $\0\%\approx ra*b$.
 We have the decomposition
 $$
 a*b:
 X
 \xto{\mu}
 X\bou X
 \xto{a\Bou b}
 Y.
 $$
 By Corollary~7.3,
 $\0^{X\bou X}_Y\%\approx ra\Bou b$.
 By Corollary~5.2,
 $\0\%\approx ra*b$.

 Take $a:X\to Y$,
 $\0\%\approx ra$.
 Check that
 $\0\%\approx ra^\dag$.
 We have the decomposition
 $$
 a^\dag:
 X
 \xto{\nu}
 X
 \xto{a}
 Y.
 $$
 By Corollary~5.2,
 $\0\%\approx ra^\dag$.

 Take $a,b:X\to Y$,
 $\0\%\approx ra$.
 Check that
 $\0\%\approx rb^\dag*(a*b)$.
 Consider the primitive transform
 $$
 Y^X
 \to
 Y^X,
 \qquad
 d
 \mapsto
 b^\dag*(d*b).
 $$
 We have
 $$
 \0^X_Y
 \sim
 b^\dag*(\0^X_Y*b)
 \%\approx r
 b^\dag*(a*b),
 $$
 where $\%\approx r$ holds by Lemma~4.3.
 By
 (the trivial part of)
 Theorem~6.1,
 $\0\%\approx rb^\dag*(a*b)$.
 \qed
 \end {demo}

 We do not know whether
 the subgroups $[X,Y]^{((s))}$ form an N-series.

 \begin {claim} [8.2. Theorem.]
 For $\+a,\+b\in[X,Y]$,
 we have
 $$
 \+a\%\approx r\+b
 \quad
 \Leftrightarrow
 \quad
 \+a^{-1}\+b\in[X,Y]^{((r+1))}.
 $$
 \end {claim}

 \begin {demo} [Proof.]
 It suffices to check that,
 for maps $a,b,c:X\to Y$,
 $a\%\approx rb$
 implies
 $c*a\%\approx rc*b$.
 This follows from Lemma~4.3 for the primitive transform
 $$
 Y^X
 \to
 Y^X,
 \qquad
 d
 \mapsto
 c*d.
 \QED
 $$
 \end {demo}

 It follows from Theorems 8.1 and 8.2 that,
 for an admissible couple $(X,\mu,\nu;Y)$,
 the relation $\%\approx r$ on $[X,Y]$ is an equivalence.


 \head {\S~9. Presheaves and extenders}


 Let
 $P$ be a finite partially ordered set
 and
 $C$ be a concrete category.
 (Concreteness is not essential;
 we assume it for convenience of notation only.)
 A cofunctor $S:P\to C$ is called a {\it presheaf}.
 For $p,q\in P$, $p\ge q$,
 we have the induced morphism
 $$
 \?|_q:
 S(p)
 \to
 S(q)
 $$
 (the {\it restriction\/} morphism).

 For a preasheaf $U:P\to\Ab$,
 we have the isomorphism
 $$
 \nabla_P:
 \bigoplus_{p\in P}
 U(p)
 \to
 \bigoplus_{p\in P}
 U(p),
 \qquad
 \ins_p(u)
 \mapsto
 \sum_{q\in P\lceil p\rceil}
 \ins_q(u|_q),
 \qquad
 u\in U(p),
 \
 p\in P.
 $$
 Hereafter,
 $$
 P\lceil p\rceil
 =
 \{\,q\in P\mid p\ge q\,\}
 $$
 and
 $$
 \ins_q:
 U(q)
 \to
 \bigoplus_{p\in P}
 U(p)
 $$
 are the canonical insertions.

 Suppose that
 $P$ has
 the infimum operation $\mee$ 
 and
 the greatest element $\top$.
 It follows that
 $P$ is a lattice.
 We put $P^\x=P\setminus\{\top\}$.
 An {\it extender\/} $\lambda$ for the preasheaf $S$ is a
 collection of morphisms
 $$
 \lambda^q_p:
 S(q)
 \to
 S(p),
 \qquad
 p,q\in P,
 \
 p\ge q,
 $$
 such that,
 for
 $p,q\in P$
 and
 $s\in S(q)$,
 $$
 \lambda^q_p(s)|_q
 =
 s
 \qquad
 \text{if $p\ge q$}
 $$
 and
 $$
 \lambda^q_\top(s)|_p
 =
 \lambda^{p\mee q}_p(s|_{p\mee q}).
 $$
 In particular,
 $$
 \lambda^q_p(s)
 =
 \lambda^q_\top(s)|_{p}.
 $$
 (The extenders we deal with satisfy the identity
 $\lambda^q_p\circ\lambda^r_q=\lambda^r_p$.
 We
 neither check
 nor use
 this.)

 Consider a preaheaf $U:P\to\Ab$ with an extender $\lambda$.
 The symbol $\Oplus$ below denotes the homomorphism of a direct
 sum given by its restrictions to the summands.

 \begin {claim} [9.1. Lemma.]
 For $q\in P$,
 the diagram
 $$
 \xymatrix {
 \bigoplus\limits_{p\in P}
 U(p)
 \ar[d]_-{\pro}
 &&
 \bigoplus\limits_{p\in P}
 U(p)
 \ar[ll]_-{\nabla_P}^-{\cong}
 \ar[rr]^-{\bigOplus\limits_{p\in P}
           \lambda^p_\top}
 \ar[d]_-{R_q}
 &&
 U(\top)
 \ar[d]^-{\?|_q}
 \\
 \bigoplus\limits_{p\in P\lceil q\rceil}
 U(p)
 &&
 \bigoplus\limits_{p\in P\lceil q\rceil}
 U(p)
 \ar[ll]_-{\nabla_{P\lceil q\rceil}}^-{\cong}
 \ar[rr]^-{\bigOplus\limits_{p\in P\lceil q\rceil}
           \lambda^p_q}
 &&
 U(q),
 }
 $$
 where $R_q$ is the homomorphism defined by the rule
 $$
 \ins_p(u)
 \mapsto
 \ins_{p\mee q}(u|_{p\mee q}),
 $$
 is commutative.
 \end {claim}

 \begin {dem}
 Direct check.
 \qed
 \end {dem}

 \begin {claim} [9.2. Lemma.]
 The homomorphism
 $$
 U(\top)
 \to
 \lim_{p\in P^\x}
 U(p),
 \qquad
 u
 \mapsto
 (u|_p)_{p\in P^\x},
 $$
 is surjective.
 \end {claim}

 \begin {demo} [Proof.]
 Take a collection
 $$
 (u_p)_{p\in P^\x}
 \in
 \lim_{p\in P^\x}
 U(p)
 \subseteq
 \bigoplus_{p\in P^\x}
 U(p).
 $$
 Define
 a collection $(v_p)_{p\in P^\x}$
 and
 a section $u$
 by the diagram
 $$
 \xymatrix @R1.4ex {
 \bigoplus\limits_{p\in P^\x}
 U(p)
 &&
 \bigoplus\limits_{p\in P^\x}
 U(p)
 \ar[ll]_-{\nabla_{P^\x}}^-{\cong}
 \ar[rr]^-{\bigOplus\limits_{p\in P^\x}
           \lambda^p_\top}
 &&
 U(\top).
 \\
 {\scriptstyle(u_p)_{
 p\in P^\x}
 }
 &&
 {\scriptstyle
 (v_p)_{p\in P^\x}
 }
 \ar@{|->}[ll]
 \ar@{|->}[rr]
 &&
 {\scriptstyle
 u
 }
 }
 $$
 Take $q\in P^\x$.
 We show that
 $u|_q=u_q$,
 which will suffice.
 In the diagram of Lemma~9.1,
 we have
 $$
 \xymatrix {
 {\scriptstyle
 (u_p)_{p\in P}
 }
 \ar[d]_-{(3)}
 &&
 {\scriptstyle
 (v_p)_{p\in P}
 }
 \ar@{|->}[ll]_-{(1)}
 \ar@{|->}[rr]^-{(2)}
 \ar@{|->}[d]_-{(4)}
 &&
 {\scriptstyle
 u
 }
 \ar@{|->}[d]^-{(5)}
 \\
 {\scriptstyle
 (u_p)_{p\in P\lceil q\rceil}
 }
 &&
 {\scriptstyle
 \ins_q(u_q)
 }
 \ar@{|->}[ll]_-{(6)}
 \ar@{|->}[rr]^-{(7)}
 &&
 {\scriptstyle
 u_q,
 }
 }
 $$
 where we put $u_\top=v_\top=0$ in $U(\top)$.
 The sendings $(1)$ and $(2)$ follow from the construction of
 the collections.
 The sending $(6)$ expresses the equalities $u_q|_p=u_p$,
 $p\in P\lceil q\rceil$,
 which hold by the definition of limit.
 The sending $(3)$ is obvious.
 The sending $(4)$ follows because
 the left square is commutative
 and
 $\nabla_{P\lceil q\rceil}$ is injective.
 The sending $(7)$ follows from the equality $\lambda^q_q=\id$.
 By commutativity of the right square,
 the sending $(5)$ holds,
 which is what was to be checked.
 \qed
 \end {demo}


 \head {\S~10. The abstract fissilizer $\Phi_E$ on $\<\_M(E)\>$}


 Fix a nonempty finite set $E$.
 The set $\PP_\x(E)$ is partially ordered by inclusion.

 For $A,B\in\AA(E)$,
 we say $A\ge B$ if,
 for each $G\in B$,
 there is $F\in A$ such that
 $F\supseteq G$.
 The set $\AA(E)$ becomes a lattice with
 the infimum operation
 $$
 A\mee B
 =
 \{\,F\cap G\mid F\in A,\ G\in B\,\}
 \setminus
 \{\varnothing\}
 $$
 and
 the greatest element $\top=\{E\}$.

 Let
 $\Mg$ be the category of sets
 and
 $M:\PP_\x(E)\to\Mg$ be a presheaf.
 We define a presheaf $\_M:\AA(E)\to\Mg$.
 For $A\in\AA(E)$,
 put
 $$
 \_M(A)
 =
 \prod_{F\in A}
 M(F).
 $$
 For $A,B\in\AA(E)$, $A\ge B$,
 define the restriction function
 $$
 \_M(A)
 \to
 \_M(B),
 \qquad
 \_m
 \mapsto
 \_m|_B,
 $$
 by putting,
 for $\_m=(m_F)_{F\in A}$,
 $$
 \_m|_B
 =
 (m_{(A)G}|_G)_{G\in B},
 $$
 where $(A)G\in A$ is the unique set
 that includes $G$.
 Clearly,
 $\_M(\{E\})=M(E)$.

 Taking composition with the functor $\<\?\>:\Mg\to\Ab$,
 we get the presheaves
 $$
 \PP_\x(E)
 \to
 \Ab,
 \qquad
 F
 \mapsto
 \<M(F)\>,
 $$
 and
 \begin {equation} \label {ab}
 \AA(E)
 \to
 \Ab,
 \qquad
 A
 \mapsto
 \<\_M(A)\>.
 \end {equation}
 For $A\in\AA(E)$,
 we have the $\ZZ$-multilinear operation
 \begin {equation} \label {bigecol}
 \bigecol_{F\in A}:
 \prod_{F\in A}
 \<M(F)\>
 \to
 \<\_M(A)\>,
 \qquad
 \bigecol_{F\in A}
 \`m_F\'
 =
 \`(m_F)_{F\in A}\',
 \end {equation}
 (combining product,
 cf.\ \eqref{bigeCop}).
 For
 $Q\in\<M(E)\>$
 and
 $A\in\AA(E)$,
 put
 $$
 Q^\secol(A)
 =
 \bigecol_{F\in A}
 Q|_F
 \in
 \<\_M(A)\>.
 $$
 Note that
 $Q^\secol(\varnothing)=\`\atom\'$,
 where $\atom\in\_M(\varnothing)$ is the unique element.
 We call an ensemble $R\in\<M(E)\>$
 {\it fissile\/} if,
 for any layout $A\in\AA(E)$,
 $$
 R|_A
 =
 R^\secol(A)
 $$
 in $\<\_M(A)\>$
 (cf.\ \S~2).

 We suppose that
 the presheaf $\_M$ has an extender
 $$
 \lambda^B_A:
 \_M(B)
 \to
 \_M(A),
 \qquad
 A,B\in\AA(E),
 \
 A\ge B.
 $$
 Then
 the preasheaf \eqref{ab} has the extender
 $$
 \<\lambda^B_A\>:
 \<\_M(B)\>
 \to
 \<\_M(A)\>,
 \qquad
 A,B\in\AA(E),
 \
 A\ge B.
 $$

 For $Q\in\<M(E)\>$,
 define an ensemble $\Phi_E(Q)\in\<M(E)\>$
 by the rule
 $$
 \xymatrix @R1.4ex {
 \bigoplus\limits_{A\in\AA(E)}
 \<\_M(A)\>
 &&
 \bigoplus\limits_{A\in\AA(E)}
 \<\_M(A)\>
 \ar[ll]_-{\nabla_{\AA(E)}}^-{\cong}
 \ar[rr]^-{
 \bigOplus\limits_{A\in\AA(E)}
 \<\lambda^A_{\{E\}}\>
 }
 &&
 \<M(E)\>.
 \\
 {\scriptstyle
 Q^\secol
 }
 &&
 {\scriptstyle
 \nabla_{\AA(E)}^{-1}(Q^\secol)
 }
 \ar@{|->}[ll]
 \ar@{|->}[rr]
 &&
 {\scriptstyle
 \Phi_E(Q)
 }
 }
 $$
 We get a function
 (not a homomorphism)
 $$
 \Phi_E:\<M(E)\>\to\<M(E)\>,
 $$
 which we call the {\it fissilizer}.

 \begin {claim} [10.1. Lemma.]
 For any ensemble $Q\in\<M(E)\>$,
 the ensemble $\Phi_E(Q)$ is fissile.
 \end {claim}

 \begin {demo} [Proof.]
 Take $A\in\AA(E)$.
 We have the commutative diagram
 $$
 \xymatrix @R1.4ex {
 **{!<8ex,0ex>}
 {\scriptstyle
 Q^\secol
 }
 \ar@{|->}@/_1ex/[]!/l8ex/;[dddddd]!/l8ex/
 &&
 {\scriptstyle
 \nabla_{\AA(E)}^{-1}(Q^\secol)
 }
 \ar@{|->}[ll]!/l8ex/
 \ar@{|->}[rr]!/r3ex/
 &&
 **{!<-3ex,0ex>}
 {\scriptstyle
 \Phi_E(Q)
 }
 \ar@{|->}@/^1ex/[]!/r5ex/;[dddddd]!/r5ex/
 \\
 \bigoplus\limits_{a\in\AA(E)}
 \<\_M(a)\>
 \ar[dd]_-{\pro}
 &&
 \bigoplus\limits_{a\in\AA(E)}
 \<\_M(a)\>
 \ar[ll]_-{\nabla_{\AA(E)}}^-{\cong}
 \ar[rr]^-{
 \bigOplus\limits_{a\in\AA(E)}
 \<\lambda^a_{\{E\}}\>
 }
 \ar[dd]_-{R_A}
 &&
 \<M(E)\>
 \ar[dd]^-{\?|_A}
 \\ \\
 \bigoplus\limits_{a\in\AA(E)\lceil A\rceil}
 \<\_M(a)\>
 \ar[dd]_-{J_A}^-{\cong}
 &&
 \bigoplus\limits_{a\in\AA(E)\lceil A\rceil}
 \<\_M(a)\>
 \ar[ll]_-{\nabla_{\AA(E)\lceil A\rceil}}^-{\cong}
 \ar[rr]^-{
 \bigOplus\limits_{a\in\AA(E)\lceil A\rceil}
 \<\lambda^a_A\>
 }
 \ar[dd]_-{J_A}^-{\cong}
 &&
 \<\_M(A)\>
 \ar[dd]^-{I_A}_-{\cong}
 \\ \\
 \bigotimes\limits_{F\in A}
 \
 \bigoplus\limits_{b\in\AA(F)}
 \<\_M(b)\>
 &&
 \bigotimes\limits_{F\in A}
 \
 \bigoplus\limits_{b\in\AA(F)}
 \<\_M(b)\>
 \ar[ll]_-{
 \bigotimes\limits_{F\in A}
 \nabla_{\AA(F)}
 }^-{\cong}
 \ar[rr]^*!/r2ex/{\labelstyle
 \bigotimes\limits_{F\in A}
 \
 \bigOplus\limits_{b\in\AA(F)}
 \<\lambda^b_{\{F\}}\>
 }
 &&
 \bigotimes\limits_{F\in A}
 \<M(F)\>,
 \\
 **{!<6ex,0ex>}
 {\scriptstyle
 \bigotimes\limits_{F\in A}
 Q|_F^\secol
 }
 &&
 {\scriptstyle
 \bigotimes\limits_{F\in A}
 \nabla_{\AA(F)}^{-1}(Q|_F^\secol)
 }
 \ar@{|->}[ll]!/l6ex/
 \ar@{|->}[rr]!/r1ex/
 &&
 **{!<-1ex,0ex>}
 {\scriptstyle
 \bigotimes\limits_{F\in A}
 \Phi_F(Q|_F)
 }
 }
 $$
 where
 the upper half comes from Lemma~9.1,
 $I_A$ is the isomorphism defined by the rule
 $$
 \`(m_F)_{F\in A}\'
 \mapsto
 \bigotimes_{F\in A}
 \`m_F\',
 $$
 and
 $J_A$ is the isomorphism defined by the rule
 $$
 \ins_a(\`\_m\')
 \mapsto
 \bigotimes_{F\in A}
 \ins_{a\mee\{F\}}(\`\_m|_{a\mee\{F\}}\')
 $$
 (note that
 $a\mee\{F\}\in\AA(F)\subseteq\AA(E)$).
 Commutativity of the lower half is checked directly.
 The sendings in the upper row hold
 by the definition of $\Phi_E$.
 The sendings in the lower row hold
 by the definition of $\Phi_F:\<M(F)\>\to\<M(F)\>$.
 The sending in the left column is checked directly.
 The sending in the right column follows.
 Since
 $$
 I_A:
 \bigecol_{F\in A}
 q_F
 \mapsto
 \bigotimes_{F\in A}
 q_F
 $$
 for $q_F\in\<M(F)\>$,
 $F\in A$,
 we get
 $$
 \Phi_E(Q)|_A
 =
 \bigecol_{F\in A}
 \Phi_F(Q|_F).
 $$
 In particular,
 for $A=\{F\}$,
 this gives
 $$
 \Phi_E(Q)|_F
 =
 \Phi_F(Q|_F).
 $$
 Thus,
 for arbitrary $A$,
 $$
 \Phi_E(Q)|_A
 =
 \bigecol_{F\in A}
 \Phi_E(Q)|_F.
 $$
 Thus
 $\Phi_E(Q)$ is fissile.
 \qed
 \end {demo}

 Let $N(A)\subseteq\<\_M(A)\>$, $A\in\AA(E)$, be a collection
 of subgroups preserved by
 the restriction homomorphisms
 and
 the homomorphisms $\<\lambda^B_A\>$.

 \begin {claim} [10.2. Lemma.]
 Let an ensemble $Q\in\<M(E)\>$ satisfy
 $$
 Q^\secol(A)
 -
 Q|_A
 \in
 N(A)
 $$
 for all $A\in\AA(E)$.
 Then
 $$
 \Phi_E(Q)-Q
 \in
 N(\{E\}).
 $$
 \end {claim}

 \begin {demo} [Proof.]
 We have the presheaf
 $$
 \AA(E)
 \to
 \Ab,
 \qquad
 A
 \mapsto
 \<\_M(A)\>/N(A),
 $$
 with the induced restriction homomorphisms.
 We have the commutative diagram
 $$
 \xymatrix @R1.4ex {
 {\scriptstyle
 (Q|_A)_{A\in\AA(E)}
 }
 &
 {\scriptstyle
 \ins_{\{E\}}(Q)
 }
 \ar@{|->}[l]
 \ar@{|->}[r]
 &
 {\scriptstyle
 Q
 }
 \\
 {\scriptstyle
 Q^\secol
 }
 &
 {\scriptstyle
 \nabla_{\AA(E)}^{-1}(Q^\secol)
 }
 \ar@{|->}[l]
 \ar@{|->}[r]
 &
 {\scriptstyle
 \Phi_E(Q)
 }
 \\
 \bigoplus\limits_{A\in\AA(E)}
 \<\_M(A)\>
 \ar[dd]_-{\pro}
 &
 \bigoplus\limits_{A\in\AA(E)}
 \<\_M(A)\>
 \ar[l]_-{\nabla_{\AA(E)}}^-{\cong}
 \ar[r]^-{
 \bigOplus\limits_{A\in\AA(E)}
 \<\lambda^A_{\{E\}}\>
 }
 \ar[dd]_-{\pro}
 &
 \<M(E)\>
 \ar[dd]^-{\pro}
 \\ \\
 \bigoplus\limits_{A\in\AA(E)}
 \<\_M(A)\>/N(A)
 &
 \bigoplus\limits_{A\in\AA(E)}
 \<\_M(A)\>/N(A)
 \ar[l]_-{\nabla_{\AA(E)}}^-{\cong}
 \ar[r]
 &
 \<M(E)\>/N(\{E\}).
 }
 $$
 The upper line of sendings is obvious.
 The lower line of sendings holds by the definition of
 $\Phi_E$.
 By hypothesis,
 the difference of the elements in the upper-left corner
 descends to zero.
 Since $\nabla_{\AA(E)}$ in the lower row is an isomorphism,
 the difference of elements in the upper-right corner
 also descends to zero.
 \qed
 \end {demo}


 \head {\S~11. Topological and simplicial constructions}


 \subhead {Topological cones.}
 Given an unbased space $U$,
 we have the space  $U_+=U\sqcup\{\0\}$.
 Take $s\in\{0,1\}$
 and
 form the space
 $$
 \oC^sU
 =
 (U\cro[0,1])/(U\cro\{s\}),
 $$
 the {\it cone\/} over $U$.
 The innate basepoint
 (where $U\cro\{s\}$ is projected)
 is called the {\it apex}.
 Using the ``base'' embedding
 $$
 U
 \xto{u\mapsto(u,1-s)}
 U\cro[0,1]
 \xto{\pro}
 \oC^sU,
 $$
 we adopt
 the inclusion $U\subseteq\oC^sU$
 and
 the based one $U_+\subseteq\oC^sU$.
 A path of the form
 $$
 [0,1]
 \xto{t\mapsto(u,t)}
 U\cro[0,1]
 \xto{\pro}
 \oC^sU
 $$
 is called a {\it generating path}.
 For an unbased subspace $V\subseteq U$,
 we have $\oC^sV\subseteq\oC^sU$.

 Notation:
 $\vC=\oC^0$,
 $\hC=\oC^1$.


 \subhead {Topological suspensions.}
 For an unbased space $U$,
 the {\it unreduced suspension\/} $\bS U$ is the colimit
 of the diagram
 $$
 \{0,1\}
 \xfrom{\pro}
 U\cro\{0,1\}
 \xto{\inc}
 U\cro[0,1].
 $$
 Let $s_{\bS U}\in\bS U$ be the point coming from
 $s\in\{0,1\}$.
 We let $0_{\bS U}$ be the basepoint of $\bS U$.

 We use also the usual reduced suspension $\Sigma$.


 \subhead {Simplicial notation.}
 The simplex category consists of the sets
 $[n]=\{0,\dotsc,n\}$.
 Let
 $\U_{[n]}$ denote the $n$th term of a simplicial set $\U$
 and
 $\f_{[n]}:\U_{[n]}\to\V_{[n]}$ be the $n$th term of a morphism
 $\f:\U\to\V$.
 Let $\D^n$ be ``the $n$-simplex'',
 the simplicial set represented
 (as a cofunctor)
 by the object $[n]$.
 We have the based simplicial set $\U_+=\U\sqcup\D^0$.


 \subhead {Unreduced Kan cones.}
 Take $s\in\{0,1\}$.
 Let $\dg^s:\D^0\to\D^1$ be the morphism induced by the
 function $\delta^s:[0]\to[1]$,
 $0\mapsto1-s$.
 Given a simplicial set $\U$,
 we define its {\it cone\/} $\OC^s\U$.
 There is a unique
 (up to an isomorphism)
 Cartesian square
 $$
 \xymatrix {
 \U
 \ar[r]^-{\i}
 \ar[d]
 &
 \OC^s\U
 \ar[d]^-{\p}
 \\
 \D^0
 \ar[r]^-{\dg^s}
 &
 \D^1
 }
 $$
 with the universal property expressed by the diagram
 $$
 \xymatrix @R1.4ex {
 \U
 \ar[rrr]^-{\i}
 \ar[dd]
 &&&
 \OC^s\U
 \ar[dd]^-{\p}
 \\
 &
 \A
 \ar[r]
 \ar[ul]
 \ar[dl]
 &
 \B
 \ar@{-->}[ru]
 \ar[dr]
 &
 \\
 \D^0
 \ar[rrr]^-{\dg^s}
 &&&
 \D^1,
 }
 $$
 where the lower trapeze is assumed to be Cartesian\footnote
 {Thoroughly,
 the cone $\OC^s$ is the functor $\sMg\to\sMg/\D^1$
 ($\sMg$ is the category of simplicial sets)
 right adjoint right inverse to
 the functor of pullback along $\dg^s:\D^0\to\D^1$.
 (Suggested by I.~S.~Baskov.)}.
 The morphism $\dg^{1-s}:\D^0\to\D^1$ lifts along $\p$
 uniquely.
 This yields a morphism $\D^0\to\OC^s\U$,
 which makes $\OC^s\U$ a based simplicial set.
 The basepoint is called the {\it apex}.
 The morphism $\i$ is injective.
 Using it,
 we adopt
 the inclusion $\U\subseteq\OC^s\U$
 and
 the based one $\U_+\subseteq\OC^s\U$.
 We call $\p$ the {\it projection}.

 All constructions are covariant/natural in $\U$.
 The functor $\OC^s$ preserves injective morphisms.
 Using this,
 we adopt the inclusion $\OC^s\V\subseteq\OC^s\U$
 for a simplicial subset $\V\subseteq\U$.

 Notation:
 $\VC=\OC^0$,
 $\HC=\OC^1$.

 There is a unique natural map $r:\oC^s|\U|\to|\OC^s\U|$ such that
 the diagram
 $$
 \xymatrix {
 &&
 \oC^s|\U|
 \ar[d]^-{r}
 \\
 |\U|
 \ar[rr]_-{\text{$\inc$ ($=|\i|$)}}
 \ar[urr]^-{\inc}
 &&
 |\OC^s\U|
 }
 $$
 is commutative
 and
 each generating path of $\oC^s|\U|$ is sent to an affine path
 in some simplex of $|\OC^s\U|$.
 The map $r$ is a homeomorphism.
 Using it,
 we adopt that
 $|\OC^s\U|=\oC^s|\U|$.


 \subhead {The reduced Kan cone.}
 For a based simplicial set $\T$,
 introduce the based simplicial set
 $\Vc\T=\VC\T/\VC(\0)$,
 where $(\0)\subseteq\T$ is the simplicial subset
 generated by the basepoint $\0\in\T_{[0]}$
 (so,
 $(\0)\cong\D^0$).
 We adopt the obvious
 inclusion $\T\subseteq\Vc\T$
 and
 identification $\Vc(\U_+)=\VC\U$.
 $\Vc$ is a functor;
 it preserves wedges.


 \subhead {The unreduced Kan suspension.}
 For a simplicial set $\U$,
 introduce the based simplicial set
 $\HS\U=\HC\U/\U$.
 It has two vertices:
 the top $1_{\HS\U}$,
 which is the image
 of the apex of the cone $\HC\U$
 under the projection $\HC\U\to\HS\U$,
 and
 the basepoint $0_{\HS\U}$
 (where the base $\U\subseteq\HC\U$ is sent).
 We have
 $$
 |\HS\U|
 =
 |\HC\U|/|\U|
 =
 \hC|\U|/|\U|
 =
 \bS|\U|.
 $$


 \subhead {The thick simplex.}
 For a set $A$,
 let $\E A$ be the simplicial set with
 $(\E A)_{[n]}=A^{[n]}$
 ($=A^{n+1}$)
 and obvious structure functions.

 For each $a\in A$,
 there is a unique retraction $\~\sg_a:\VC\HC\E A\to\HC\E A$
 sending
 the apex
 to the vertex $a\in A=(\E A)_{[0]}\subseteq(\HC\E A)_{[0]}$.
 Define retractions
 $\-\sg_a$
 and
 $\sg_a$
 by the commutative diagram
 $$
 \xymatrix {
 \VC\HC\E A
 \ar[r]^-{\VC\q}
 \ar[d]_-{\~\sg_a}
 &
 \VC\HS\E A
 \ar[r]^-{\r}
 \ar[d]_-{\-\sg_a}
 &
 \Vc\HS\E A
 \ar[dl]^-{\sg_a}
 \\
 \HC\E A
 \ar[r]^-{\q}
 &
 \HS\E A,
 &
 }
 $$
 where $\q$ and $\r$ are projections.
 We call $\sg_a$ the {\it canonical contraction}.

 \begin {claim} [11.1. Lemma.]
 Let
 $B\subseteq A$ be a subset.
 Then,
 for $a\in B$,
 the diagram
 $$
 \xymatrix {
 \Vc\HS\E B
 \ar[r]
 \ar[d]_-{\sg_a}
 &
 \Vc\HS\E A
 \ar[d]^-{\sg_a}
 \\
 \HS\E B
 \ar[r]
 &
 \HS\E A,
 }
 $$
 where the horizontal arrows are induced by the inclusion $B\to A$,
 is commutative.
 \qed
 \end {claim}

 If $A$ is finite,
 let
 \begin {equation} \label {theta}
 \theta_A:|\E A|\to\Delta A
 \end {equation}
 be the unbased map
 that sends,
 for each $a\in A$,
 the corresponding vertex $|a|$ of $|\E A|$ to
 the corresponding vertex $\`a\'$ of $\Delta A$
 and
 is affine on simplices.
 Hereafter,
 we put $\Delta\varnothing=\varnothing$.


 \subhead {Barycentric subdivision.}
 Let $\K$ be an
 (abstract simplicial)
 complex.
 We order the set of simplices of $\K$
 by reverse inclusion.
 Define the simplicial set $\b\K$ as
 the nerve of this partially ordered set.
 For a subcomplex $\L\subseteq\K$,
 we have $\b\L\subseteq\b\K$.

 There is a unique homeomorphism $|\b\K|\to|\K|$
 that
 sends
 the vertex of $|\b\K|$
 corresponding to a simplex $k$ of $\K$
 to
 the barycentre of the simplex $|k|\subseteq|\K|$
 and
 takes
 each simplex of $|\b\K|$
 to
 some simplex of $|\K|$
 affinely.
 Using it,
 we adopt that
 $|\b\K|=|\K|$.
 

 \subhead {Canonical retractions.}
 Given
 a complex $\K$
 and
 a subcomplex $\L\subseteq\K$,
 we
 have $\VC\b\L\subseteq\VC\b\K$
 and
 define the based morphism
 $$
 \rh^\K_\L:
 \VC\b\K
 \to
 \VC\b\L
 $$
 as the retraction
 that sends all vertices outside $\VC\b\L$ to the apex of
 $\VC\b\L$.
 We call $\rh^\K_\L$ the {\it canonical retraction}.

 \begin {claim} [11.2. Lemma.]
 For two subcomplexes $\L,\M\subseteq\K$,
 the diagram
 $$
 \xymatrix {
 \VC\b\L
 \ar[r]^-{\Inc}
 \ar[d]_-{\rh^\L_{\L\cap\M}}
 &
 \VC\b\K
 \ar[d]^-{\rh^\K_\M}
 \\
 \VC\b(\L\cap\M)
 \ar[r]^-{\Inc}
 &
 \VC\b\M
 }
 $$
 is commutative.
 \qed
 \end {claim}


 \head {\S~12. Canonical retractions in the cones $\VC\b\d E$
 and $\vC\Delta E$}


 Fix a nonempty finite set $E$.


 \subhead {The simplex $\d E$ and its subcomplexes.}
 Let the $\d E$ be the complex whose
 set of vertices is $E$
 and
 set of simplices is $\PP_\x(E)$.
 For $F\in\PP_\x(E)$,
 we have the subcomplex $\d F\subseteq\d E$.
 For $A\in\AA(E)$,
 introduce the subcomplex
 $$
 \d[A]
 =
 \coprod_{F\in A}
 \d F
 \subseteq
 \d E.
 $$
 For $A,B\in\AA(E)$,
 we have
 $$
 A\ge B
 \ \Rightarrow\
 \d[A]\supseteq\d[B]
 $$
 and
 $\d[A\mee B]=\d[A]\cap\d[B]$.
 Moreover,
 $\d[\{E\}]=\d E$.

 For $A,B\in\AA(E)$, $A\ge B$,
 we have
 the canonical retraction
 $$
 \rh^A_B
 =
 \rh^{\d[A]}_{\d[B]}:
 \VC\b\d[A]
 \to
 \VC\b\d[B].
 $$

 \begin {claim} [12.1. Corollary.]
 For two layouts $A,B\in\AA(E)$,
 the diagram
 $$
 \xymatrix {
 \VC\b\d[A]
 \ar[r]^-{\Inc}
 \ar[d]_-{\rh^A_{A\mee B}}
 &
 \VC\b\d E
 \ar[d]^-{\rh^{\{E\}}_B}
 \\
 \VC\b\d[A\mee B]
 \ar[r]^-{\Inc}
 &
 \VC\b\d[B]
 }
 $$
 is commutative.
 \end {claim}

 \begin {dem}
 Follows from Lemma~11.2.
 \qed
 \end {dem}


 \subhead {Geometric realization.}
 We adopt the obvious identification $|\d E|=\Delta E$.
 For $F\in\PP_\x(E)$,
 $|\d F|=\Delta F$ as subsets of $\Delta E$.
 For $A\in\AA(E)$,
 $|\d[A]|=\Delta[A]$ in the same sense.
 For $A,B\in\AA(E)$, $A\ge B$,
 we have
 $\Delta[A]\supseteq\Delta[B]$
 and
 the retraction $\rho^A_B$,
 $$
 \xymatrix @R1.4ex {
 \vC\Delta[A]
 \ar[r]^-{\rho^A_B}
 \ar@{=}[d]
 &
 \vC\Delta[B]
 \ar@{=}[d]
 \\
 |\VC\b\d[A]|
 \ar[r]^-{|\rh^A_B|}
 &
 |\VC\b\d[B]|.
 }
 $$
 We call $\rho^A_B$ the {\it canonical retraction},
 too.

 \begin {claim} [12.2. Corollary.]
 For two layouts $A,B\in\AA(E)$,
 the diagram
 $$
 \xymatrix {
 \vC\Delta[A]
 \ar[r]^-{\inc}
 \ar[d]_-{\rho^A_{A\mee B}}
 &
 \vC\Delta E
 \ar[d]^-{\rho^{\{E\}}_B}
 \\
 \vC\Delta[A\mee B]
 \ar[r]^-{\inc}
 &
 \vC\Delta[B]
 }
 $$
 is commutative.
 \end {claim}

 \begin {dem}
 Follows from Corollary~12.1.
 \qed
 \end {dem}


 \head {\S~13. The fissilizer $\Phi_E$ on
 $\<(Y^X_a)^{\vC\Delta E}\>$}


 Fix
 a space $Z$
 and
 a finite set $E$.
 Consider the presheaf $M:\PP_\x(E)\to\Mg$,
 $F\mapsto Z^{\vC\Delta F}$
 (with the obvious restriction functions).
 For $A\in\AA(E)$,
 we have
 $$
 \vC\Delta[A]
 =
 \vC
 \bigl(
 \coprod_{F\in A}
 \Delta F
 \bigr)
 =
 \bigbou_{F\in A}
 \vC\Delta F.
 $$
 We identify the presheaf $\_M:\AA(E)\to\Mg$
 (see \S~10)
 with the presheaf $A\mapsto Z^{\vC\Delta[A]}$
 by the chain of equalities/obvious identifications
 \begin {equation} \label {MM}
 \_M(A)
 =
 \prod_{F\in A}
 M(F)
 =
 \prod_{F\in A}
 Z^{\vC\Delta F}
 =
 Z^{
 \bigbou_{F\in A}
 \vC\Delta F
 }
 =
 Z^{\vC\Delta[A]}.
 \end {equation}
 In our case,
 the combining product \eqref{bigecol} coincides with
 the $\ZZ$-multilinear operation
 \begin {equation} \label {bigeBou}
 \bigeBou_{F\in A}:
 \prod_{F\in A}
 \<Z^{\vC\Delta F}\>
 \to
 \<Z^{\vC\Delta[A]}\>,
 \qquad
 \bigeBou_{F\in A}
 \`v_F\'
 =
 \`
 \bigBou_{F\in A}
 v_F
 \',
 \end {equation}
 (cf.\ also \eqref{eBou}).
 We will need the following formulas:
 \begin {equation} \label {epsilon-hash}
 \epsilon
 \bigl(
 \bigeBou_{F\in A}
 Q_F
 \bigr)
 =
 \prod_{F\in A}
 \epsilon(Q_F)
 \end {equation}
 and
 \begin {equation} \label {hash-restricted}
 \bigl(
 \bigeBou_{F\in A}
 Q_F
 \bigr)
 \big|_{\vC\Delta G}
 =
 \bigl(
 \prod_{F\in A\setminus\{G\}}
 \epsilon(Q_F)
 \bigr)
 Q_G,
 \qquad
 G\in A.
 \end {equation}

 For $A,B\in\AA(E)$, $A\ge B$,
 let $\lambda^B_A:\_M(B)\to\_M(A)$ be the function
 $$
 Z^{\rho^A_B}:
 Z^{\vC\Delta[B]}
 \to
 Z^{\vC\Delta[A]},
 $$
 where
 $\rho^A_B:\vC\Delta[A]\to\vC\Delta[B]$ is
 the canonical retraction.
 It follows from Corollary~12.2 that
 the functions $\lambda^B_A$ form an extender.
 By \S~10,
 we get the fissilizer
 \begin {equation} \label {Phi-Z}
 \Phi_E:
 \<Z^{\vC\Delta E}\>
 \to
 \<Z^{\vC\Delta E}\>.
 \end {equation}

 \begin {claim} [13.1. Corollary.]
 For any ensemble $Q\in\<Z^{\vC\Delta E}\>$,
 the ensemble $\Phi_E(Q)$ is fissile.
 \end {claim}

 \begin {dem}
 Follows from Lemma~10.1.
 \qed
 \end {dem}

 \subhead {$(X,r)$-almost fissile ensembles.}
 Let
 $X$
 and
 $Y$
 be spaces,
 $X$ compact Hausdorff,
 and
 $a:X\to Y$ be a map.
 Consider the space $(Y^X,a)$,
 which is $Y^X$ with $a$ as the basepoint.
 For a space $T$,
 we have the inclusion $\<(Y^X,a)^T\>\subseteq\<(Y^X)^{(T)}\>$.
 
 An ensemble $Q\in\<(Y^X,a)^{\vC\Delta E}\>$ is called
 {\it $(X,r)$-almost fissile\/} if,
 for any layout $A\in\AA(E)$,
 $$
 \bigeBou_{F\in A}
 Q|_{\vC\Delta F}
 -
 Q|_{\vC\Delta[A]}
 \in
 \<(Y^X)^{(\vC\Delta[A])}\>^{(r+1)}_X.
 $$

 \begin {claim} [13.2. Lemma.]
 Any affine ensemble $Q\in\<(Y^X,a)^{\vC\Delta E}\>$ is
 $(X,1)$-almost fissile.
 \end {claim}

 \begin {demo} [Proof.]
 Take $A\in\AA(E)$.
 Consider the quantity $D\in\<(Y^X,a)^{\vC\Delta[A]}\>$,
 $$
 D
 =
 \bigeBou_{F\in A}
 Q|_{\vC\Delta F}
 -
 Q|_{\vC\Delta[A]}.
 $$
 We should show that
 $D\in\<(Y^X)^{(\vC\Delta[A])}\>^{(2)}_X$.
 Consider the homomorphism
 $$
 \<\uc^X\>:
 \<(Y^X)^{(\vC\Delta[A])}\>
 \to
 \<Y^{\vC\Delta[A]\hsm X}\>.
 $$
 We should show that
 $\<\uc^X\>(D)\in\<Y^{\vC\Delta[A]\hsm X}\>^{(2)}$.
 Take $R\in\F_1(\vC\Delta[A]\hsm X)$.
 We check that
 $\<\uc^X\>(D)|_R=0$.
 We are in (at least) one of the two following cases.

 {\it Case 0:\/}
 $R=\{\0\}$.
 We have
 \begin {multline*}
 \epsilon(\<\uc^X\>(D))
 =
 \epsilon(D)
 =
 \hfill
 \text{(using \eqref{epsilon-hash})}
 \hfill
 =
 \prod_{F\in A}
 \epsilon(Q|_{\vC\Delta F})
 -
 \epsilon(Q|_{\vC\Delta[A]})
 =
 \qquad
 \\
 \qquad
 =
 \prod_{F\in A}
 \epsilon(Q)
 -
 \epsilon(Q)
 =
 \hfill
 \text{(since $\epsilon(Q)=1$)}
 \hfill
 =0,
 \end {multline*}
 which suffices in this case.
 
 {\it Case 1:\/}
 $R\subseteq\vC\Delta G\hsm X$
 for some $G\in A$.
 It suffices to check that
 $\<\uc^X\>(D)|_{\vC\Delta G\hsm X}=0$.
 We have the commutative diagram
 $$
 \xymatrix {
 \<(Y^X)^{(\vC\Delta[A])}\>
 \ar[rr]^-{\<\uc^X\>}
 \ar[d]_-{\?|_{\vC\Delta G}}
 &&
 \<Y^{\vC\Delta[A]\hsm X}\>
 \ar[d]^-{\?|_{\vC\Delta G\hsm X}}
 \\
 \<(Y^X)^{(\vC\Delta G)}\>
 \ar[rr]^-{\<\uc^X\>}
 &&
 \<Y^{\vC\Delta G\hsm X}\>.
 }
 $$
 Thus
 it suffices to check that
 $D|_{\vC\Delta G}=0$.
 We have
 \begin {multline*}
 D|_{\vC\Delta G}
 =
 \hfill
 \text{(using \eqref{hash-restricted})}
 \hfill
 =
 \bigl(
 \prod_{F\in A\setminus\{G\}}
 \epsilon(Q|_{\vC\Delta F})
 \bigr)
 Q|_{\vC\Delta G}
 -
 Q|_{\vC\Delta G}
 =
 \qquad
 \\
 \qquad
 \text{(since $\epsilon(Q|_{\vC\Delta F})=\epsilon(Q)=1$)}
 \hfill
 =0.
 \QED
 \end {multline*}
 \end {demo}

 We let $Z$ be a subspace of $(Y^X,a)$.
 For a space $T$,
 we have the inclusion $\<Z^T\>\subseteq\<(Y^X,a)^T\>$.
 We have the fissilizer \eqref{Phi-Z}.

 \begin {claim} [13.3. Corollary.]
 Let $Q\in\<Z^{\vC\Delta E}\>$ be an $(X,r)$-almost
 fissile ensemble.
 Then
 $$
 \Phi_E(Q)
 -
 Q
 \in
 \<(Y^X)^{(\vC\Delta E)}\>^{(r+1)}_X.
 $$
 \end {claim}

 \begin {demo} [Proof.]
 For $A\in\AA(E)$,
 put
 $$
 N(A)
 =
 \<Z^{\vC\Delta[A]}\>
 \cap
 \<(Y^X)^{(\vC\Delta[A])}\>^{(r+1)}_X
 \subseteq
 \<(Y^X)^{(\vC\Delta[A])}\>.
 $$
 We have
 \begin {multline*}
 N(A)
 \subseteq
 \<Z^{\vC\Delta[A]}\>
 =
 \hfill
 \text{(by \eqref{MM})}
 \hfill
 =
 \<\_M(A)\>.
 \end {multline*}
 By Lemma~3.1,
 this family is preserved by
 the restriction homomorphisms of the presheaf
 $A\mapsto\<\_M(A)\>$
 and
 the homomorphisms $\<\lambda^B_A\>$.
 Since $Q$ is $(X,r)$-almost fissile,
 it satisfies the hypothesis of Lemma~10.2.
 Thus
 $\Phi_E(Q)-Q\in N(\{E\})$.
 Clearly,
 $N(\{E\})\subseteq\<(Y^X)^{(\vC\Delta E)}\>^{(r+1)}_X$.
 \qed
 \end {demo}

 \subhead {A sufficient condition of strong $r$-similarity.}
 Suppose that
 $X$
 and
 $Y$
 are cellular.
 Let $a,b:X\to Y$ be maps.
 We
 let $a$ be the basepoint of the component $Y^X_a\subseteq Y^X$
 and
 set $Z=Y^X_a$.
 We say that
 $a$ is {\it firmly $r$-similar to $b$},
 $a\%\approxident rb$,
 if,
 for any nonempty finite set $E$,
 there is a fissile ensemble $R\in\<(Y^X_a)^{\vC\Delta E}\>$
 such that
 \begin {equation} \label {R}
 \`\Xi^{\Delta E}(b)\'
 -
 R|_{\Delta E}
 \in
 \<(Y^X)^{(\Delta E)}\>^{(r+1)}_X.
 \end {equation}

 \begin {claim} [13.4. Lemma.]
 Let $a,b:X\to Y$ be maps.
 Then
 $a\%\approxident rb$ implies $a\%\approx rb$.
 \end {claim}

 We do not know whether the converse holds.

 \begin {demo} [Proof.]
 Take a nonempty finite set $E$.
 We have a fissile ensemble $R\in\<(Y^X_a)^{\vC\Delta E}\>$
 satisfying \eqref{R}.
 We seek a fissile ensemble $S\in\<(Y^X)^{\Delta E}\>$
 such that
 \begin {equation} \label {S}
 \`\Xi^{\Delta E}(b)\'
 -
 S
 \in
 \<(Y^X)^{(\Delta E)}\>^{(r+1)}_X.
 \end {equation}
 Put $S= R|_{\Delta E}$.

 For a layout $A\in\AA(E)$,
 we have
 \begin {multline*}
 S|_{\Delta[A]}
 =
 R|_{\vC\Delta[A]}|_{\Delta[A]}
 =
 \hfill
 \text{(since $R$ is fissile)}
 \hfill
 =
 \bigl(
 \bigeBou_{F\in A}
 R|_{\vC\Delta F}
 \bigr)
 \big|_{\Delta[A]}
 =
 \qquad
 \\
 \qquad
 \text{(by the definitions of
 $\textstyle\bigeBou$
 (\eqref{bigeBou})
 and
 $\textstyle\bigeCop$
 (\eqref{bigeCop}))}
 \hfill
 =
 \bigeCop_{F\in A}
 R|_{\Delta F}
 =
 \bigeCop_{F\in A}
 S|_{\Delta F}.
 \end {multline*}
 Thus $S$ is fissile.

 The condition \eqref{S} is just the equality \eqref{R}.
 \qed
 \end {demo}

 \begin {claim} [13.5. Proposition.]
 Let $a,b:X\to Y$ be maps.
 Suppose that,
 for any nonempty finite set $E$,
 there is an $(X,r)$-almost fissile ensemble
 $Q_E\in\<(Y^X_a)^{\vC\Delta E}\>$ such that
 \begin {equation} \label {QE}
 \`\Xi^{\Delta E}(b)\'
 -
 Q_E|_{\Delta E}
 \in
 \<(Y^X)^{(\Delta E)}\>^{(r+1)}_X.
 \end {equation}
 Then
 $a\%\approx rb$
 and,
 moreover,
 $a\%\approxident rb$.
 \end {claim}

 \begin {demo} [Proof.]
 Take a nonempty finite set $E$.
 Put
 $Q=Q_E$
 and
 $R=\Phi_E(Q)\in\<(Y^X_a)^{\vC\Delta E}\>$.
 By Corollary~13.1,
 $R$ is fissile.
 By Corollary~13.3,
 $$
 R-Q
 \in
 \<(Y^X)^{(\vC\Delta E)}\>^{(r+1)}_X.
 $$
 By Lemma~3.1,
 $$
 R|_{\Delta E}
 -
 Q|_{\Delta E}
 \in
 \<(Y^X)^{(\Delta E)}\>^{(r+1)}_X.
 $$
 Using \eqref{QE},
 we get
 $$
 \`\Xi^{\Delta E}(b)\'
 -
 R|_{\Delta E}
 \in
 \<(Y^X)^{(\Delta E)}\>^{(r+1)}_X.
 $$
 Thus
 $a\%\approxident rb$.
 By Lemma~13.4,
 $a\%\approx rb$.
 \qed
 \end {demo}


 \head {\S~14. Strong $1$-similarity}


 Let
 $X$
 and
 $Y$
 be spaces,
 $X$ compact Hausdorff.

 \begin {claim} [14.1. Lemma.]
 Let $U$ be an unbased space.
 Then
 the homomorphism
 $$
 \<\Xi^U\>:
 \<Y^X\>
 \to
 \<(Y^X)^{(U)}\>
 $$
 takes $\<Y^X\>^{(s)}$ to $\<(Y^X)^{(U)}\>^{(s)}_X$.
 \end {claim}

 \begin {demo} [Proof.]
 Consider the projection
 $$
 p:
 U\hsm X
 \to
 X,
 \qquad
 u\shs x\mapsto x.
 $$
 We have the commutative diagram
 $$
 \xymatrix {
 \<Y^X\>
 \ar[rr]^-{\<\Xi^U\>}
 \ar[drr]_-{\<Y^p\>}
 &&
 \<(Y^X)^{(U)}\>
 \ar[d]^-{\<\uc^X\>}
 \\
 &&
 \<Y^{U\hsm X}\>.
 }
 $$

 The homomorphism $\<Y^p\>$ takes $\<Y^X\>^{(s)}$ to
 $\<Y^{U\hsm X}\>^{(s)}$
 by \cite[Lemma~2.1]{sim-1-cir}.
 Thus,
 by the diagram,
 $\<\Xi^U\>$ takes $\<Y^X\>^{(s)}$ to
 $\<\uc^X\>^{-1}(\<Y^{U\hsm X}\>^{(s)})$,
 which is $\<(Y^X)^{(U)}\>^{(r+1)}_X$
 by the definition of the latter.
 \qed
 \end {demo}

 Suppose that
 $X$
 and
 $Y$
 are cellular.

 \begin {claim} [14.2. Theorem.]
 Let $a,b:X\to Y$ be maps such that
 $a\%\sim1b$.
 Then
 $a\%\approx1b$.
 \end {claim}

 \begin {demo} [Proof.]
 We have an ensemble $A\in\<Y^X_a\>$,
 $$
 A=
 \sum_i
 u_i\`a_i\',
 $$
 such that
 $\`b\'-A\in\<Y^X\>^{(2)}$.
 For each $i$,
 choose a path $h_i:[0,1]\to Y^X_a$ from $a$ to $a_i$
 and
 consider the composition
 $$
 q_i:
 \vC\Delta E
 \xto{\mathrm{projection}}
 [0,1]
 \xto{h_i}
 Y^X_a.
 $$
 Consider the ensemble $Q\in\<(Y^X_a)^{\vC\Delta E}\>$,
 $$
 Q=
 \sum_i
 u_i\`q_i\'.
 $$
 We have
 \begin {multline*}
 \epsilon(Q)
 =
 \epsilon(A)
 =
 \hfill
 \text{(since $\`b\'-A\in\<Y^X\>^{(1)}$)}
 \hfill
 =
 \epsilon(\`b\')
 =
 1.
 \end {multline*}
 By Lemma~13.2,
 $Q$ is $(X,1)$-almost fissile.
 Clearly,
 $q_i|_{\Delta E}=\Xi^{\Delta E}(a_i)$.
 Thus
 $Q|_{\Delta E}=\<\Xi^{\Delta E}\>(A)$.
 We get
 $$
 \`\Xi^{\Delta E}(b)\'-Q|_{\Delta E}
 =
 \<\Xi^{\Delta E}\>(\`b\'-A)
 \in
 \<(Y^X)^{(\Delta E)}\>^{(2)},
 $$
 where $\in$ holds by Lemma~14.1.
 By Proposition~13.5,
 $a\%\approx1b$.
 \qed
 \end {demo}


 \head {\S~15. Two identities}


 Let
 $A$
 and
 $I$
 be finite sets.
 Let $\PP(I)$ be the set of subsets of $I$.
 Consider the set $\PP(I)^A$ of functions $k:A\to\PP(I)$.
 For $k\in\PP(I)^A$,
 put
 $$
 U(k)
 =
 \bigcup_{a\in A}
 k(a)
 \in
 \PP(I).
 $$
 Let $\RR(A,I)$ be the set of $k\in\PP(I)^A$ such that
 $U(k)=I$
 (covers).

 \begin {claim} [15.1. Lemma.]
 In the group $\<\PP(I)\>^{\otimes A}$,
 the equality holds
 $$
 \sum_{J\in\PP(I)}
 (-1)^{|I|-|J|}
 \bigotimes_{a\in A}
 \`J\'
 =
 \sum_{k\in\RR(A,I)}
 \
 \bigotimes_{a\in A}
 \
 \sum_{J\in\PP(k(a))}
 (-1)^{|k(a)|-|J|}
 \`J\'.
 $$
 \end {claim}

 \begin {demo} [Proof.]
 We have
 \begin {multline*}
 \sum_{J\in\PP(I)}
 (-1)^{|I|-|J|}
 \bigotimes_{a\in A}
 \
 \sum_{K\in\PP(J)}
 \`K\'
 =
 \\
 =
 \sum_{J\in\PP(I)}
 (-1)^{|I|-|J|}
 \sum_{\substack{
 k\in\PP(I)^A:
 \\
 U(k)\subseteq J
 }}
 \
 \bigotimes_{a\in A}
 \`k(a)\'
 =
 \\
 =
 \sum_{k\in\PP(I)^A}
 \bigl(
 \sum_{\substack{
 J\in\PP(I):
 \\
 J\supseteq U(k)
 }}
 (-1)^{|I|-|J|}
 \bigr)
 \bigotimes_{a\in A}
 \`k(a)\'
 \overset{(*)}
 =
 \sum_{k\in\RR(A,I)}
 \
 \bigotimes_{a\in A}
 \`k(a)\',
 \end {multline*}
 where $(*)$ holds because
 the inner sum on the left equals
 $1$
 if $U(k)=I$
 and
 $0$
 otherwise.
 The set $\PP(I)$ is partially ordered by inclusion.
 We have the isomorphism
 $$
 \nabla_{\PP(I)}^{-1}:
 \<\PP(I)\>
 \to
 \<\PP(I)\>
 $$
 (see \S~9),
 under which
 $$
 \sum_{K\in\PP(J)}
 \`K\'
 \mapsto
 \`J\',
 \qquad
 J\in\PP(I),
 $$
 and
 $$
 \`K\'
 \mapsto
 \sum_{J\in\PP(K)}
 (-1)^{|K|-|J|}
 \`J\',
 \qquad
 K\in\PP(I).
 $$
 Applying it to each factor of the summands in the
 left
 and
 right
 sides of the calculation,
 we get the required equality.
 \qed
 \end {demo}

 Put $\PP^\x(I)=\PP(I)\setminus\{I\}$.
 We adopt the inclusion $\PP^\x(I)^A\subseteq\PP(I)^A$.
 Let $\RR'(A,I)$ be the set of $k\in\PP^\x(I)^A$ such that
 $U(k)=I$.

 \begin {claim} [15.2. Lemma.]
 In the group $\<\PP^\x(I)\>^{\otimes A}$,
 the equality holds
 \begin {multline*}
 \bigotimes_{a\in A}
 \
 \sum_{J\in\PP^\x(I)}
 (-1)^{|I|-1-|J|}
 \`J\'
 -
 \sum_{J\in\PP^\x(I)}
 (-1)^{|I|-1-|J|}
 \bigotimes_{a\in A}
 \`J\'
 =
 \\
 =
 \sum_{k\in\RR'(A,I)}
 \
 \bigotimes_{a\in A}
 \
 \sum_{J\in\PP(k(a))}
 (-1)^{|k(a)|-|J|}
 \`J\'.
 \end {multline*}
 \end {claim}

 \begin {demo} [Proof.]
 We use the inclusion
 $\<\PP^\x(I)\>^{\otimes A}\subseteq\<\PP(I)\>^{\otimes A}$.
 Put
 $$
 T(k)
 =
 \bigotimes_{a\in A}
 \
 \sum_{J\in\PP(k(a))}
 (-1)^{|k(a)|-|J|}
 \`J\',
 \qquad
 k\in\PP(I)^A.
 $$
 We have
 \begin {multline} \label {P}
 \sum_{k\in\PP(I)^A}
 T(k)
 =
 \bigotimes_{a\in A}
 \
 \sum_{K\in\PP(I)}
 \
 \sum_{J\in\PP(K)}
 (-1)^{|K|-|J|}
 \`J\'
 =
 \\
 =
 \bigotimes_{a\in A}
 \
 \sum_{J\in\PP(I)}
 \bigl(
 \sum_{\substack{
 K\in\PP(I):
 \\
 K\supseteq J
 }}
 (-1)^{|K|-|J|}
 \bigr)
 \`J\'
 \overset{(*)}
 =
 \bigotimes_{a\in A}
 \`I\',
 \end {multline}
 where $(*)$ holds because
 the inner sum on the left equals
 $1$
 if $J=I$
 and
 $0$
 otherwise.
 We have also
 \begin {multline} \label {Px}
 \sum_{k\in\PP^\x(I)^A}
 T(k)
 =
 \bigotimes_{a\in A}
 \
 \sum_{K\in\PP^\x(I)}
 \
 \sum_{J\in\PP(K)}
 (-1)^{|K|-|J|}
 \`J\'
 =
 \\
 =
 \bigotimes_{a\in A}
 \
 \sum_{J\in\PP^\x(I)}
 \bigl(
 \sum_{\substack{
 K\in\PP^\x(I):
 \\
 K\supseteq J
 }}
 (-1)^{|K|-|J|}
 \bigr)
 \`J\'
 =
 \\
 =
 \bigotimes_{a\in A}
 \
 \sum_{J\in\PP^\x(I)}
 (-1)^{|I|-1-|J|}
 \`J\'.
 \end {multline}
 Note that
 $$
 \RR(A,I)\supseteq\RR'(A,I),
 \qquad
 \PP(I)^A\supseteq\PP^\x(I)^A,
 $$
 and
 $$
 \RR(A,I)\setminus\RR'(A,I)
 =
 \PP(I)^A\setminus\PP^\x(I)^A
 $$
 as subsets of $\PP(I)^A$.
 Thus
 $$
 \sum_{k\in\RR'(A,I)}
 T(k)
 =
 \sum_{k\in\RR(A,I)}
 T(k)
 -
 \sum_{k\in\PP(I)^A}
 T(k)
 +
 \sum_{k\in\PP^\x(I)^A}
 T(k)
 =
 $$
 (by
 Lemma~15.1
 and
 equalities
 \eqref{P}
 and
 \eqref{Px})
 \begin {multline*}
 =
 \sum_{J\in\PP(I)}
 (-1)^{|I|-|J|}
 \bigotimes_{a\in A}
 \`J\'
 - 
 \bigotimes_{a\in A}
 \`I\'
 +
 \bigotimes_{a\in A}
 \
 \sum_{J\in\PP^\x(I)}
 (-1)^{|I|-1-|J|}
 \`J\'
 =
 \\
 =
 -
 \sum_{J\in\PP^\x(I)}
 (-1)^{|I|-1-|J|}
 \bigotimes_{a\in A}
 \`J\'
 +
 \bigotimes_{a\in A}
 \
 \sum_{J\in\PP^\x(I)}
 (-1)^{|I|-1-|J|}
 \`J\',
 \end {multline*}
 as required.
 \qed
 \end {demo}


 \head {\S~16. Chained monoids}


 Let $P$ be a monoid.
 Then $\<P\>$ is its monoid ring.
 We call the monoid $P$ {\it chained\/} if
 $\<P\>$ is equipped with a chain of left ideals $\<P\>^{[s]}$,
 $$
 \<P\>
 =
 \<P\>^{[0]}
 \supseteq
 \<P\>^{[1]}
 \supseteq
 \dotso.
 $$

 Given a finite set $I$,
 we
 consider $\PP(I)$ as a monoid with respect to intersection
 and
 chain it by letting $\<\PP(I)\>^{[s]}$ be
 the subgroup generated by elements
 $$
 \omega_J
 =
 \sum_{K\in\PP(J)}
 (-1)^{|J|-|K|}
 \`K\',
 $$
 where $J\in\PP(I)$, $|J|\ge s$.


 \head {\S~17. The filtration $\<\Z^\T\>^{[s]}$}


 Let $P$ be a chained monoid.
 Let
 $\T$
 and
 $\Z$
 be based simplicial sets.
 Let $\Z^\T$ denote the set of based morphisms $\T\to\Z$.
 Let $P$ act on $\Z$
 (on the left;
 preserving the basepoint).
 For an element $p\in P$,
 let
 $p_{(\Z)}:\Z\to\Z$
 be its action.
 (We will use this notation for all actions.)
 The set $\Z^\T$ carries the induced action of $P$.
 Thus
 the abelian group $\<\Z^\T\>$ becomes a
 (left)
 module over $\<P\>$.
 We define a filtration
 $$
 \<\Z^\T\>
 =
 \<\Z^\T\>^{[0]}
 \supseteq
 \<\Z^\T\>^{[1]}
 \supseteq
 \dotso.
 $$
 Let
 $\T_j$, $j\in(n)$,
 be based simplicial sets
 and
 $$
 \f:
 \T
 \to
 \bigbou_{j\in(n)}
 \T_j
 $$
 be a based morphism.
 We have
 the $\ZZ$-multilinear operation
 $$
 \bigeBou_{j\in(n)}:
 \prod_{j\in(n)}
 \<\Z^{\T_j}\>
 \to
 \<
 \Z^{
 \bigbou_{j\in(n)}
 \T_j
 }\>,
 \qquad
 \bigeBou_{j\in(n)}
 \`\v_j\'
 =
 \`
 \bigBou_{j\in(n)}
 \v_j
 \',
 $$
 (combining product,
 cf.\
 \eqref{bigeCop},
 \eqref{eBou}
 and
 \eqref{bigeBou}),
 and
 the homomorphism
 $$
 \<\Z^\f\>:
 \<
 \Z^{
 \bigbou_{j\in(n)}
 \T_j
 }\>
 \to
 \<\Z^\T\>.
 $$
 Take ensembles $v_j\in\<P\>^{[s_j]}\<\Z^{\T_j}\>$, $j\in(n)$,
 and
 consider the ensemble $v\in\<\Z^\T\>$,
 \begin {equation} \label {block}
 v
 =
 \<\Z^\f\>
 \bigl(
 \bigeBou_{j\in(n)}
 v_j
 \bigr).
 \end {equation}
 We call $v$ a {\it block\/} of {\it rank\/} $s_1+\dotso+s_n$.
 We let $\<\Z^\T\>^{[s]}\subseteq\<\Z^\T\>$ be the subgroup
 generated by all blocks of rank at least $s$.
 One easily sees that
 it is a submodule.

 \begin {claim} [17.1. Lemma.]
 Let
 $\tT$ be a based simplicial set
 and
 $\k:\tT\to\T$ be a based simplicial morphism.
 Then
 the homomorphism
 $$
 \<\Z^\k\>:
 \<\Z^\T\>
 \to
 \<\Z^{\tT}\>
 $$
 takes $\<\Z^\T\>^{[s]}$ to $\<\Z^{\tT}\>^{[s]}$.
 \qed
 \end {claim}

 \begin {claim} [17.2. Lemma.]
 Let
 $\tZ$ be a based simplicial set
 with an action of $P$
 and
 $\h:\Z\to\tZ$ be a $P$-equivariant based simplicial morphism.
 Then
 the homomorphism
 $$
 \<\h^\T\>:
 \<\Z^\T\>
 \to
 \<\tZ^\T\>
 $$
 takes $\<\Z^\T\>^{[s]}$ to $\<\tZ^\T\>^{[s]}$.
 \qed
 \end {claim}

 The cone $\Vc\Z$ carries the induced action of $P$.
 We have the function
 $$
 \Vc^\T_\Z:
 \Z^\T
 \to
 (\Vc\Z)^{\Vc\T},
 \qquad
 \v\mapsto\Vc\v.
 $$

 \begin {claim} [17.3. Lemma.]
 The homomorphism
 $$
 \<\Vc^\T_\Z\>:
 \<\Z^\T\>
 \to
 \<(\Vc\Z)^{\Vc\T}\>
 $$
 takes $\<\Z^\T\>^{[s]}$ to $\<(\Vc\Z)^{\Vc\T}\>^{[s]}$.
 \end {claim}

 \begin {demo} [Proof.]
 It suffices to show that
 $\<\Vc^\T_\Z\>$ sends any block to a block of the same rank.
 Consider
 the block \eqref{block}.
 Since
 $v_j\in\<P\>^{[s_j]}\<\Z^{\T_j}\>$
 and
 the functions
 $$
 \Vc^{\T_j}_\Z:
 \Z^{\T_j}
 \to
 (\Vc\Z)^{\Vc\T_j}
 $$
 preserve the action of $P$,
 we have
 $$
 \<\Vc^{\T_j}_\Z\>(v_j)
 \in
 \<P\>^{[s_j]}\<(\Vc\Z)^{\Vc\T_j}\>.
 $$
 Let
 $$
 \Ins_k:
 \T_k
 \to
 \bigbou_{j\in(n)}
 \T_j
 $$
 be the canonical insertions.
 We have the commutative diagram
 $$
 \xymatrix {
 &&
 \bigbou\limits_{j\in(n)}
 \Vc\T_j
 \ar[d]^-{
 \e
 :=
 \bigBou\limits_{j\in(n)}
 \Vc\,\Ins_j
 }
 \\
 \Vc\T
 \ar[urr]^-{\g}
 \ar[rr]_-{\Vc\f}
 &&
 \Vc(
 \bigbou\limits_{j\in(n)}
 \T_j
 ),
 }
 $$
 where
 $\e$ is an isomorphism
 (since $\Vc$ preserves wedges)
 and
 $\g$ is the unique lift of $\Vc\f$.
 For arbitrary based morphisms $\v_j:\T_j\to\Z$,
 we have the commutative diagram with sendings
 $$
 \xymatrix {
 &
 (\Vc\Z)^{
 \bigbou_{j\in(n)}
 \Vc\T_j
 }
 \ar[dl]_-{(\Vc\Z)^\g}
 \\
 (\Vc\Z)^{\Vc\T}
 &
 (\Vc\Z)^{
 \Vc(
 \bigbou_{j\in(n)}
 \T_j
 )
 },
 \ar[l]^-{(\Vc\Z)^{\Vc\f}}
 \ar[u]_-{(\Vc\Z)^\e}
 }
 \qquad
 \xymatrix {
 &
 {\scriptstyle
 \bigBou\limits_{j\in(n)}
 \Vc\v_j
 }
 \ar@{|->}[dl]
 \\
 {\scriptstyle
 \Vc(
 \Z^\f(
 \bigBou\limits_{j\in(n)}
 \v_j
 )
 )
 }
 &
 {\scriptstyle
 \Vc(
 \bigBou\limits_{j\in(n)}
 \v_j
 ).
 }
 \ar@{|->}[l]
 \ar@{|->}[u]
 }
 $$
 Thus
 we have
 the commutative diagram
 $$
 \xymatrix {
 &&
 \<(\Vc\Z)^{
 \bigbou_{j\in(n)}
 \Vc\T_j
 }\>
 \ar[dll]_-{\<(\Vc\Z)^\g\>}
 \\
 \<(\Vc\Z)^{\Vc\T}\>
 &&
 \<(\Vc\Z)^{
 \Vc
 (
 \bigbou_{j\in(n)}
 \T_j
 )
 }\>
 \ar[ll]^-{\<(\Vc\Z)^{\Vc\f}\>}
 \ar[u]_-{(\<\Vc\Z)^\e\>}
 }
 $$
 and the sendings
 $$
 \xymatrix {
 &&
 {\scriptstyle
 \bigeBou\limits_{j\in(n)}
 \<\Vc^{\T_j}_\Z\>(v_j)
 }
 \ar@{|->}[dll]
 \\
 {\scriptstyle
 \<\Vc^\T_\Z\>(
 \<\Z^\f\>(
 \bigeBou\limits_{j\in(n)}
 v_j
 )
 )
 }
 &&
 {\scriptstyle
 \<\Vc^{
 \bigbou_{j\in(n)}
 \T_j
 }_\Z\>(
 \bigeBou\limits_{j\in(n)}
 v_j
 )
 }
 \ar@{|->}[ll]
 \ar@{|->}[u]
 }
 $$
 for our
 (and
 arbitrary)
 ensembles $v_j$.
 We get
 $$
 \<\Vc^\T_\Z\>(v)
 =
 \<\Vc^\T_\Z\>
 \bigl(
 \<\Z^\f\>
 \bigl(
 \bigeBou\limits_{j\in(n)}
 v_j
 \bigr)
 \bigr)
 =
 \<(\Vc\Z)^\g\>
 \bigl(
 \bigeBou\limits_{j\in(n)}
 \<\Vc^{\T_j}_\Z\>(v_j)
 \bigr),
 $$
 as promised.
 \qed
 \end {demo}

 \begin {claim} [17.4. Lemma.]
 Let
 $\T_i$,
 $i\in(m)$,
 be based simplicial sets
 and
 $v_i\in\<\Z^{\T_i}\>^{[s_i]}$ be ensembles.
 Then
 $$
 \bigeBou_{i\in(m)}
 v_i
 \in
 \<\Z^{
 \bigbou_{i\in(m)}
 \T_i
 }\>^{[s_1+\dotso+s_m]}.
 \QED
 $$
 \end {claim}


 \subhead {Fissile and almost fissile ensembles.}
 Let $E$ be a nonempy finite set.
 For a layout $A\in\AA(E)$,
 we have
 $$
 \VC\b\d[A]
 =
 \bigvee_{F\in A}
 \VC\b\d F
 \subseteq
 \VC\b\d E.
 $$
 An ensemble $q\in\<\Z^{\VC\b\d E}\>$ is called
 {\it fissile\/} if,
 for any $A\in\AA(E)$,
 $$
 q|_{\VC\b\d[A]}
 =
 \bigeBou_{F\in A}
 q|_{\VC\b\d F}
 $$
 in $\<\Z^{\VC\b\d[A]}\>$
 (cf.\ \S\S\ 2, 10).
 It is called {\it $r$-almost fissile\/} if,
 for any $A\in\AA(E)$,
 $$
 \bigeBou_{F\in A}
 q|_{\VC\b\d F}
 -
 q|_{\VC\b\d[A]}
 \in
 \<\Z^{\VC\b\d[A]}\>^{[r+1]}
 $$
 (cf.\ \S~13).


 \head {\S~18. The wedge $\W(I)$}


 Fix a finite set $I$.
 Consider the based simplicial set
 $$
 \W(I)
 =
 \bigbou_{J\in\PP(I)}
 \HS\E(I\setminus J).
 $$
 Let
 $$
 \Ins_J:
 \HS\E(I\setminus J)
 \to
 \W(I)
 $$
 be the canonical insertions.
 The {\it lead\/} vertex
 $$
 \top_{\W(I)}
 =
 (\Ins_I)_{[0]}(
 1_{\HS\E\varnothing}
 )
 \in
 \W(I)_{[0]}
 $$
 is isolated.
 $\W(I)$ has the based simplicial subsets
 $$
 \W^\x(I)
 =
 \bigbou_{J\in\PP^\x(I)}
 \HS\E(I\setminus J)
 $$
 and
 $$
 \W^L(I)
 =
 \bigbou_{J\in\PP(L)}
 \HS\E(I\setminus J),
 \qquad
 L\in\PP^\x(I).
 $$

 For $J,K\in\PP(I)$, $J\supseteq K$,
 let
 $$
 \tg^J_K:
 \HS\E(I\setminus J)
 \to
 \HS\E(I\setminus K)
 $$
 be the morphism induced by the inclusion
 $I\setminus J\to I\setminus K$.

 Let the monoid $\PP(I)$ act on $\W(I)$ by the rule
 $$
 \xymatrix {
 \HS\E(I\setminus J)
 \ar[rr]^-{\tg^J_{K\cap J}}
 \ar[d]_-{\Ins_J}
 &&
 \HS\E(I\setminus(K\cap J))
 \ar[d]^-{\Ins_{K\cap J}}
 \\
 \W(I)
 \ar[rr]^-{K_{(\W(I))}}
 &&
 \W(I),
 }
 $$
 $K\in\PP(I)$.
 The simplicial subsets
 $\W^\x(I)$
 and
 $\W^L(I)$
 are $\PP(I)$-invariant.

 For
 $L\in\PP^\x(I)$
 and
 $i\in I\setminus L$,
 we define a retraction $\sg^L_i$ by the commutative diagram
 $$
 \xymatrix {
 \Vc\HS\E(I\setminus J)
 \ar[r]^-{\Vc\,\Ins^L_J}
 \ar[d]_-{\sg_i}
 &
 \Vc\W^L(I)
 \ar[d]^-{\sg^L_i}
 \\
 \HS\E(I\setminus J)
 \ar[r]^-{\Ins^L_J}
 &
 \W^L(I),
 }
 $$
 where
 $\Ins^L_J$ are the canonical insertions
 and
 $\sg_i$ are the canonical contractions (see \S~11).
 We call $\sg^L_i$ the canonical contraction,
 too.
 It follows from Lemma~11.1 that
 $\sg^L_i$ is $\PP(I)$-equivariant.

 Given
 a based simplicial set $\T$,
 introduce the {\it filling\/} function
 $$
 \chi^L_i:
 \W^L(I)^\T
 \to
 \W^L(I)^{\Vc\T},
 \qquad
 \v
 \mapsto
 (
 \Vc\T
 \xto{\Vc\v}
 \Vc\W^L(I)
 \xto{\sg^L_i}
 \W^L(I)
 ).
 $$
 Since $\sg^L_i$ is a retraction,
 \begin {equation} \label {fill}
 \chi^L_i(\v)|_\T=\v.
 \end {equation}


 \head {\S~19. The module $\<\W(I)^{\VC\b\d E}\>$}


 Fix a finite set $I$.
 We consider the $\<\PP(I)\>$-modules $\<\W(I)^\T\>$
 for a number of based simplicial sets $\T$.
 For a $\PP(I)$-invariant based simplicial subset
 $\Z\subseteq\W(I)$,
 the subgroup $\<\Z^\T\>\subseteq\<\W(I)^\T\>$ is a
 $\<\PP(I)\>$-submodule.
 If $\Z\subseteq\tZ$ for two such subsets,
 then
 $\<\Z^\T\>^{[s]}\subseteq\<\tZ^\T\>^{[s]}$
 by Lemma~17.2.

 \begin {claim} [19.1. Lemma.]
 For
 $L\in\PP^\x(I)$,
 $i\in I\setminus L$,
 and
 a based simplicial set $\T$,
 the filling homomorphism
 $$
 \<\chi^L_i\>:
 \<\W^L(I)^\T\>
 \to
 \<\W^L(I)^{\Vc\T}\>
 $$
 takes
 $\<\W^L(I)^\T\>^{[s]}$
 to
 $\<\W^L(I)^{\Vc\T}\>^{[s]}$.
 \end {claim}

 \begin {demo} [Proof.]
 By construction of $\chi^L_i$,
 we have the decomposition
 $$
 \<\chi^L_i\>:
 \<\W^L(I)^\T\>
 \xto{\<\Vc^\T_{\W^L(I)}\>}
 \<(\Vc\W^L(I))^{\Vc\T}\>
 \xto{\<(\sg^L_i)^{\Vc\T}\>}
 \<\W^L(I)^{\Vc\T}\>.
 $$
 By Lemma~17.3,
 $\<\Vc^\T_{\W^L(I)}\>$ takes
 $\<\W^L(I)^\T\>^{[s]}$
 to
 $\<(\Vc\W^L(I))^{\Vc\T}\>^{[s]}$.
 Since $\sg^L_i$ is $\PP(I)$-equivariant,
 $\<(\sg^L_i)^{\Vc\T}\>$ takes
 the latter
 to
 $\<\W^L(I)^{\Vc\T}\>^{[s]}$
 by Lemma~17.2.
 \qed
 \end {demo}

 Fix a nonempty finite set $E$.
 For
 $F\in\PP_\x(E)$
 and
 $J\in\PP(I)$,
 introduce the based morphism
 $$
 \xg^F_J:
 (\b\d F)_+
 \to
 \W(I)
 $$
 that takes $\b\d F$ to the vertex
 $(\Ins_J)_{[0]}(1_{\HS\E(I\setminus J)})$.

 \begin {claim} [19.2. Lemma.]
 For
 $F\in\PP_\x(E)$
 and
 $J\in\PP(I)$,
 $$
 \sum_{K\in\PP(J)}
 (-1)^{|J|-|K|}
 \`\xg^F_K\'
 \in
 \<\W^J(I)^{(\b\d F)_+}\>^{[|J|]}.
 $$
 \end {claim}

 \begin {demo} [Proof.]
 Since
 $$
 \xg^F_K
 =
 K_{(\W^J(I)^{(\b\d F)_+})}(\xg^F_J),
 $$
 the ensemble in question
 equals $\omega_J\`\xg^F_J\'$
 and
 thus
 belongs to $\<\PP(I)\>^{[|J|]}\<\W^J(I)^{(\b\d F)_+}\>$,
 which is contained in $\<\W^J(I)^{(\b\d F)_+}\>^{[|J|]}$
 by the definition of the latter.
 \qed
 \end {demo}

 \begin {claim} [19.3. Lemma.]
 There exist fissile ensembles
 $$
 p_J\in\<\W^\x(I)^{\VC\b\d E}\>,
 \qquad
 J\in\PP^\x(I),
 $$
 satisfying the following conditions
 for each $J\in\PP^\x(I)$:

 \smallskip
 \noindent
 $(1)$
 one has
 $$
 p_J|_{(\b\d E)_+}=\`\xg^E_J\'
 $$
 in $\<\W^\x(I)^{(\b\d E)_+}\>$;

 \smallskip
 \noindent
 $(2)$
 one has
 $$
 \sum_{K\in\PP(J)}
 (-1)^{|J|-|K|}
 p_K
 \in
 \<\W^\x(I)^{\VC\b\d E}\>^{[|J|]}.
 $$
 \end {claim}

 \begin {demo} [Proof.]
 We will construct ensembles
 $$
 p^F_J\in\<\W^J(I)^{\VC\b\d F}\>,
 \qquad
 (F,J)\in\PP_\x(I)\cro\PP^\x(I),
 $$
 satisfying the following conditions
 $(0^F_J)$,
 $(1^F_J)$,
 and
 $(2^F_J)$
 for each pair $(F,J)\in\PP_\x(I)\cro\PP^\x(I)$:

 \smallskip
 \noindent
 $(0^F_J)$
 one has
 $$
 p^F_J|_{\VC\b\d[B]}
 =
 \bigeBou_{G\in B}
 p^G_J
 $$
 in $\<\W^J(I)^{\VC\b\d[B]}\>$
 for all $B\in\AA(F)$;

 \smallskip
 \noindent
 $(1^F_J)$
 one has
 $$
 p^F_J|_{(\b\d F)_+}=\`\xg^F_J\'
 $$
 in $\<\W^J(I)^{(\b\d F)_+}\>$;

 \smallskip
 \noindent
 $(2^F_J)$
 one has
 $$
 \sum_{K\in\PP(J)}
 (-1)^{|J|-|K|}
 p^F_K
 \in
 \<\W^J(I)^{\VC\b\d F}\>^{[|J|]}.
 $$

 \smallskip
 Note that
 $(0^F_J)$ implies
 $$
 p^F_J|_{\VC\b\d G}
 =
 p^G_J
 $$
 for $G\in\PP_\x(F)$.
 Thus
 $(0^F_J)$ will yield
 $$
 p^F_J|_{\VC\b\d[B]}
 =
 \bigeBou_{G\in B}
 p^F_J|_{\VC\b\d G}
 $$
 for all $B\in\AA(F)$,
 which means that
 $p^F_J$ is fissile.
 Thus
 it will remain to put $p_J=p^E_J$.

 Induction on $(F,J)\in\PP_\x(E)\cro\PP^\x(I)$.
 Take a pair $(F,J)$.
 We assume that
 $p^G_K$ are defined
 and
 the conditions $(0^G_K)$--$(2^G_K)$ are satisfied
 for
 $$
 (G,K)\in\PP_\x(F)\cro\PP(J)\setminus\{(F,J)\}.
 $$
 We
 construct $p^F_J$
 and
 check the conditions $(0^F_J)$--$(2^F_J)$.

 For $B\in\AA(F)$,
 put
 $$
 U(B)=\<\W^J(I)^{\VC\b\d[B]}\>^{[|J|]}.
 $$
 For $B,C\in\AA(F)$, $B\ge C$,
 we have,
 by Lemma~17.1,
 the restriction homomorphism
 $$
 \?|_{\VC\b\d[C]}:
 U(B)
 \to
 U(C).
 $$
 Thus we have a presheaf
 $$
 U:\AA(F)\to\Ab.
 $$
 By Lemma~17.1,
 the canonical retractions
 $$
 \rh^B_C:
 \VC\b\d[B]
 \to
 \VC\b\d[C]
 $$
 induce homomorphisms
 $$
 \lambda^C_B
 =
 \<\W^J(I)^{\rh^B_C}\>|_{ U(C)\to U(B)}:
 U(C)
 \to
 U(B),
 $$
 which form an extender for $U$,
 as follows from Corollary~12.1.
 For $B\in\AA^\x(F)=\AA(F)\setminus\{\{F\}\}$,
 introduce the ensemble $u_B\in\<\W^J(I)^{\VC\b\d[B]}\>$,
 $$
 u_B
 =
 \sum_{K\in\PP(J)}
 (-1)^{|J|-|K|}
 \bigeBou_{G\in B}
 p^G_K.
 $$
 By Lemma~15.1,
 $$
 u_B
 =
 \sum_{l\in\RR(B,J)}
 \
 \bigeBou_{G\in B}
 \
 \sum_{K\in\PP(l(G))}
 (-1)^{|l(G)|-|K|}
 p^G_K.
 $$
 By $(2^G_{l(G)})$,
 the inner sum belongs to
 $\<\W^J(I)^{\VC\b\d G}\>^{[|l(G)|]}$.
 Using
 Lemma~17.4
 and
 the inequality
 $$
 \sum_{G\in B}
 |l(G)|
 \ge
 |J|,
 $$
 we get that
 the combining product
 and
 thus
 $u_B$
 belong to $\<\W^J(I)^{\VC\b\d[B]}\>^{[|J|]}$.
 We have got
 $u_B\in U(B)$.
 For
 $B,C\in\AA^\x(F)$, $B\ge C$,
 and
 $K\in\PP(J)$,
 we have
 \begin {multline*}
 \bigl(
 \bigeBou_{G\in B}
 p^G_K
 \bigr)
 \big|_{\VC\b\d[C]}
 =
 \hfill
 \text{(by naturality of $\textstyle\bigeBou$)}
 \hfill
 =
 \bigeBou_{G\in B}
 p^G_K|_{\VC\b\d[C\mee\{G\}]}
 =
 \qquad
 \\
 \qquad
 \text{(by $(0^G_K)$)}
 \hfill
 =
 \bigeBou_{G\in B}
 \
 \bigeBou_{H\in C\mee\{G\}}
 p^H_K
 =
 \bigeBou_{H\in C}
 p^H_K.
 \end {multline*}
 It follows that $u_B|_{\VC\b\d[C]}=u_C$,
 that is,
 $$
 (u_B)_{B\in\AA^\x(F)}
 \in
 \lim_{B\in\AA^\x(F)}
 U(B).
 $$
 By Lemma~9.2,
 there exists an ensemble
 \begin {equation} \label {u}
 u
 \in
 U(\{F\})
 =
 \<\W^J(I)^{\VC\b\d F}\>^{[|J|]}
 \end {equation}
 such that
 $$
 u|_{\VC\b\d[B]}=u_B,
 \qquad
 B\in\AA^\x(F).
 $$

 Consider the ensembles $q,r\in\<\W^J(I)^{\VC\b\d F}\>$,
 $$
 q
 =
 \sum_{K\in\PP^\x(J)}
 (-1)^{|J|-1-|K|}
 p^F_K,
 \qquad
 r=q+u.
 $$

 For $B\in\AA^\x(F)$,
 we have
 \begin {multline} \label {qB}
 q|_{\VC\b\d[B]}
 =
 \sum_{K\in\PP^\x(J)}
 (-1)^{|J|-1-|K|}
 p^F_K|_{\VC\b\d[B]}
 =
 \hfill
 \text{(by $(0^F_K)$)}
 \qquad
 \\
 =
 \sum_{K\in\PP^\x(J)}
 (-1)^{|J|-1-|K|}
 \bigeBou_{G\in B}
 p^G_K
 \end {multline}
 and
 \begin {multline} \label {rB}
 r|_{\VC\b\d[B]}
 =
 q|_{\VC\b\d[B]}
 +
 u|_{\VC\b\d[B]}
 =
 q|_{\VC\b\d[B]}
 +
 u_B
 =
 \hfill
 \text{(by \eqref{qB})}
 \qquad
 \\
 =
 \sum_{K\in\PP^\x(J)}
 (-1)^{|J|-1-|K|}
 \bigeBou_{G\in B}
 p^G_K
 +
 \sum_{K\in\PP(J)}
 (-1)^{|J|-|K|}
 \bigeBou_{G\in B}
 p^G_K
 =
 \bigeBou_{G\in B}
 p^G_J.
 \end {multline}

 We have
 \begin {multline*}
 r|_{(\b\d F)_+}
 +
 \sum_{K\in\PP^\x(J)}
 (-1)^{|J|-|K|}
 \`\xg^F_K\'
 =
 \hfill
 \text{(by $(1^F_K)$)}
 \qquad
 \\
 =
 r|_{(\b\d F)_+}
 +
 \sum_{K\in\PP^\x(J)}
 (-1)^{|J|-|K|}
 p^F_K|_{(\b\d F)_+}
 =
 r|_{(\b\d F)_+}
 -
 q|_{(\b\d F)_+}
 =
 \\
 \qquad
 =
 u|_{(\b\d F)_+}
 \in
 \hfill
 \text{(by Lemma~17.1)}
 \hfill
 \in
 \<\W^J(I)^{(\b\d F)_+}\>^{[|J|]}.
 \end {multline*}
 From
 this
 and
 Lemma~19.2,
 \begin {equation} \label {r-base}
 \`\xg^F_J\'
 -
 r|_{(\b\d F)_+}
 \in
 \<\W^J(I)^{(\b\d F)_+}\>^{[|J|]}.
 \end {equation}

 Choose $i\in I\setminus J$.
 We have the filling homomorphism
 $$
 \<\chi^J_i\>:
 \<\W^J(I)^{(\b\d F)_+}\>
 \to
 \<\W^J(I)^{\VC\b\d F}\>.
 $$
 Put
 $$
 p^F_J
 =
 r
 +
 \<\chi^J_i\>(
 \`\xg^F_J\'
 -
 r|_{(\b\d F)_+}
 ).
 $$

 {\it Check of $(0^F_J)$.}
 For $B=\{F\}$,
 the condition is satisfied trivially.
 Take $B\in\AA^\x(F)$.
 We have
 \begin {multline} \label {rB-base}
 r|_{(\b\d[B])_+}
 =
 \hfill
 \text{(by \eqref{rB} and naturality of $\textstyle\bigeBou$)}
 \hfill
 =
 \bigeBou_{G\in B}
 p^G_J|_{(\b\d G)_+}
 =
 \qquad
 \\
 \qquad
 \text{(by $(1^G_J)$)}
 \hfill
 =
 \bigeBou_{G\in B}
 \`\xg^G_J\'
 =
 \`\xg^F_J\'|_{(\b\d[B])_+}.
 \end {multline}
 By construction of $\chi^J_i$,
 we have the commutative diagram
 $$
 \xymatrix {
 \W^J(I)^{(\b\d F)_+}
 \ar[r]^-{\chi^J_i}
 \ar[d]_-{\?|_{(\b\d[B])_+}}
 &
 \W^J(I)^{\VC\b\d F}
 \ar[d]^-{\?|_{\VC\b\d[B]}}
 \\
 \W^J(I)^{(\b\d[B])_+}
 \ar[r]^-{\chi^J_i}
 &
 \W^J(I)^{\VC\b\d[B]}.
 }
 $$
 We get
 \begin {multline*}
 p^F_J|_{\VC\b\d[B]}
 =
 r|_{\VC\b\d[B]}
 +
 \<\chi^J_i\>(
 \`\xg^F_J\'
 -
 r|_{(\b\d F)_+}
 )|_{\VC\b\d[B]}
 =
 \\
 \qquad
 \text{(by the diagram)}
 \hfill
 =
 r|_{\VC\b\d[B]}
 +
 \<\chi^J_i\>(
 \`\xg^F_J\'|_{(\b\d[B])_+}
 -
 r|_{(\b\d[B])_+}
 )
 =
 \qquad
 \\
 \qquad
 \text{(by \eqref{rB-base})}
 \hfill
 =
 r|_{\VC\b\d[B]}
 =
 \hfill
 \text{(by \eqref{rB})}
 \hfill
 =
  \bigeBou_{G\in B}
 p^G_J.
 \end {multline*}

 {\it Check of $(1^F_J)$.}
 We have
 \begin {multline*}
 p^F_J|_{(\b\d F)_+}
 -
 r|_{(\b\d F)_+}
 =
 \<\chi^J_i\>(
 \`\xg^F_J\'
 -
 r|_{(\b\d F)_+}
 )|_{(\b\d F)_+}
 =
 \\
 \qquad
 \text{(by \eqref{fill})}
 \hfill
 =
 \`\xg^F_J\'
 -
 r|_{(\b\d F)_+}.
 \end {multline*}
 Thus
 $p^F_J|_{(\b\d F)_+}=\`\xg^F_J\'$.

 {\it Check of $(2^F_J)$.}
 It follows from \eqref{r-base}
 by Lemma~19.1,
 that
 \begin {equation} \label {chi}
 \<\chi^J_i\>(
 \`\xg^F_J\'
 -
 r|_{(\b\d F)_+}
 )
 \in
 \<\W^J(I)^{\VC\b\d F}\>^{[|J|]}.
 \end {equation}
 We have
 \begin {multline*}
 \sum_{K\in\PP(J)}
 (-1)^{|J|-|K|}
 p^F_K
 =
 p^F_J-q
 =
 r
 +
 \<\chi^J_i\>(
 \`\xg^F_J\'
 -
 r|_{(\b\d F)_+}
 )
 -
 q
 =
 \\
 =
 u
 +
 \<\chi^J_i\>(
 \`\xg^F_J\'
 -
 r|_{(\b\d F)_+}
 )
 \in
 \<\W^J(I)^{\VC\b\d F}\>^{[|J|]},
 \end {multline*}
 where $\in$ follows from \eqref{u} and \eqref{chi}.
 \qed
 \end {demo}

 \begin {claim} [19.4. Corollary.]
 There exists an $(|I|-1)$-almost fissile ensemble
 $q\in\<\W^\x(I)^{\VC\b\d E}\>$ such that
 $$
 \`\xg^E_I\'
 -
 q|_{(\b\d E)_+}
 \in
 \<\W(I)^{(\b\d E)_+}\>^{[|I|]}.
 $$
 \end {claim}

 \begin {demo} [Proof.]
 Lemma~19.3 gives fissile ensembles
 $p_J\in\<\W^\x(I)^{\VC\b\d E}\>$
 satisfying the conditions
 $(1)$
 and
 $(2)$
 thereof.
 Put
 \begin {equation*}
 q
 =
 \sum_{J\in\PP^\x(I)}
 (-1)^{|I|-1-|J|}
 p_J.
 \end {equation*}

 Check that
 $q$ is $(|I|-1)$-almost fissile.
 Take $A\in\AA(E)$.
 We have
 \begin {multline*}
 \bigeBou_{F\in A}
 q|_{\VC\b\d F}
 -
 q|_{\VC\b\d[A]}
 =
 \bigeBou_{F\in A}
 \
 \sum_{J\in\PP^\x(I)}
 (-1)^{|I|-1-|J|}
 p_J|_{\VC\b\d F}
 -
 \\
 \qquad
 -
 \sum_{J\in\PP^\x(I)}
 (-1)^{|I|-1-|J|}
 p_J|_{\VC\b\d[A]}
 =
 \hfill
 \text{(since $p_J$ are fissile)}
 \qquad
 \\
 =
 \bigeBou_{F\in A}
 \
 \sum_{J\in\PP^\x(I)}
 (-1)^{|I|-1-|J|}
 p_J|_{\VC\b\d F}
 -
 \sum_{J\in\PP^\x(I)}
 (-1)^{|I|-1-|J|}
 \bigeBou_{F\in A}
 p_J|_{\VC\b\d F}
 =
 \\
 \qquad
 \text{(by Lemma~15.2)}
 \hfill
 =
 \sum_{k\in\RR'(A,I)}
 \
 \bigeBou_{F\in A}
 \
 \sum_{J\in\PP(k(F))}
 (-1)^{|k(F)|-|J|}
 p_J|_{\VC\b\d F}
 =
 \qquad
 \\
 =
 \sum_{k\in\RR'(A,I)}
 \
 \bigeBou_{F\in A}
 \bigl(
 \sum_{J\in\PP(k(F))}
 (-1)^{|k(F)|-|J|}
 p_J
 \bigr)
 \big|_{\VC\b\d F}.
 \end {multline*}
 By condition $(2)$,
 the inner sum of the last expression belongs to
 $\<\W^\x(I)^{\VC\b\d E}\>^{[|k(F)|]}$.
 By Lemma~17.1,
 its restriction to $\VC\b\d F$ belongs to
 $\<\W^\x(I)^{\VC\b\d F}\>^{[|k(F)|]}$.
 Using
 Lemma~17.4
 and
 the inequality
 $$
 \sum_{F\in A}
 |k(F)|
 \ge
 |I|,
 $$
 we get that
 the combining product
 and
 thus
 the whole expression
 belong to $\<\W^\x(I)^{\VC\b\d[A]}\>^{[|I|]}$,
 as required.

 We have
 \begin {multline*}
 \`\xg^E_I\'
 -
 q|_{(\b\d E)_+}
 =
 \`\xg^E_I\'
 -
 \sum_{J\in\PP^\x(I)}
 (-1)^{|I|-1-|J|}
 p_J|_{(\b\d E)_+}
 =
 \\
 \qquad
 \text{(by condition $(1)$)}
 \hfill
 =
 \`\xg^E_I\'
 -
 \sum_{J\in\PP^\x(I)}
 (-1)^{|I|-1-|J|}
 \`\xg^E_J\'
 =
 \qquad
 \\
 \qquad
 =
 \sum_{J\in\PP(I)}
 (-1)^{|I|-|J|}
 \`\xg^E_J\'
 \in
 \hfill
 \text{(by Lemma~19.2)}
 \hfill
 \in
 \<\W(I)^{(\b\d E)_+}\>^{[|I|]}.
 \end {multline*}
 \qed
 \end {demo}


 \head {\S~20. The filtration $\<(Y^X)^T\>^{[s]}$}


 We give a topological version of the definition of \S~17.
 Let
 $T$
 and
 $Z$
 be spaces.
 Let a chained monoid $P$ act on $Z$
 (preserving the basepoint).
 The set $Z^T$ carries the induced action of $P$.
 Thus
 the abelian group $\<Z^T\>$ becomes a module over $\<P\>$.
 We define a filtration $\<Z^T\>^{[s]}$.
 Let
 $T_j$, $j\in(n)$,
 be spaces
 and
 $$
 f:
 T
 \to
 \bigbou_{j\in(n)}
 T_j
 $$
 be a map.
 Take ensembles $V_j\in\<P\>^{[s_j]}\<Z^{T_j}\>$,
 $j\in(n)$,
 and
 consider the ensemble $V\in\<Z^T\>$,
 \begin {equation} \label {V}
 V
 =
 \<Z^f\>
 \bigl(
 \bigeBou_{j\in(n)}
 V_j
 \bigr).
 \end {equation}
 We call $V$ a {\it block\/} of {\it rank\/} $s_1+\dotso+s_n$.
 We let $\<Z^T\>^{[s]}\subseteq\<Z^T\>$ be the subgroup
 generated by all blocks of rank at least $s$.
 One easily sees that
 it is a submodule.

 \begin {claim} [20.1. Lemma.]
 Let
 $\~Z$ be a space
 with an action of $P$
 and
 $h:Z\to\~Z$ be a $P$-equivariant map.
 Then
 the homomorphism
 $$
 \<h^T\>:
 \<Z^T\>
 \to
 \<\~Z^T\>
 $$
 takes $\<Z^T\>^{[s]}$ to $\<\~Z^T\>^{[s]}$.
 \qed
 \end {claim}

 \begin {claim} [20.2. Lemma.]
 Let
 $\T$
 and
 $\Z$
 be based simplicial sets.
 Let
 $P$ act on $\Z$
 and
 thus
 on $|\Z|$.
 Consider
 the geometric realization function
 $$
 \gamma:
 \Z^\T
 \to
 |\Z|^{|\T|},
 \qquad
 \v\mapsto|\v|,
 $$
 and
 the homomorphism
 $$
 \<\gamma\>:
 \<\Z^\T\>
 \to
 \<|\Z|^{|\T|}\>.
 $$
 Then
 $\<\gamma\>$ takes $\<\Z^\T\>^{[s]}$ to
 $\<|\Z|^{|\T|}\>^{[s]}$.
 \qed
 \end {claim}


 \subhead {The case $Z=Y^X$.}
 Let
 $I$ be a finite set
 and
 $Y$ be a space
 with an action of the chained monoid $P=\PP(I)$.
 We suppose that
 the action is {\it special\/}:
 $$
 Y
 =
 \bigcup_{i\in I}
 \Fix\{i\}_{(Y)}.
 $$

 \begin {claim} [20.3. Lemma.]
 Let $T$ be a space.
 The set $Y^T$ carries the induced action of $\PP(I)$.
 Then,
 in the $\<\PP(I)\>$-module $\<Y^T\>$,
 the inclusion holds
 $$
 \<\PP(I)\>^{[s]}\<Y^T\>
 \subseteq
 \<Y^T\>^{(s)}.
 $$
 \end {claim}

 \begin {demo} [Proof.]
 Take
 a map $u\in Y^T$
 and
 a subset $J\in\PP(I)$, $|J|\ge s$.
 The ensembles of the form $\omega_J\`u\'$ generate the subgroup
 $\<\PP(I)\>^{[s]}\<Y^T\>$.
 Thus
 we should show that
 $\omega_J\`u\'\in\<Y^T\>^{(s)}$.
 Take a subspace $R\in\F_{s-1}(T)$.
 We should check that
 $\omega_J\`u\'|_R=0$
 in $\<Y^R\>$.
 Since the action is special,
 for each $t\in T$,
 there is $i_t\in I$ such that
 $u(t)\in\Fix\{i_t\}_{(Y)}$.
 Consider the subset
 $$
 K
 =
 \{\,i_t\mid t\in R\setminus\{\0\}\}
 \in
 \PP(I).
 $$
 Clearly,
 $|K|<s$.
 For $t\in R\setminus\{\0\}$,
 we have
 \begin {multline*}
 K_{(Y)}(u(t))
 =
 K_{(Y)}(\{i_t\}_{(Y)}(u(t)))
 =
 \\
 =
 (K\cap\{i_t\})_{(Y)}(u(t))
 =
 \{i_t\}_{(Y)}(u(t))
 =
 u(t).
 \end {multline*}
 Thus
 $K_{(Y)}\circ u\|=Ru$.
 Thus
 $\`K\'\`u\'\|=R\`u\'$
 in $\<Y^R\>$.
 Since $|K|<s\le|J|$,
 we have $K\nsupseteq J$.
 It follows that
 $\omega_J\`K\'=0$
 in $\<\PP(I)\>$.
 We get
 $$
 \omega_J\`u\'
 \|=R
 \omega_J\`K\'\`u\'
 =
 0.
 \QED
 $$
 \end {demo}

 Let $X$ be a compact Hausdorff space.
 Consider the space $Z=Y^X$.
 It carries the induced action of $\PP(I)$.

 \begin {claim} [20.4. Lemma.]
 Let $T$ be a space.
 Then
 $$
 \<(Y^X)^T\>^{[s]}
 \subseteq
 \<(Y^X)^T\>^{(s)}_X.
 $$
 \end {claim}
 (See \eqref{YXT} for the filtration on the right.)

 \begin {demo} [Proof.]
 Take a block $V\in\<(Y^X)^T\>$ of rank at least $s$.
 We should show that
 $V\in\<(Y^X)^T\>^{(s)}_X$.
 Consider the isomorphism
 $$
 \<(Y^X)^T\>
 \xto{\<\huc^X\>}
 \<Y^{T\sma X}\>.
 $$
 By Lemma~3.3,
 we should show that
 $\<\huc^X\>(V)\in\<Y^{T\sma X}\>^{(s)}$.
 We have the equality \eqref{V}
 for some
 spaces $T_j$,
 map $f$
 and
 ensembles $V_j\in\<\PP(I)\>^{[s_j]}\<(Y^X)^{T_j}\>$,
 where $s_1+\dots s_n\ge s$.
 Since the function
 $$
 \huc^X:
 (Y^X)^{T_j}
 \to
 Y^{T_j\sma X}
 $$
 is $\PP(I)$-equivariant,
 $\<\huc^X\>(V_j)\in\<\PP(I)\>^{[s_j]}\<Y^{T_j\sma X}\>$.
 By Lemma~20.3,
 $$
 \<\huc^X\>(V_j)\in\<Y^{T_j\sma X}\>^{(s_j)}.
 $$
 Consider the commutative diagram
 $$
 \xymatrix @C1.4ex @R1.4ex {
 {\scriptstyle
 (V_j)_{j\in(n)}
 }
 \ar@{|->}[rrrrr]
 \ar@{|->}[dddd]
 &&&&&
 {\scriptstyle
 (\<\huc^X\>(V_j))_{j\in(n)}
 }
 \ar@{|->}[ddd]
 \\
 &
 \prod\limits_{j\in(n)}
 \<(Y^X)^{T_j}\>
 \ar[rrr]^-{
 \prod\limits_{j\in(n)}
 \<\huc^X\>
 }
 \ar[ddd]_-{
 \bigeBou\limits_{j\in(n)}
 }
 &&&
 \prod\limits_{j\in(n)}
 \<Y^{T_j\sma X}\>
 \ar[dd]^-{
 \bigeBou\limits_{j\in(n)}
 }
 &
 \\ \\
 &&&&
 \<Y^{
 \bigbou_{j\in(n)}
 T_j\sma X
 }\>
 \ar@{=}[d]
 &
 {\scriptstyle
 \bigeBou\limits_{j\in(n)}
 \<\huc^X\>(V_j)
 }
 \ar@{=}[d]
 \\
 {\scriptstyle
 \bigeBou\limits_{j\in(n)}
 V_j
 }
 \ar@{|->}[ddd]
 &
 \<(Y^X)^{
 \bigbou_{j\in(n)}
 T_j
 }\>
 \ar[rrr]^-{\<\huc^X\>}
 \ar[dd]_-{\<(Y^X)^f\>}
 &&&
 \<Y^{
 (
 \bigbou_{j\in(n)}
 T_j
 )\sma X
 }\>
 \ar[dd]^-{\<Y^{f\sma\id_X}\>}
 &
 {\scriptstyle
 W
 }
 \ar@{|->}[ddd]^-{(*)}
 \\ \\
 &
 \<(Y^X)^T\>
 \ar[rrr]^-{\<\huc^X\>}
 &&&
 \<Y^{T\sma X}\>.
 &
 \\
 {\scriptstyle
 V
 }
 \ar@{|->}[rrrrr]
 &&&&&
 {\scriptstyle
 \<\huc^X\>(V)
 }
 }
 $$
 Here the ensemble $W$ is defined by the equality shown.
 All the sendings are obvious
 except $(*)$,
 which follows by commutativity of the diagram.
 By \cite[Lemma~3.1]{sim-1-cir},
 $$
 \bigeBou_{j\in(n)}
 \<\huc^X\>(V_j)
 \in
 \<Y^{
 \bigbou_{j\in(n)}
 T_j\sma X
 }\>^{(s)}.
 $$
 Equivalently,
 $$
 W
 \in
 \<Y^{
 (
 \bigbou_{j\in(n)}
 T_j
 )
 \sma X
 }\>^{(s)}.
 $$
 By \cite[Lemma~2.1]{sim-1-cir},
 $\<\huc^X\>(V)\in\<Y^{T\sma X}\>^{(s)}$,
 as was to be shown.
 \qed
 \end {demo}


 \head {\S~21. The wedge $V(I)$ and a $\PP(I)$-equivariant map
 $h:V(I)\to Z$}


 Let $I$ be a finite set.
 We give a topological version of $\W(I)$.
 Consider the space
 $$
 V(I)
 =
 \bigbou_{J\in\PP(I)}
 \bS\Delta(I\setminus J).
 $$
 Let
 $$
 \ins_J:
 \bS\Delta(I\setminus J)
 \to
 V(I)
 $$
 be the canonical insertions.
 $V(I)$ consists of
 the isolated {\it lead\/} point
 $$
 \top_{V(I)}
 =
 \ins_I(1_{\bS\Delta\varnothing})
 $$
 and
 the subspace
 $$
 V^\x(I)
 =
 \bigbou_{J\in\PP^\x(I)}
 \bS\Delta(I\setminus J),
 $$
 which is
 contractible.

 For $J,K\in\PP(I)$, $J\supseteq K$,
 let
 $$
 \tau^J_K:
 \bS\Delta(I\setminus J)
 \to
 \bS\Delta(I\setminus K)
 $$
 be the map induced by the inclusion
 $I\setminus J\to I\setminus K$.

 Let the monoid $\PP(I)$ act on $V(I)$ by the rule
 $$
 \xymatrix {
 \bS\Delta(I\setminus J)
 \ar[rr]^-{\tau^J_{K\cap J}}
 \ar[d]_-{\ins_J}
 &&
 \bS\Delta(I\setminus(K\cap J))
 \ar[d]^-{\ins_{K\cap J}}
 \\
 V(I)
 \ar[rr]^-{K_{(V(I))}}
 &&
 V(I),
 }
 $$
 $K\in\PP(I)$.
 The subspace $V^\x(I)$ is $\PP(I)$-invariant.

 For $J\in\PP(I)$,
 we have the map
 $$
 e_J:
 |\HS\E(I\setminus J)|
 =
 \bS|\E(I\setminus J)|
 \xto{\bS\theta_{I\setminus J}}
 \bS\Delta(I\setminus J)
 $$
 (see \eqref{theta} for $\theta_{I\setminus J}$).
 These $e_J$ form the map
 \begin {equation} \label {e}
 e
 =
 \bigbou_{J\in\PP(I)}
 e_J:
 |\W(I)|
 \to
 V(I).
 \end {equation}
 It
 is $\PP(I)$-equivariant,
 sends the point $|\top_{\W(I)}|$ to $\top_{V(I)}$,
 and
 takes the subspace $|\W^\x(I)|$ to $V^\x(I)$.

 \begin {claim} [21.1. Lemma.]
 Let $Z$ be a space
 with an action of $\PP(I)$.
 Suppose that
 the basepoint path component $Z_\0\subseteq Z$ is weakly
 contractible.
 Let $\top_Z\in Z$ be a point such that
 $$
 K_{(Z)}(\top_Z)
 \in
 Z_\0
 $$
 for all $K\in\PP^\x(I)$.
 Then
 there exists a $\PP(I)$-equivariant map $h:V(I)\to Z$ such
 that
 $h(\top_{V(I)})=\top_Z$.
 \end {claim}

 \begin {demo} [Proof.]
 We
 crop $Z$
 and
 assume that
 $Z=Z_\0\cup\{\top_Z\}$.
 We will construct maps
 $$
 h^J:
 \bS\Delta(I\setminus J)
 \to
 Z,
 \qquad
 J\in\PP(I),
 $$
 satisfying the following conditions
 $(\top)$
 and
 $(*^K_J)$
 for $J,K\in\PP(I)$, $J\subseteq K$:

 \smallskip
 \noindent
 $(\top)$
 one has $h^I(1_{\bS\Delta\varnothing})=\top_Z$;

 \smallskip
 \noindent
 $(*^K_J)$
 the diagram
 $$
 \xymatrix {
 \bS\Delta(I\setminus K)
 \ar[r]^-{h^K}
 \ar[d]_-{\tau^K_J}
 &
 Z
 \ar[d]^-{J_{(Z)}}
 \\
 \bS\Delta(I\setminus J)
 \ar[r]^-{h^J}
 &
 Z
 }
 $$
 is commutative.

 \smallskip
 Note that
 the condition $(*^J_J)$
 is the equality $J_{(Z)}\circ h^J=h^J$.

 Induction on $J\in\PP(I)$.
 We define the map $h^I$ by the condition $(\top)$.
 The condition $(*^I_I)$ is satisfied trivially.
 Take $J\in\PP^\x(I)$.
 We assume that
 the maps $h^K$ are defined
 for $K\supsetneq J$
 and
 the conditions $(*^L_K)$ are satisfied
 for $L\supseteq K\supsetneq J$.
 We
 construct $h^J$
 and check $(*^K_J)$
 for $K\supseteq J$.

 For $K\supsetneq J$,
 put
 $$
 B_K
 =
 \Im(
 \bS\Delta(I\setminus K)
 \xto{\tau^K_J}
 \bS\Delta(I\setminus J)
 ).
 $$
 Since $\tau^K_J$ is an embedding,
 there is a map $f^K:B_K\to Z_\0$ such that
 $$
 f^K(\tau^K_J(t))
 =
 J_{(Z)}(h^K(t)),
 \qquad
 t\in\bS\Delta(I\setminus K),
 $$
 (we use here that
 $\Im J_{(Z)}\subseteq Z_\0$).
 We show that
 $$
 f^K
 \|={B_K\cap B_L}
 f^L
 $$
 for $K,L\supsetneq J$.
 Take $s\in B_K\cap B_L$.
 Since $B_K\cap B_L=B_{K\cup L}$,
 we have $s=\tau^{K\cup L}_J(t)$
 for some $t\in\bS\Delta(I\setminus(K\cup L))$.
 We have the commutative diagram
 $$
 \xymatrix {
 &
 {\scriptstyle
 t
 }
 \ar@{|->}[dl]
 &
 \bS\Delta(I\setminus(K\cup L))
 \ar[r]^-{h^{K\cup L}}
 \ar[dl]_-{\tau^{K\cup L}_J}
 \ar[d]^-{\tau^{{K\cup L}}_K}
 &
 Z
 \ar[d]^-{K_{(Z)}}
 \\
 {\scriptstyle
 s
 }
 &
 \bS\Delta(I\setminus J)
 &
 \bS\Delta(I\setminus K)
 \ar[l]_-{\tau^K_J}
 \ar[r]^-{h^K}
 &
 Z
 }
 $$
 (the square is commutative by $(*^{K\cup L}_K)$).
 Using the diagram,
 we get
 \begin {multline*}
 f^K(s)
 =
 f^K(\tau^{K\cup L}_J(t))
 =
 f^K(\tau^K_J(\tau^{K\cup L}_K(t)))
 =
 J_{(Z)}(h^K(\tau^{K\cup L}_K(t)))
 =
 \\
 =
 J_{(Z)}(K_{(Z)}(h^{K\cup L}(t)))
 =
 (J\cap K)_{(Z)}(h^{K\cup L}(t))
 =
 J_{(Z)}(h^{K\cup L}(t)).
 \end {multline*}
 Similarly,
 $f^L(s)=J_{(Z)}(h^{K\cup L}(t))$.
 Thus
 $f^K(s)=f^L(s)$,
 as promised.

 We have
 $$
 \bigcup_{K\supsetneq J}
 B_K
 =
 \bS\bd\Delta(I\setminus J)
 \subseteq
 \bS\Delta(I\setminus J),
 $$
 where $\bd\Delta(I\setminus J)$ denotes
 the boundary of the simplex $\Delta(I\setminus J)$.
 Since $B_K$ are closed,
 there is a map
 $$
 f:
 \bS\bd\Delta(I\setminus J)
 \to
 Z_\0
 $$
 such that
 $f|_{B_K}=f^K$
 for all $K\supsetneq J$.
 Since
 $\bS\bd\Delta(I\setminus J)$ is the boundary
 of the ball $\bS\Delta(I\setminus J)$
 and
 $Z_\0$ is weakly contractible,
 $f$ extends to a map
 $$
 g:
 \bS\Delta(I\setminus J)
 \to
 Z_\0.
 $$
 We put
 $$
 h^J(s)
 =
 J_{(Z)}(g(s)),
 \qquad
 s\in\bS\Delta(I\setminus J).
 $$
 
 Clearly,
 $J_{(Z)}\circ h^J=h^J$,
 which is the condition $(*^J_J)$.
 We check the condition $(*^K_J)$
 for $K\supsetneq J$.
 For $t\in\bS\Delta(I\setminus K)$,
 we have
 \begin {multline*}
  h^J(\tau^K_J(t))
 =
 J_{(Z)}(g(\tau^K_J(t)))
 =
 J_{(Z)}(f(\tau^K_J(t)))
 =
 J_{(Z)}(f^K(\tau^K_J(t)))
 =
 \\
 =
 J_{(Z)}(J_{(Z)}(h^K(t)))
 =
 J_{(Z)}(h^K(t)),
 \end {multline*}
 as required.

 We union all the $h^J$ into the desired $h$:
 $$
 h
 =
 \bigBou_{J\in\PP(I)}
 h^J.
 $$
 Since $\top_{V(I)}=\ins_I(1_{\bS\Delta\varnothing})$,
 we have
 \begin {multline*}
 h(\top_{V(I)})
 =
 h^I(1_{\bS\Delta\varnothing})
 =
 \hfill
 \text{(by $(\top)$)}
 \hfill
 =
 \top_Z.
 \end {multline*}
 To show that
 $h$ is $\PP(I)$-equivariant,
 we should check that,
 for $K,J\in\PP(I)$,
 the diagram
 $$
 \xymatrix {
 \bS\Delta(I\setminus J)
 \ar[r]^-{h^J}
 \ar[d]_-{\tau^J_{K\cap J}}
 &
 Z
 \ar[d]^-{K_{(Z)}}
 \\
 \bS\Delta(I\setminus(K\cap J))
 \ar[r]^-{h^{K\cap J}}
 &
 Z
 }
 $$
 is commutative.
 Indeed,
 \begin {multline*}
 K_{(Z)}\circ h^J
 =
 \hfill
 \text{(by $(*^J_J)$)}
 \hfill
 =
 K_{(Z)}\circ J_{(Z)}\circ h^J
 =
 \qquad
 \\
 \qquad
 =
 (K\cap J)_{(Z)}\circ h^J
 =
 \hfill
 \text{(by $(*^J_{K\cap J})$)}
 \hfill
 =
 h^{K\cap J}\circ\tau^J_{K\cap J}.
 \QED
 \end {multline*}
 \end {demo}


 \head {\S~22. The realization
 $\Upsilon_h:\W(I)^\T\to(Y^X)^{|\T|}$}


 Let
 $X$
 and
 $Y$
 be spaces,
 $X$ compact Hausdorff.
 Let
 $I$ be a finite set
 and
 $Y$ carry a special action of the monoid $\PP(I)$.
 Let $h:V(I)\to Y^X$ be a $\PP(I)$-equivariant map.
 Let $\T$ be a based simplicial set.
 Introduce the function
 $$
 \Upsilon_h
 =
 \Upsilon^\T_h:
 \W(I)^\T
 \to
 (Y^X)^{|\T|},
 \qquad
 \v
 \mapsto
 (
 |\T|
 \xto{|v|}
 |\W(I)|
 \xto{e}
 V(I)
 \xto{h}
 Y^X
 ),
 $$
 (see \eqref{e} for $e$),
 the {\it realization}.

 \begin {claim} [22.1. Lemma.]
 The function $\Upsilon_h$ takes
 $\W^\x(I)^\T$
 to
 $(Y^X_\0)^{|\T|}$.
 \end {claim}

 \begin {demo} [Proof.]
 The map $e$ takes $|\W^\x(I)|$ to $V^\x(I)$.
 Since $V^\x(I)$ is path connected,
 $h$ takes it to $Y^X_\0$.
 \qed
 \end {demo}

 Consider the homomorphism
 $$
 \<\Upsilon_h\>:
 \<\W(I)^\T\>
 \to
 \<(Y^X)^{|\T|}\>.
 $$

 \begin {claim} [22.2. Lemma.]
 The homomorphism $\<\Upsilon_h\>$ takes
 $\<\W(I)^\T\>^{[s]}$
 to
 $\<(Y^X)^{|\T|}\>^{(s)}_X$.
 \end {claim}

 \begin {demo} [Proof.]
 We have the decomposition
 $$
 \<\Upsilon_h\>:
 \<\W(I)^\T\>
 \xto{\<\gamma\>}
 \<|\W(I)|^{|\T|}\>
 \xto{\<(h\circ e)^{|\T|}\>}
 \<(Y^X)^{|\T|}\>,
 $$
 where $\gamma:\W(I)^\T\to|\W(I)|^{|\T|}$ is the geometric
 realization function.
 By Lemma~20.2,
 $\<\gamma\>$ takes
 $\<\W(I)^\T\>^{[s]}$
 to
 $\<|\W(I)|^{|\T|}\>^{[s]}$.
 By Lemma~20.1,
 $\<(h\circ e)^{|\T|}\>$ takes
 the latter
 to
 $\<(Y^X)^{|\T|}\>^{[s]}$,
 which is contained in
 $\<(Y^X)^{|\T|}\>^{(s)}_X$
 by Lemma~20.4.
 \qed
 \end {demo}

 \begin {claim} [22.3. Lemma.]
 Let
 $E$ be a nonempty finite set
 and
 $q\in\<\W(I)^{\VC\b\d E}\>$ be an $r$-almost fissile ensemble.
 Then
 the ensemble
 $\<\Upsilon_h\>(q)\in\<(Y^X)^{\vC\Delta E}\>$ is
 $(X,r)$-almost fissile.
 \end {claim}

 \begin {demo} [Proof.]
 Take $A\in\AA(E)$.
 The diagram
 $$
 \xymatrix {
 \prod\limits_{F\in A}
 \W(I)^{\VC\b\d E}
 \ar[rr]^-{
 \prod\limits_{F\in A}
 \Upsilon_h
 }
 \ar[d]_-{
 \prod\limits_{F\in A}
 \?|_{\VC\b\d F}
 }
 &&
 \prod\limits_{F\in A}
 (Y^X)^{\vC\Delta E}
 \ar[d]^-{
 \prod\limits_{F\in A}
 \?|_{\vC\Delta F}
 }
 \\
 \prod\limits_{F\in A}
 \W(I)^{\VC\b\d F}
 \ar[rr]^-{
 \prod\limits_{F\in A}
 \Upsilon_h
 }
 \ar[d]_-{
 \bigBou\limits_{F\in A}
 }
 &&
 \prod\limits_{F\in A}
 (Y^X)^{\vC\Delta F}
 \ar[d]^-{
 \bigBou\limits_{F\in A}
 }
 \\
 \W(I)^{\VC\b\d[A]}
 \ar[rr]^-{
 \Upsilon_h
 }
 &&
 (Y^X)^{\vC\Delta[A]}
 \\
 \W(I)^{\VC\b\d E}
 \ar[u]^-{
 \?|_{\VC\b\d[A]}
 }
 \ar[rr]^-{
 \Upsilon_h
 }
 &&
 (Y^X)^{\vC\Delta E}
 \ar[u]_-{
 \?|_{\vC\Delta[A]}
 }
 }
 $$
 is commutative because
 $\Upsilon^\T_h$ is natural with respect to $\T$.
 Thus the diagram
 $$
 \xymatrix {
 \prod\limits_{F\in A}
 \<\W(I)^{\VC\b\d E}\>
 \ar[rr]^-{
 \prod\limits_{F\in A}
 \<\Upsilon_h\>
 }
 \ar[d]_-{
 \prod\limits_{F\in A}
 \?|_{\VC\b\d F}
 }
 &&
 \prod\limits_{F\in A}
 \<(Y^X)^{\vC\Delta E}\>
 \ar[d]^-{
 \prod\limits_{F\in A}
 \?|_{\vC\Delta F}
 }
 \\
 \prod\limits_{F\in A}
 \<\W(I)^{\VC\b\d F}\>
 \ar[rr]^-{
 \prod\limits_{F\in A}
 \<\Upsilon_h\>
 }
 \ar[d]_-{
 \bigeBou\limits_{F\in A}
 }
 &&
 \prod\limits_{F\in A}
 \<(Y^X)^{\vC\Delta F}\>
 \ar[d]^-{
 \bigeBou\limits_{F\in A}
 }
 \\
 \<\W(I)^{\VC\b\d[A]}\>
 \ar[rr]^-{
 \<\Upsilon_h\>
 }
 &&
 \<(Y^X)^{\vC\Delta[A]}\>
 \\
 \<\W(I)^{\VC\b\d E}\>
 \ar[u]^-{
 \?|_{\VC\b\d[A]}
 }
 \ar[rr]^-{
 \<\Upsilon_h\>
 }
 &&
 \<(Y^X)^{\vC\Delta E}\>
 \ar[u]_-{
 \?|_{\vC\Delta[A]}
 }
 }
 $$
 is also commutative.
 In it,
 we have
 $$
 \xymatrix @R1.4ex {
 {\scriptstyle
 (q)_{F\in A}
 }
 \ar@{|->}[rr]
 \ar@{|->}[dd]
 &&
 {\scriptstyle
 (Q)_{F\in A}
 }
 \ar@{|->}[dd]
 \\ \\
 {\scriptstyle
 (q|_{\VC\b\d F})_{F\in A}
 }
 \ar@{|->}[dd]
 &&
 {\scriptstyle
 (Q|_{\vC\Delta F})_{F\in A}
 }
 \ar@{|->}[dd]
 \\ \\
 {\scriptstyle
 \bigeBou_{F\in A}
 q|_{\VC\b\d F}
 }
 \ar@{|->}[rr]^-{(1)}
 &&
 {\scriptstyle
 \bigeBou\limits_{F\in A}
 Q|_{\vC\Delta F}
 }
 \\
 {\scriptstyle
 q|_{\VC\b\d[A]}
 }
 \ar@{|->}[rr]^-{(2)}
 &&
 {\scriptstyle
 Q|_{\vC\Delta[A]}
 }
 \\ \\
 {\scriptstyle
 q
 }
 \ar@{|->}[uu]
 \ar@{|->}[rr]
 &&
 {\scriptstyle
 Q,
 }
 \ar@{|->}[uu]
 }
 $$
 where $Q=\<\Upsilon_h\>(q)$.
 All the sendings are obvious
 except
 $(1)$
 and
 $(2)$,
 which follow by commutativity of the diagram.
 Since $q$ is $r$-almost fissile,
 $$
 \bigeBou_{F\in A}
 q|_{\VC\b\d F}
 -
 q|_{\VC\b\d[A]}
 \in
 \<\W(I)^{\VC\b\d[A]}\>^{[r+1]}.
 $$
 By Lemma~22.2,
 $$
 \bigeBou_{F\in A}
 Q|_{\vC\Delta F}
 -
 Q|_{\vC\Delta[A]}
 \in
 \<(Y^X)^{\vC\Delta[A]}\>^{(r+1)}_X.
 $$
 Thus $Q$ is $(X,r)$-almost fissile.
 \qed
 \end {demo}


 \head {\S~23. Brunnian loops in a wedge of circles}


 Fix a finite set $I$ of cardinality $s$.
 Put
 $B(I)=I_+\sma S^1$
 (a wedge of $s$ circles).
 Let the monoid $\PP(I)$ act on the space $I_+$ by putting
 $$
 J_{(I_+)}(i)=
 \begin{cases}
 i
 &
 \text{if $i\in J$,}
 \\
 \0
 &
 \text{otherwise}
 \end{cases}
 $$
 for
 $i\in I_+$,
 $J\in\PP(I)$.
 This action induces one on $B(I)$.
 A map $v:S^1\to B(I)$
 (a loop)
 is called {\it Brunnian} if
 the composition
 $$
 S^1
 \xto{v}
 B(I)
 \xto{J_{(B(I))}}
 B(I)
 $$
 is null-homotopic
 for all $J\in\PP^\x(I)$.

 \begin {claim} [23.1. Lemma.]
 Let $w:S^1\to B(I)$ be a Brunnian loop.
 Then $\0\%\approx{s-1}w$.
 \end {claim}

 \begin {demo} [Proof.]
 Take a finite set $E$.
 Consider the loop space $B(I)^{S^1}$.
 It carries the induced action of the monoid $\PP(I)$.
 The path component $B(I)^{S^1}_\0$ is weakly contractible.
 Since $w$ is Brunnian,
 $J_{(B(I)^{S^1})}(w)$
 ($=J_{(B(I))}\circ w)$
 $\in B(I)^{S^1}_\0$
 for all $J\in\PP^\x(I)$.
 Therefore,
 Lemma~21.1 yields
 a $\PP(I)$-equivariant map $h:V(I)\to B(I)^{S^1}$ such that
 $h(\top_{V(I)})=w$.
 Consider the commutative diagram
 $$
 \xymatrix @C1.4ex {
 {\scriptstyle
 q
 }
 \ar@{|->}@/^3ex/[rrrrr]
 &
 \<\W(I)^{\VC\b\d E}\>
 \ar[rrr]^-{\<\Upsilon_h\>}
 \ar[d]_-{\?|_{(\b\d E)_+}}
 &&&
 \<(B(I)^{S^1})^{\vC\Delta E}\>
 \ar[d]^-{\?|_{(\Delta E)_+}}
 &
 {\scriptstyle
 Q
 }
 \\
 {\scriptstyle
 \`\xg^E_I\'-q|_{(\b\d E)_+}
 }
 \ar@{|->}@/_3ex/[rrrrr]
 &
 \<\W(I)^{(\b\d E)_+}\>
 \ar[rrr]^-{\<\Upsilon_h\>}
 &&&
 \<(B(I)^{S^1})^{(\Delta E)_+}\>
 \ar[d]^-{\?|_{\Delta E}}
 &
 {\scriptstyle
 \`\Upsilon_h(\xg^E_I)\'-Q|_{(\Delta E)_+}
 }
 \ar@{|->}[d]
 \\
 &
 &&&
 \<(B(I)^{S^1})^{(\Delta E)}\>,
 &
 {\scriptstyle
 \`\Xi^{\Delta E}(w)\'-Q|_{\Delta E}
 }
 }
 $$
 where $\Upsilon_h$ are the realizations.
 By Corollary~19.4,
 there is an $(s-1)$-almost fissile ensemble
 $q\in\<\W^\x(I)^{\VC\b\d E}\>$ such that
 \begin {equation} \label {q}
 \`\xg^E_I\'-q|_{(\b\d E)_+}
 \in
 \<\W(I)^{(\b\d E)_+}\>^{[s]}.
 \end {equation}
 Put $Q=\<\Upsilon_h\>(q)$
 (see the diagram).
 By Lemma~22.1,
 $Q\in\<(B(I)^{S^1}_\0)^{\vC\Delta E}\>$.
 By Lemma~22.3,
 $Q$ is $({S^1},s-1)$-almost fissile.
 Using the diagram,
 we get
 \begin {equation} \label {Q}
 \`\Upsilon_h(\xg^E_I)\'
 -
 Q|_{(\Delta E)_+}
 =
 \<\Upsilon_h\>(
 \`\xg^E_I\'
 -
 q|_{(\b\d E)_+}
 )
 \in
 \<(B(I)^{S^1})^{(\Delta E)_+}\>^{(s)}_{S^1},
 \end {equation}
 where $\in$ follows from \eqref{q} by Lemma~22.2.
 Let us drop the disjoint basepoint.
 Since
 \begin {multline*}
 \Upsilon_h(\xg^E_I)|_{\Delta E}
 =
 \hfill
 \text{(by construction of $\Upsilon_h$)}
 \hfill
 =
 \Xi^{\Delta E}(h(\top_{V(I)}))=
 \Xi^{\Delta E}(w),
 \end {multline*}
 we have
 $$
 \`\Xi^{\Delta E}(w)\'-Q|_{\Delta E}
 =
 (
 \`\Upsilon_h(\xg^E_I)\'-Q|_{(\Delta E)_+}
 )|_{\Delta E}
 \in
 \<(B(I)^{S^1})^{(\Delta E)}\>^{(s)}_{S^1},
 $$
 where $\in$ follows from \eqref{Q}
 (by
 \eqref{YXT}
 and
 Lemma~3.1).
 By Proposition~13.5,
 $\0\%\approx{s-1}w$.
 \qed
 \end {demo}


 \head {\S~24. Loops in an arbitrary space}


 \subhead {Nested commutators.}
 A {\it nesting\/} $t$ of {\it weight\/} $|t|\ge 1$ is
 either
 the atom $\atom$
 if $|t|=1$,
 or
 a pair $(t',t'')$ of nestings with $|t'|+|t''|=|t|$.
 Given
 elements $g_1,\dotsc,g_s$ of a group $G$
 and
 a nesting $t$ of weight $s$,
 the {\it $t$-nested commutator\/}
 $$
 \=t\(g_i\)_{i=1}^s\in G
 $$
 is defined
 (by induction on $s$)
 to be
 either
 $g_1$
 if $s=1$,
 or
 the commutator
 $$
 \(
 \,
 \={t'\!}\(g_i\)_{i=1}^{|t'|},
 \,
 \={t''\!}\(g_i\)_{i=|t'|+1}^{s}
 \,
 \)
 $$
 if $t=(t',t'')$.
 The nested commutators of weight $s$ in $G$ generate
 $\gamma^sG$,
 the $s$th term of the lower central series of $G$.


 \subhead {Loops.}
 Let $Y$ be a cellular space.
 We consider the group $\pi_1(Y)=[S^1,Y]$
 with the filtration $\pi_1(Y)^{((s))}=[S^1,Y]^{((s))}$
 (see \S~8).

 \begin {claim} [24.1. Theorem.]
 One has
 $$
 \pi_1(Y)^{((s))}
 =
 \gamma^s\pi_1(Y).
 $$
 \end {claim}

 Recall \cite[Theorem~11.2]{sim-1-cir}:
 \begin {equation} \label {pi1s}
 \pi_1(Y)^{(s)}
 =
 \gamma^s\pi_1(Y).
 \end {equation}
 Thus,
 by
 Theorem~8.2
 and
 \cite[Theorem~4.2]{sim-1-cir},
 the strong $r$-similarity on $\pi_1(Y)$
 coincides with the $r$-similarity.

 \begin {demo} [Proof.]
 The inclusion
 $\pi_1(Y)^{((s))}\subseteq\gamma^s\pi_1(Y)$ follows from
 the inclusion $\pi_1(Y)^{((s))}\subseteq\pi_1(Y)^{(s)}$
 (immediate from the definitions)
 and
 \eqref{pi1s}.

 Check that
 $\gamma^s\pi_1(Y)\subseteq\pi_1(Y)^{((s))}$.
 Since $\pi_1(Y)^{((s))}$ is a subgroup
 (by Theorem~8.1),
 it suffices to show that,
 for any
 nesting $t$ of weight $s$
 and
 loops $a_1,\dotsc,a_s:S^1\to Y$,
 one has
 $$
 \=t\([a_i]\)_{i=1}^s\in\pi_1(Y)^{((s))}.
 $$
 Put
 $$
 B_s
 =
 \bigbou_{i\in(s)}
 S^1
 \qquad
 \text{and}
 \qquad
 a
 =
 \bigBou_{i\in(s)}
 a_i:
 B_s
 \to
 Y.
 $$
 Let $\ins_i:S^1\to B_s$ be the canonical insertions.
 Choose a loop $e:S^1\to B_s$ with
 $$
 [e]
 =
 \=t\([\ins_i]\)_{i=1}^s
 $$
 in $\pi_1(B_s)$.
 So $\=t\([a_i]\)_{i=1}^s=[a\circ e]$.
 Clearly,
 the loop $e$ is Brunnian.
 By Lemma~23.1,
 $[e]\in\pi_1(B_s)^{((s))}$.
 By Corollary~5.2,
 $[a\circ e]\in\pi_1(Y)^{((s))}$,
 as was to be shown.
 \qed
 \end {demo}


 \head {\S~25. Whitehead products}


 Let $T_i$,
 $i=1,2$,
 be compact cellular spaces
 and
 $$
 T_i
 \xfrom{p_i}
 T_1\cro T_2
 \xto{k}
 T_1\sma T_2
 $$
 be the projections.
 The map
 $$
 \Sigma(T_1\cro T_2)
 \xto{\Sigma k}
 \Sigma(T_1\sma T_2)
 $$
 is homotopy right-invertible
 (because
 there is a canonical map $r$
 of the join $T_1*T_2$
 to $\Sigma(T_1\cro T_2)$
 such that
 $\Sigma k\circ r$ is a homotopy equivalence).
 Let $Y$ be a space.
 Given homotopy classes $\+a_i\in[\Sigma T_i,Y]$,
 $i=1,2$,
 consider the homotopy classes
 $$
 \+a_i\circ\Sigma p_i:
 \Sigma(T_1\cro T_2)
 \xto{\Sigma p_i}
 \Sigma T_i
 \xto{\+a_i}
 Y,
 \qquad
 i=1,2,
 $$
 and
 their commutator
 $$
 \(\+a_1\circ\Sigma p_1,\+a_2\circ\Sigma p_2\)
 \in
 [\Sigma(T_1\cro T_2),Y].
 $$
 The Whitehead product
 $$
 \[\+a_1,\+a_2\]
 \in
 [\Sigma(T_1\sma T_2),Y]
 $$
 is uniquely defined by (homotopy) commutativity of the diagram
 $$
 \xymatrix {
 \Sigma(T_1\cro T_2)
 \ar[rrr]^-{\(\+a_1\circ\Sigma p_1,\+a_2\circ\Sigma p_2\)}
 \ar[d]_-{\Sigma k}
 &&&
 Y
 \\
 \Sigma(T_1\sma T_2),
 \ar[urrr]_-{\[\+a_1,\+a_2\]}
 &&&
 }
 $$
 see \cite[Section~7.8]{S}.
 

 \subhead {Nested Whitehead products.}
 Let
 $T_i$,
 $i\in(s)$,
 be compact cellular spaces
 and
 $$
 T_i
 \xfrom{p_i}
 T_1\cro\dotso\cro T_s
 \xto{k}
 T_1\sma\dotso\sma T_s
 $$
 be the projections.

 \begin {claim} [25.1. Lemma.]
 The map
 $$
 \Sigma(T_1\cro\dotso\cro T_s)
 \xto{\Sigma k}
 \Sigma(T_1\sma\dotso\sma T_s)
 $$
 is homotopy right-invertible.
 \end {claim}

 \begin {demo} [Proof.]
 Induction on $s$.
 If $s=1$,
 $k$ is the identity.
 Take $s>1$.
 Put
 $$
 T'
 =
 T_1\cro\dotso\cro T_{s-1},
 \qquad
 Z'
 =
 T_1\sma\dotso\sma T_{s-1}.
 $$
 Let
 $$
 T'\cro T_s
 \xto{K}
 T'\sma T_s
 \qquad
 \text{and}
 \qquad
 T'
 \xto{k'}
 Z'
 $$
 be the projections.
 We have the decomposition
 $$
 \Sigma k:
 \Sigma(T'\cro T_s)
 \xto{\Sigma K}
 \Sigma(T'\sma T_s)
 \xto{\Sigma(k'\sma\id_{T_s})}
 \Sigma(Z'\sma T_s),
 $$
 where
 $\Sigma K$ is right-invertible
 (as noted in the beginning of \S~25)
 and
 the second arrow is right-invertible because
 it coincides with
 $$
 \Sigma T'\sma T_s
 \xto{\Sigma k'\sma\id_{T_s}}
 \Sigma Z'\sma T_s,
 $$
 which is right-invertible because
 $\Sigma k'$ is
 by the induction hypothesis.
 \qed
 \end {demo}

 Let
 $Y$
 be space
 and
 $\+a_i\in[\Sigma T_i,Y]$,
 $i\in(s)$,
 be homotopy classes.
 Given a nesting $t$ of weight $s$,
 define the $t$-nested Whitehead product
 $$
 \=t\[\+a_i\]_{i=1}^s
 \in
 \[\Sigma(T_1\sma\dotso\sma T_s),Y]
 $$
 (by induction on $s$)
 to be
 either
 $\+a_1$
 if $s=1$,
 or
 $$
 \[
 \,
 \={t'\!}\[\+a_i\]_{i=1}^{|t'|},
 \,
 \={t''\!}\[\+a_i\]_{i=|t'|+1}^s
 \,
 \]
 $$
 if $t=(t',t'')$.
 
 Consider the homotopy classes
 $$
 \+a_i\circ\Sigma p_i:
 \Sigma(T_1\cro\dotso\cro T_s)
 \xto{\Sigma p_i}
 \Sigma T_i
 \xto{\+a_i}
 Y,
 \qquad
 i\in(s).
 $$

 \begin {claim} [25.2. Lemma.]
 For a nesting $t$ of weght $s$,
 the diagram
 $$
 \xymatrix {
 \Sigma(T_1\cro\dotso\cro T_s)
 \ar[rrr]^-{
 \+c
 :=
 \=t\(\+a_i\circ\Sigma p_i\)_{i=1}^s
 }
 \ar[d]_-{\Sigma k}
 &&&
 Y
 \\
 \Sigma(T_1\sma\dotso\sma T_s)
 \ar[urrr]_-{
 \+w
 :=
 \=t\[\+a_i\]_{i=1}^s
 }
 &&&
 }
 $$
 is (homotopy) commutative.
 \end {claim}

 \begin {demo} [Proof.]
 Induction on $s$.
 If $s=1$,
 $\Sigma k$ is the identity
 and
 $\+c=\+w=\+a_1$.
 Take $s>1$.
 We have $t=(t',t'')$.
 Put
 $s'=|t'|$,
 $s''=|t''|$,
 and
 \begin {align*}
 T'
 &=
 T_1\cro\dotso\cro T_{s'},
 &
 T''
 &=
 T_{s'+1}\cro\dotso\cro T_s,
 \\
 Z'
 &=
 T_1\sma\dotso\sma T_{s'},
 &
 Z''
 &=
 T_{s'+1}\sma\dotso\sma T_s.
 \end {align*}
 We have the commutative diagrams of projections
 $$
 \xymatrix {
 &
 T'\cro T''
 \ar[d]^-{P'}
 \ar[dl]_-{p_i}
 \\
 T_i
 &
 T'
 \ar[l]_-{p'_i}
 \ar[r]^-{k'}
 &
 Z',
 }
 \qquad
 \xymatrix {
 &
 T'\cro T''
 \ar[d]^-{P''}
 \ar[dl]_-{p_i}
 \\
 T_i
 &
 T''
 \ar[l]_-{p''_i}
 \ar[r]^-{k''}
 &
 Z''
 }
 $$
 ($i\le s'$ on the left,
 $i\ge s'+1$ on the right).
 Consider the diagram
 \begin {equation} \label {w'}
 \xymatrix {
 \Sigma(T'\cro T'')
 \ar[d]_-{\Sigma P'}
 \ar[drrr]^*!/r2ex/{\labelstyle
 \~{\+c}'
 :=
 \={t'\!}\(\+a_i\circ\Sigma p_i\)_{i=1}^{s'}
 }
 &&&
 \\
 \Sigma T'
 \ar[rrr]^(0.4){
 \={t'\!}\(\+a_i\circ\Sigma p'_i\)_{i=1}^{s'}
 }
 \ar[d]_-{\Sigma k'}
 &&&
 Y
 \\
 \Sigma Z'.
 \ar[urrr]_-{
 \+w'
 :=
 \={t'\!}\[\+a_i\]_{i=1}^{s'}
 }
 &&&
 }
 \end {equation}
 The upper triangle is commutative because
 the function
 $$
 [\Sigma T',Y]
 \to
 [\Sigma(T'\cro T''),Y]
 $$
 induced by $\Sigma P'$
 is a homomorphism
 and sends $\+a_i\circ\Sigma p'_i$ to $\+a_i\circ\Sigma p_i$.
 The lower triangle is commutative
 by the induction hypothesis.
 Similarly,
 we have the commutative diagram
 \begin {equation} \label {w''}
 \xymatrix {
 \Sigma(T'\cro T'')
 \ar[d]_-{\Sigma P''}
 \ar[drrr]^*!/r2ex/{\labelstyle
 \~{\+c}''
 :=
 \={t''\!}\(\+a_i\circ\Sigma p_i\)_{i=s'+1}^s
 }
 &&&
 \\
 \Sigma T''
 \ar[rrr]^(0.4){
 \={t''\!}\(\+a_i\circ\Sigma p''_i\)_{i=s'+1}^s
 }
 \ar[d]_-{\Sigma k''}
 &&&
 Y
 \\
 \Sigma Z''.
 \ar[urrr]_*!/r2ex/{\labelstyle
 \+w''
 :=
 \={t''\!}\[\+a_i\]_{i=s'+1}^s
 }
 &&&
 }
 \end {equation}
 We have the commutative diagram of projections
 $$
 \xymatrix {
 T'
 \ar[d]_{k'}
 &&
 T'\cro T''
 \ar[ll]_-{P'}
 \ar[rr]^-{P''}
 \ar[d]^-{k'\cro k''}
 \ar@/_8ex/[dd]_(0.7){k}
 &&
 T''
 \ar[d]^-{k''}
 \\
 Z'
 &&
 Z'\cro Z''
 \ar[ll]_-{Q'}
 \ar[rr]^-{Q''}
 \ar[d]^-{K}
 &&
 Z''
 \\
 &&
 Z'\sma Z''.
 &&
 }
 $$
 Consider the diagram
 $$
 \xymatrix {
 \Sigma(T'\cro T'')
 \ar[d]_-{\Sigma(k'\cro k'')}
 \ar[drrrrr]^-{
 \+c
 =
 \(\~{\+c}',\~{\+c}''\)
 }
 \ar@/_16ex/[dd]_{\Sigma k}
 &&&&&
 \\
 \Sigma(Z'\cro Z'')
 \ar[rrrrr]^(0.4){
 \(
 \+w'\circ\Sigma Q',
 \+w''\circ\Sigma Q''
 \)
 }
 \ar[d]_-{\Sigma K}
 &&&&&
 Y
 \\
 \Sigma(Z'\sma Z'').
 \ar[urrrrr]_-{
 \+w=\[\+w',\+w''\]
 }
 &&&&&
 }
 $$
 The upper triangle is commutative because
 the function
 $$
 [\Sigma(Z'\cro Z''),Y]
 \to
 [\Sigma(T'\cro T''),Y]
 $$
 induced by $\Sigma(k'\cro k'')$
 is a homomorphism
 under which
 \begin {multline*}
 \+w'\circ\Sigma Q'
 \mapsto
 \+w'\circ\Sigma k'\circ\Sigma P'
 =
 \hfill
 \text{(by diagram \eqref{w'})}
 \hfill
 =
 \~{\+c}'
 \end {multline*}
 and
 \begin {multline*}
 \+w''\circ\Sigma Q''
 \mapsto
 \+w''\circ\Sigma k''\circ\Sigma P''
 =
 \hfill
 \text{(by diagram \eqref{w''})}
 \hfill
 =
 \~{\+c}''.
 \end {multline*}
 The lower triangle is commutative
 by the definition of Whitehead product.
 We are done.
 \qed
 \end {demo}

 \begin {claim} [25.3. Corollary.]
 Let $R$ be a homotopy right-inverse of $\Sigma k$:
 $$
 \xymatrix {
 \Sigma(T_1\cro\dotso\cro T_s)
 \ar[rr]^-{\Sigma k}
 &&
 \Sigma(T_1\sma\dotso\sma T_s)
 \ar@/^4ex/[ll]^-{R},
 }
 \qquad
 \Sigma k \circ R\sim\id.
 $$
 Then,
 for any nesting $t$ of weght $s$,
 the diagram
 $$
 \xymatrix {
 \Sigma(T_1\cro\dotso\cro T_s)
 \ar[rrr]^-{
 \+c
 :=
 \=t\(\+a_i\circ\Sigma p_i\)_{i=1}^s
 }
 &&&
 Y
 \\
 \Sigma(T_1\sma\dotso\sma T_s).
 \ar[u]^-{R}
 \ar[urrr]_-{
 \+w
 :=
 \=t\[\+a_i\]_{i=1}^s
 }
 &&&
 }
 $$
 is (homotopy) commutative.
 \end {claim}

 \begin {demo} [Proof.]
 We have
 \begin {multline*}
 \+c\circ R
 =
 \hfill
 \text{(by Lemma~25.2)}
 \hfill
 =
 \+w\circ\Sigma k\circ R
 =
 \hfill
 \text{(since $\Sigma k \circ R\sim\id$)}
 \hfill
 =
 \+w.
 \QED
 \end {multline*}
 \end {demo}


 \head {\S~26. Loops and Whitehead products}


 Consider the wedge
 $$
 B_s
 =
 \bigbou_{i\in(s)}
 S^1.
 $$
 Given
 a loop $v:S^1\to B_s$
 and
 a space $T$,
 introduce the map $v^\Sigma$:
 $$
 \xymatrix @R2.8ex {
 \Sigma T
 \ar[rr]^-{v^\Sigma}
 \ar@{=}[d]
 &&
 \bigbou\limits_{i\in(s)}
 \Sigma T
 \ar@{=}[d]
 \\
 S^1\sma T
 \ar[rr]^-{v\sma\id_T}
 &&
 B_s\sma T.
 }
 $$

 Let
 $$
 \ins_j:S^1
 \to
 B_s
 \quad
 \text{and}
 \quad
 \ins^T_j:\Sigma T
 \to
 \bigbou_{i\in(s)}
 \Sigma T,
 \qquad
 j\in(s),
 $$
 be the canonical insertions.

 \begin {claim} [26.1. Lemma.]
 The function
 $$
 \pi_1(B_s)
 \to
 [
 \Sigma T,
 \bigbou_{i\in(s)}
 \Sigma T
 ],
 \qquad
 [v]\mapsto[v^\Sigma],
 $$
 is a homomorhism,
 under which $[\ins_i]\mapsto[\ins^T_i]$.
 \qed
 \end {claim}

 Let $T_i$,
 $i\in(s)$,
 be spaces
 and
 $$
 T_i
 \xfrom{p_i}
 T_1\cro\dotso\cro T_s
 \xto{k}
 T_1\sma\dotso\sma T_s
 $$
 be the projections.
 Let
 $Y$
 be a space
 and
 $a_i:\Sigma T_i\to Y$
 be maps.
 We have the compositions
 $$
 a_i\circ\Sigma p_i:
 \Sigma(T_1\cro\dots\cro T_s)
 \xto{\Sigma p_i}
 \Sigma T_i
 \xto{a_i}
 Y.
 $$

 \begin {claim} [26.2. Lemma.]
 Let $t$ be a nesting of weight $s$.
 Let $e:S^1\to B_s$ be a loop with
 $$
 [e]
 =
 \=t\([\ins_i]\)_{i=1}^s
 $$
 in $\pi_1(B_s)$.
 Then
 the diagram
 $$
 \xymatrix {
 \Sigma(T_1\cro\dots\cro T_s)
 \ar[rr]^-{e^\Sigma}
 \ar[drr]_*!/l2ex/{\labelstyle
 \+c
 :=
 \=t\([a_i]\circ\Sigma p_i\)_{i=1}^s
 }
 &&
 \bigbou\limits_{i\in(s)}
 \Sigma(T_1\cro\dots\cro T_s)
 \ar[d]^-{
 A
 :=
 \bigBou\limits_{i\in(s)}
 a_i\circ\Sigma p_i
 }
 \\
 &&
 Y.
 }
 $$
 is (homotopy) commutative.
 \end {claim}

 \begin {demo} [Proof.]
 Put $T=T_1\cro\dotso T_s$.
 By Lemma~26.1,
 the function
 $$
 \pi_1(B_s)
 \to
 [
 T,
 \bigbou_{i\in(s)}
 \Sigma T
 ],
 \qquad
 [v]\mapsto[v^\Sigma],
 $$
 is a homomorhism,
 under which $[\ins_i]\mapsto[\ins^T_i]$.
 Thus
 $$
 [e^\Sigma]
 =
 \=t\([\ins^T_i]\)_{i=1}^s.
 $$
 The map $A$ induces a homomorphism
 $$
 [
 \Sigma T,
 \bigbou_{i\in(s)}
 \Sigma T
 ]
 \to
 [\Sigma T,Y],
 $$
 under which $[\ins^T_i]\mapsto[a_i]\circ\Sigma p_i$
 and
 thus
 $$
 [e^\Sigma]
 =
 \=t\([\ins^T_i]\)_{i=1}^s
 \mapsto
 \=t\([a_i]\circ\Sigma p_i\)_{i=1}^s
 =
 \+c,
 $$
 which is what was to be shown.
 \qed
 \end {demo}

 By Lemma~25.1,
 we have the diagram
 $$
 \xymatrix {
 \Sigma(T_1\cro\dotso\cro T_s)
 \ar[rr]^-{\Sigma k}
 &&
 \Sigma(T_1\sma\dotso\sma T_s)
 \ar@/^4ex/[ll]^-{R},
 }
 $$
 where $\Sigma k\circ R\sim\id$.
 For a loop $v:S^1\to B_s$,
 introduce the composition $\prescript vR\[a_i\]_{i=1}^s$:
 $$
 \xymatrix {
 \Sigma(T_1\cro\dots\cro T_s)
 \ar[rr]^-{v^\Sigma}
 &&
 \bigbou\limits_{i\in(s)}
 \Sigma(T_1\cro\dots\cro T_s)
 \ar[d]^-{
 A
 :=
 \bigBou\limits_{i\in(s)}
 a_i\circ\Sigma p_i
 }
 \\
 \Sigma(T_1\sma\dots\sma T_s)
 \ar[u]^-{R}
 \ar[rr]^-{\prescript vR\[a_i\]_{i=1}^s}
 &&
 Y.
 }
 $$

 \begin {claim} [26.3. Lemma.]
 Let
 $t$ be a nesting of weight $s$
 and
 $e:S^1\to B_s$ be a loop with
 $$
 [e]
 =
 \=t\([\ins_i]\)_{i=1}^s
 $$
 in $\pi_1(B_s)$.
 Then
 $$
 [\,\prescript eR\[a_i\]_{i=1}^s\,]
 =
 \=t\[[a_i]\]_{i=1}^s
 $$
 in $[\Sigma(T_1\sma\dotso\sma T_s),Y]$.
 \end {claim}

 \begin {demo} [Proof.]
 Recall the homotopy class
 $$
 \Sigma(T_1\cro\dots\cro T_s)
 \xto{
 \+c
 :=
 \=t\([a_i]\circ\Sigma p_i\)_{i=1}^s
 }
 Y.
 $$
 We have
 \begin {multline*}
 [\prescript eR\[a_i\]_{i=1}^s]
 =
 [A\circ e^\Sigma\circ R]
 =
 \hfill
 \text{(by Lemma~26.2)}
 \qquad
 \\
 \qquad
 =
 \+c\circ R
 =
 \hfill
 \text{(by Corollary~25.3)}
 \hfill
 =
 \=t\[[a_i]\]_{i=1}^s.
 \QED
 \end {multline*}
 \end {demo}


 \head {\S~27. Strong nullarity of Whitehead products}


 Let $a_i:\Sigma T_i\to Y$,
 etc.,
 be as in \S~26.

 \begin {claim} [27.1. Lemma.]
 Let $v:S^1\to B_s$ be a loop such that
 $\0\%\approx r v$.
 Then
 $$
 \0
 \%\approx r
 \prescript vR\[a_i\]_{i=1}^s.
 $$
 \end {claim}

 \begin {demo} [Proof.]
 We have
 $$
 \prescript vR\[a_i\]_{i=1}^s
 =
 A
 \circ
 v^\Sigma
 \circ
 R
 $$
 (see the construction).
 By Corollary~5.4,
 $\0\%\approx r v^\Sigma$.
 By Corollary~5.2,
 $\0\%\approx rA\circ v^\Sigma\circ R$.
 \qed
 \end {demo}

 Given
 a nesting $t$ of weight $s$
 and
 homotopy classes $\+a_i\in[\Sigma T_i,Y]$,
 $i\in(s)$,
 consider the $t$-nested Whitehead product
 $$
 \=t\[\+a_i\]_{i=1}^s
 \in
 [
 \Sigma(T_1\sma\dotso T_s),
 Y
 ].
 $$

 \begin {claim} [27.2. Theorem.]
 One has 
 $$
 \=t\[\+a_i\]_{i=1}^s
 \in
 [\Sigma(T_1\sma\dotso T_s),Y]^{((s))}.
 $$
 \end {claim}

 \begin {demo} [Proof.]
 For each $i$,
 choose a representative $a_i:\Sigma T_i\to Y$ of $\+a_i$.
 Choose a loop $e:S^1\to B_s$ with
 $$
 [e]
 =
 \=t\([\ins_i]\)_{i=1}^s
 $$
 in $\pi_1(B_s)$.
 Clearly,
 the loop $e$ is Brunnian.
 By Lemma~23.1,
 $\0\%\approx{s-1}e$.
 By Lemma~27.1,
 $$
 \0
 \%\approx{s-1}
 \prescript eR\[a_i\]_{i=1}^s.
 $$
 By Lemma~26.3,
 $$
 [\,\prescript eR\[a_i\]_{i=1}^s\,]
 =
 \=t\[\+a_i\]_{i=1}^s.
 $$
 Thus
 $$
 \=t\[\+a_i\]_{i=1}^s
 \in
 [\Sigma(T_1\sma\dotso T_s),Y]^{((s))}.
 \QED
 $$
 \end {demo}



 \noindent
 \href{mailto:ssp@pdmi.ras.ru}{\tt ssp@pdmi.ras.ru}

 \noindent
 \url{http://www.pdmi.ras.ru/~ssp}

 \end {document}